\documentclass[a4paper,12pt,final]{article}

\usepackage[english]{babel}
\usepackage[utf8]{inputenc}
\usepackage[T1]{fontenc}
\usepackage{lmodern}
\usepackage{graphicx}
\usepackage{setspace}
\usepackage{amsmath,amssymb}
\usepackage{hyperref}
\usepackage[french]{varioref}
\usepackage{color}
\usepackage[final]{pdfpages}
\usepackage{mdframed}
\usepackage{lipsum}
\usepackage[a4paper]{geometry}
\newcommand{\reporttitle}{A scalable Lagrange-Remap scheme for compressible multimaterial Euler equations with sharp interface reconstruction}     
\newcommand{\reportauthor}{Bastien \textsc{Chaudet}} 
\newcommand{\reportsubject}{Projet de Fin d'\'Etudes} 
\newcommand{\HRule}{\rule{\linewidth}{0.5mm}}
\hypersetup{
    pdftitle={\reporttitle}, %
    pdfauthor={\reportauthor},%
    pdfsubject={\reportsubject},%
}

\begin{document}

  \begin{titlepage}

\begin{center}

\begin{minipage}[t]{0.32\textwidth}
  \begin{flushleft}
    \includegraphics [width=33mm]{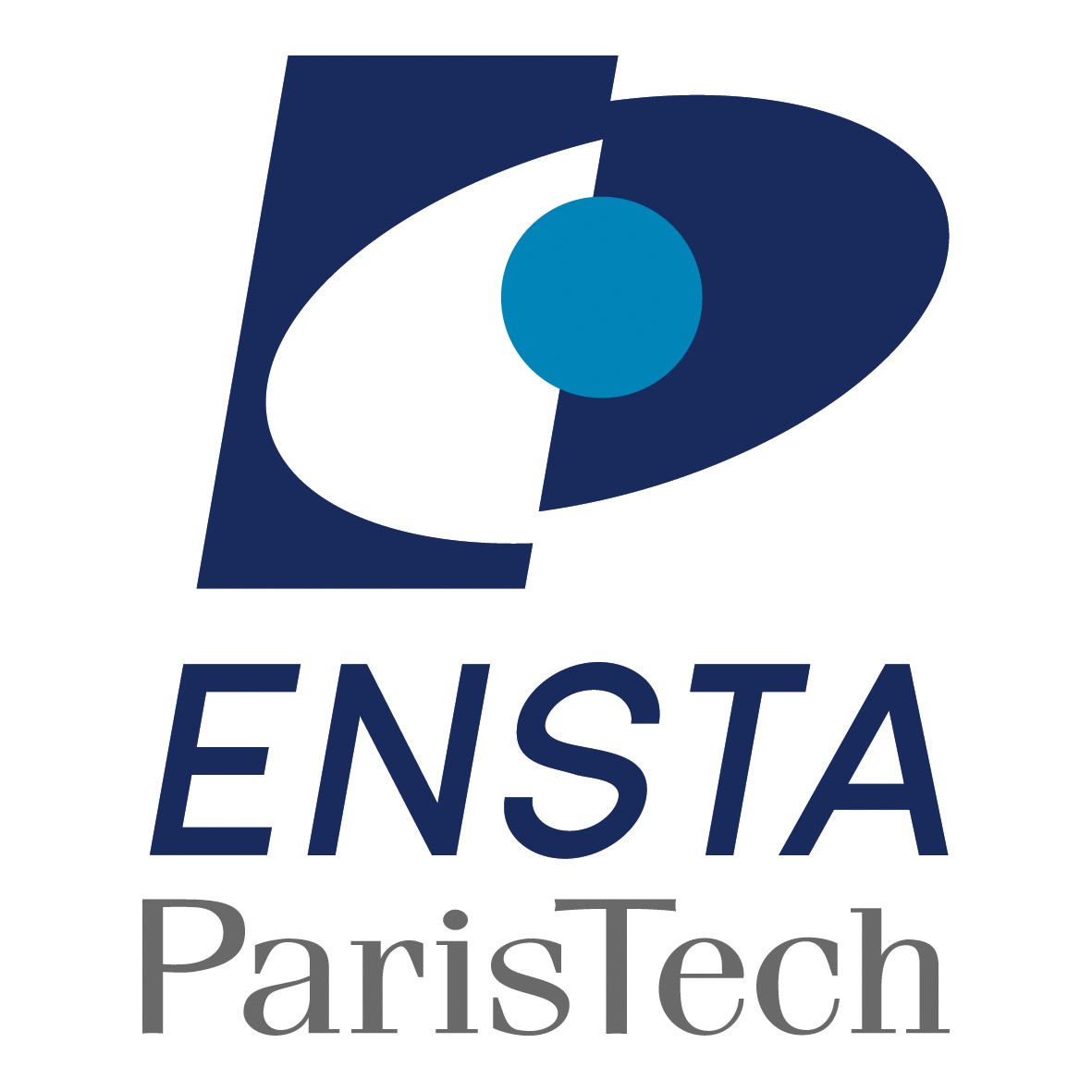} \\[2cm]
  \end{flushleft}
\end{minipage}
\begin{minipage}[t]{0.32\textwidth}
  \begin{center}
    \includegraphics [width=30mm]{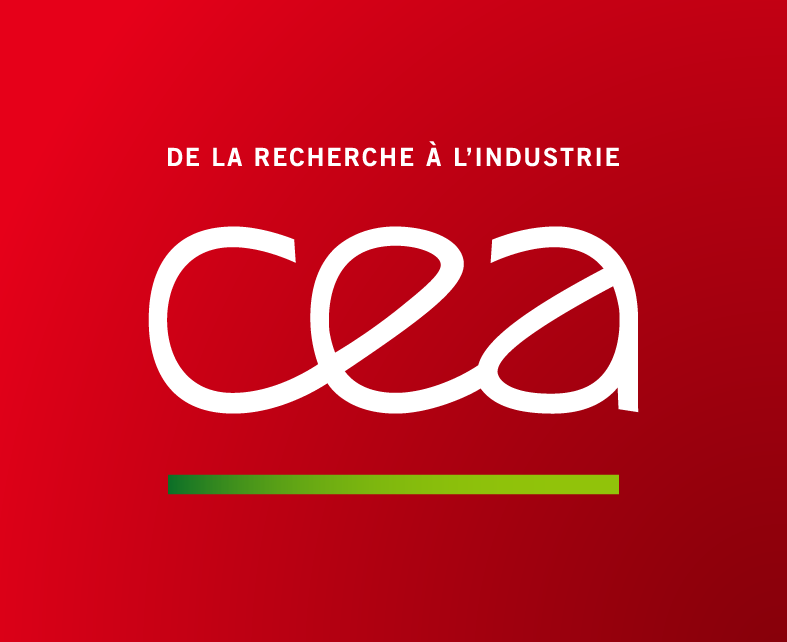} \\[2cm]
  \end{center}
\end{minipage}
\begin{minipage}[t]{0.32\textwidth}
  \begin{flushright}
    \includegraphics [width=30mm]{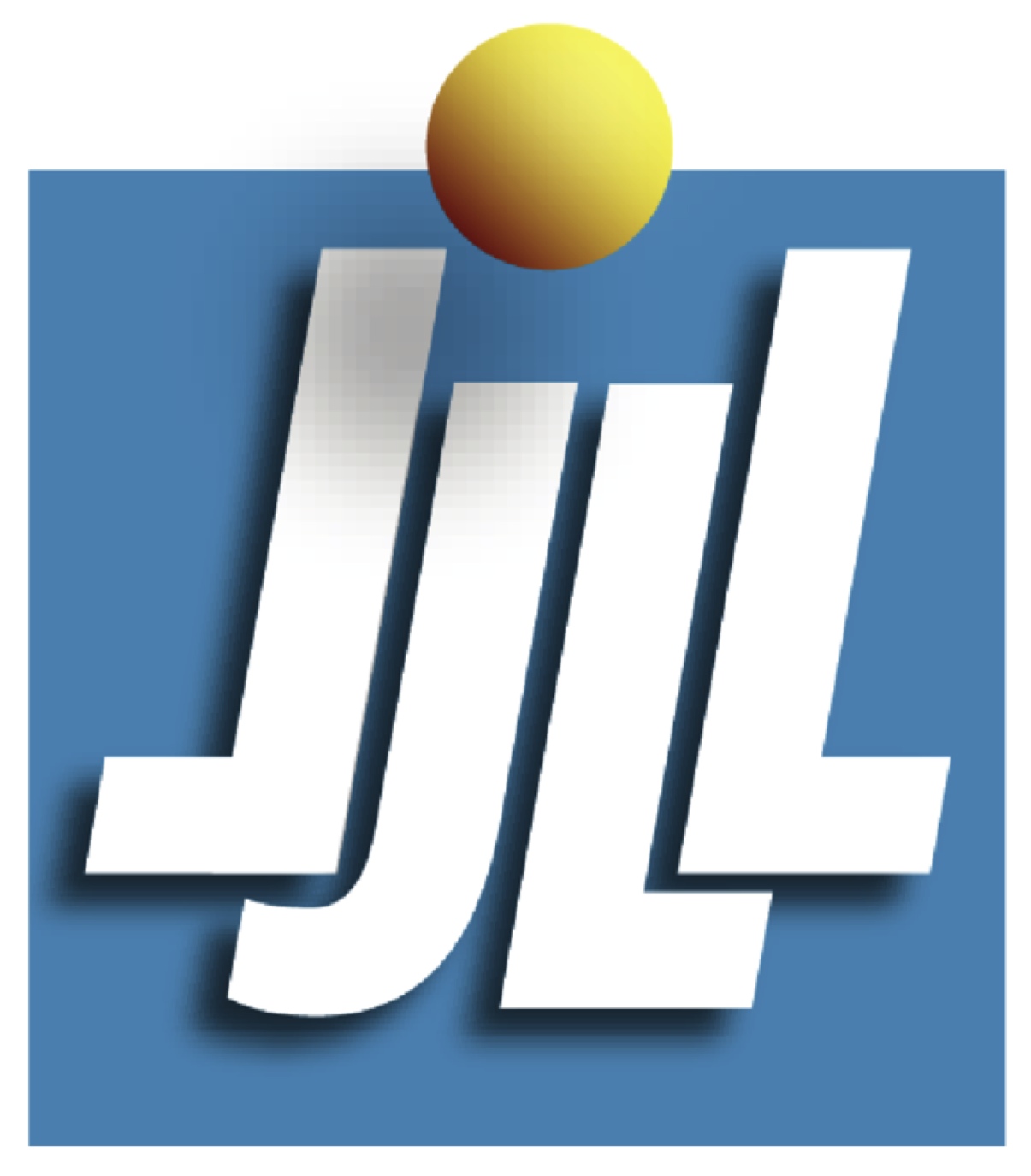} \\[2cm]
  \end{flushright}
\end{minipage} \\[2.2cm]

\textsc{\Large \reportsubject}\\[0.5cm]
\HRule \\[0.6cm]
{
{\huge \bfseries \reporttitle}\\[0.6cm]
}
\HRule \\[1.5cm]

\begin{minipage}[t]{0.3\textwidth}
  \begin{flushleft} 
    \emph{Author:}\\
    \reportauthor
  \end{flushleft}
\end{minipage}
\begin{minipage}[t]{0.6\textwidth}
  \begin{flushright} 
    \emph{Supervisors:} \\
    Dr.~Jean-Philippe \textsc{Braeunig} (CEA)\\
    Pr.~Christophe \textsc{Hazard} (ENSTA)
  \end{flushright}
\end{minipage}

\vfill

\begin{flushleft}
{\large \textcolor{red}{\textsc{This report is NOT CONFIDENTIAL.}}}
\end{flushleft}

\begin{flushleft}
{\large Internship at the CEA, DAM, DIF, F-91297 Arpajon Cedex.}
\end{flushleft}

{\large From June to November 2014}

\end{center}

\end{titlepage}

  \clearpage
  \newpage
  \thispagestyle{empty} 
  \strut
  \newpage
  \newpage
  \thispagestyle{empty}
\section*{Certificate of non confidentiality}

\hspace{1.5em}The present report is not confidential. Accordingly, any disclosure of its content in paper and/or electronic format to third parties is allowed.
  \clearpage
  \thispagestyle{empty}
\section*{Acknowledgments}

\hspace{1.5em}First of all, I would like to thank St\'ephane Bernard and Jean-Philippe Perlat for greeting me in their service and laboratory, respectively.

I also thank Patrick Le Tallec and Florian De Vuyst for their pieces of advise, and for giving time to fruitful discussions.

I am particularly grateful to my supervisor Jean-Philippe Braeunig for his unchanging help and involvement, for his availability, and his support. His ideas and the clarity of his explanations have greatly contributed to the work presented in this report.

Thanks to Mathieu Peybernes and Rapha\"el Poncet for providing such a clear numerical code, and for answering all my questions.

  \clearpage 
  \setcounter{page}{1}
\section*{Abstract}
\addcontentsline{toc}{section}{Abstract}

\hspace{1.5em}This work is in the field of multi-material compressible fluid flows simulation. The proposed scheme is eulerian and related to finite volumes methods, but in a Lagrange-Remap formalism on regular orthogonal meshes. The Lagrangian scheme is staggered and the remap phase is similar to a finite volume advection scheme. The multi-material extension uses classical VOF fluxes for sharp interface reconstruction. The originality of the scheme is in the attempt for a 9 points remap scheme without directional splitting. This strategy should allow to preserve good properties of classical multi-material staggered schemes, while saving parallel communications with the one step remap. Results will be discussed and compared to those from classical Lagrange-Remap schemes on severe benchmarks. 

\textbf{Keywords:} Lagrange-Remap, Euler equations, multi-fluid flows, VOF, directional splitting, HPC, hydrodynamics, numerical analysis
  \newpage
  \thispagestyle{empty} 
  \strut
  \newpage
  \thispagestyle{empty} 
  \strut
  \tableofcontents 
  \clearpage
  \listoffigures 
  \clearpage
  \newpage
  \strut
  \section*{Introduction}
\addcontentsline{toc}{section}{Introduction}

\hspace{1.5em}Fluid mechanics, especially hydrodynamics, has first been conceptualized by L. Euler in 1755 when he introduced the equation of perfect fluids $\rho \frac{\partial u}{\partial t} + \rho (\mathbf{u}\cdot\nabla)\mathbf{u} = - \nabla \!P$. Then it has been generalized by C. Navier and G. Stokes by taking into account the viscosity of the fluid in the equation. In the framework of this internship, we got interested in compressible Euler equations, neglecting the viscous term. L. Euler's approach was to study these equations into a fixed referential, unlike J. L. Lagrange, who took the point of view of a referential following the matter during its displacement. 

Those points of view can be extended to the numerical resolution of these equations. Indeed, Lagrangian simulation codes are implemented such that the mesh moves with the fluid. They are often preferred to Eulerian simulation codes because they naturally refine the mesh around zones of interest (high compression zones, shocks, etc) and give more faithful retranscriptions of contact discontinuities. However, they are not really robust to vorticity, since the mesh tends to deform with the fluid. An alternative can be found in ALE (\textbf{A}rbitrary  \textbf{L}agrangian \textbf{E}ulerian) schemes. The difference between ALE and pure Lagrangian schemes comes from the remap phase. Indeed, in order to fix the problems of stability and robustness of Lagrangian schemes when vortices appear in the simulation, Lagrangian mesh and variables are remapped onto a regularized mesh at each time step. This avoids to get tangled cells, while saving the interesting geometric properties of Lagrangian schemes. The specific case of remapping onto the initial Eulerian mesh at each time step is the Lagrange-Remap formalism. This formalism constitutes the basis of the work presented in this report. 

The extension to multi-material models enables to simulate much more complex flows. There are several physical and numerical models for multi-fluid flows, more or less diffusive, depending on the miscibility of the fluids in the physical case, and on the properties of the scheme. A numerical model adapted to the simulation of miscible fluids is the mixing model: mixed cells, i.e.\ cells containing two species, are represented with one single pressure, and a mass concentration per material. Another numerical multi-material model is the reconstruction of a sharp interface between the two fluids, which of course models non-miscible fluids, and multiphase flows. The (discontinuous) interface is represented as a segment in each mixed-cell, and is placed with volume fractions and the normal vector to this segment. \\

In the framework of HPC (\textbf{H}igh \textbf{P}erformance \textbf{C}omputing), a trend in numerical analysis is to build numerical schemes adapted to exascale computers constraints, i.e.\ using a lot of computational cores. Two of the main requirements for a scalable scheme are to include kernels with high arithmetic intensities, and to reduce as much as possible the number of communications needed when parallelizing with MPI (\textbf{M}essage \textbf{P}assing \textbf{I}nterface). The performances of a code using many cores can be restricted significantly if such features are not taken into account. Still regarding HPC issues,  it is often important to reduce the complexity of algorithms. Here, we consider fixed Eulerian regular orthogonal grids.\\
There exists one type of remap used in most codes: the AD (\textbf{A}lternate \textbf{D}irections) remap. It is known to be quite accurate and robust but it shows one main drawback: it does not seem very adapted to HPC. To be more specific, in the framework of MPI parallelization, since the remap is processed in two steps, the number of communications is high. An alternative remap - the Direct remap -, using the same fluxes but remapping in one step, has been proposed by Q. Debray and J.-P. Braeunig in 2013, see \cite{debray2013projection}. \\
The work presented in this report follows on from the latter, and aims at testing a new type of remap which is more adapted to HPC, and incorporating it in a 2D hydrodynamics Lagrange-Remap scheme. The main idea is still to perform the remap in one step, using the same tools as the AD remap, while recovering at least the same robustness and accuracy as the AD remap. \\
The work done during this internship includes the theoretical study of this new remap (Direct with Corner Fluxes), getting familiar with SHY code, a mono-material Lagrange-Remap research code with an AD remap provided by the CEA (M. Peybernes and R. Poncet), the implementation of Direct and Direct with Corner Fluxes remaps in SHY, and the extension to two-materials for all remaps: theoretical study and implementation of two different multi-fluid models. An additional theoretical study on the diffusion of some hydrodynamics eulerian schemes has been conducted during the first two months, that took place at the LRC-MESO (\textbf{L}aboratoire de \textbf{R}echerche \textbf{C}onventionn\'e Mod\'elisation \textbf{MESO}scopique), CMLA (\textbf{C}entre de \textbf{M}a-th\'ematiques et \textbf{L}eurs \textbf{A}pplications), ENS Cachan. A report of this study can be found in the Annexe.\\

First, the Lagrangian phase will be presented, writing both continuous and discrete systems of equations. Then, the principles of three different types of remap: the AD remap, the Direct remap and the remap studied in this work (the Direct Remap with Corner Fluxes), will be explained. The Direct Remap with Corner Fluxes will be analysed in further details, before coming to the multi-material extension of the three schemes. Finally, several results of convergence and robustness will conclude this report. 

\clearpage

 \newpage
 \strut
 \newpage

\section{Lagrangian phase}

	\subsection{Lagrangian formalism}

\hspace{1.5em}The Lagrangian formalism consists in solving the equations in the matter referential. At this aim, a special derivative has been introduced: the material derivative $\mathrm{d}_t = \partial _t+ \mathbf{u} \cdot \nabla_x$ , where $\mathbf{u}$ denotes the velocity of the fluid particle. Let the classical thermodynamical quantities be: $\rho$ (density), $e$ (specific internal energy) and $P$ (pressure). A term of artificial viscosity $Q$, called Wilkins pseudo-viscosity, see \cite{wilkins1963calculation}, is added in order to remain entropic. It  makes the scheme more stable by spreading shocks and regularizing the solution at its neighbourhoods. This continuity is used to recover local differential equations from integral ones. \\
With these notations, in this formalism, the equations to be solved on the moving domain $\Omega(t)$ write:
	
\begin{equation}
\label{continuous_system}
	\left\lbrace \hspace{0.05\textwidth}
	\begin{array}{l}
 		\dfrac{\mathrm{d}}{\mathrm{d}t} \left( \displaystyle \int_{\Omega(t)} \rho \: \mathrm{d}x \right) = 0\\
 		\:\\
 		\dfrac{\mathrm{d}}{\mathrm{d}t} \:\mathbf{x} = \mathbf{u}\\ 		
 		\: \\
 		\rho \dfrac{\mathrm{d}}{\mathrm{d}t} \:\mathbf{u} = - \nabla (P + Q)\\
 		\: \\
 		\dfrac{\mathrm{d}}{\mathrm{d}t} \: e = - (P + Q) \: \dfrac{\mathrm{d}}{\mathrm{d}t} \left( \dfrac{1}{\rho} \right)\\
 		\: \\
 		\: P = \mathcal{P}(\rho,e)
 	\end{array}
 	\right.
\end{equation}
\vspace{0.3cm}

\noindent The Wilkins pseudo-viscosity is the following:
\begin{displaymath}
Q = a_1 \rho c \: \lvert \Delta \mathbf{u} \rvert + a_2 \rho \:{\lvert \Delta \mathbf{u} \rvert}^2  
\end{displaymath}
where $\Delta \mathbf{u}$ is the velocity gap in the direction of the shock, $a_1$, $a_2$ are two real numbers called respectively the linear and quadratic pseudo-viscosity coefficients, and $c$ is the sound speed.   

Let us now proceed to both space and time discretizations of system \eqref{continuous_system}.

	\subsection{Discretization}
	
\hspace{1.5em}A SGPC (\textbf{S}taggered  \textbf{G}rid \textbf{P}redictor \textbf{C}orrector) scheme is used. Space discretization is based on Wilkins', which is second order in space, see \cite{vonneumann1950method}. It is a predictor corrector scheme in time, which is second order as well.
	
\begin{figure}[h]
\begin{center}
\includegraphics[width=0.6\textwidth, height=0.45\textwidth]{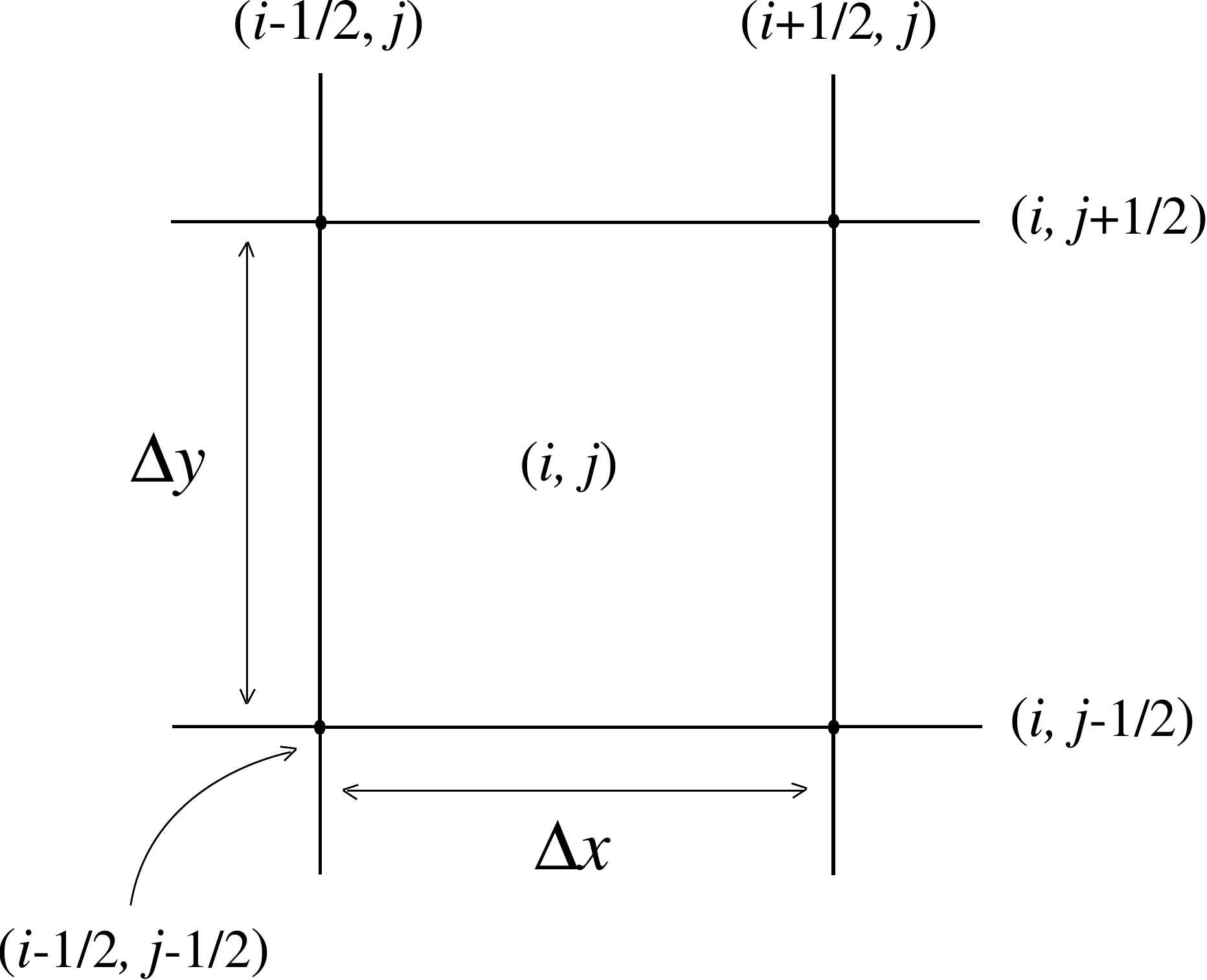}
\vspace{0.3cm}
\caption{Notations of indices in cell $c = (i,j)$.}
\label{fig:notation_indices}
\end{center}
\end{figure}

	We denote by $\Delta x$ and $\Delta y$ the space steps in the two directions, and $\Delta t^n$ the time step at $t^n$ (iteration $n$). The prediction phase computes values at $t^{n+1/2} = t^n + \Delta t^n / 2$. Cell-centred variables and node-centred variables are written respectively in cell $c = (i,j)$ in $(\Delta x , \Delta y)$ basis, and at node $p = (i \pm 1/2,j \pm 1/2)$, see Figure~\ref{fig:notation_indices}.
	
\vspace{0.5cm}

The domain $\Omega(t^n)$ is splitted up into rectangular cells like the one represented in Figure~\ref{fig:notation_indices}, i.e.\ $\Omega(t^n) = \bigcup \Omega_c^n$. The whole mesh is deformed with respect to the movement of the fluid, that is, the matter contained into $\Omega_c^n$ at $t^n$ is to be found in $\Omega_c^{n+1, \:lag}$ at the end of the Lagrangian phase. In particular, the mass in each moving subdomain $\Omega_c$ remains constant. The discrete system of equation writes, $\forall \:p$, $c$:

	\paragraph{Prediction phase}
	
\begin{equation}
\label{prediction}
	\left\lbrace \hspace{0.05\textwidth}
	\begin{array}{l}
 		m_c^{n+1/2} = m_c^n\\
 		\: \\
 		\mathbf{x}_p^{n+1/2} = \mathbf{x}_p^n + \dfrac{\Delta t^n}{2} \mathbf{u}_p^n\\
 		\: \\ 		
 		Vol_c^{n+1/2} = \mathcal{F}\big( (\mathbf{x}_p^{n+1/2})_{p \in \{c\}}\big)\\
 		\: \\
 		\mathbf{u}_p^{n+1/2} = \mathbf{u}_p^{n} - \dfrac{\Delta t^n}{2} \big( \nabla (P + Q) \big)_p^{n+1/2} / \rho_p^n\\
 		\: \\
 		e_c^{n+1/2} = e_c^n - (P_c^n + Q_c^n) \left( \dfrac{1}{\rho_c^{n+1/2}} - \dfrac{1}{\rho_c^{n}} \right)\\
 		\: \\
 		P_c^{n+1/2} = \mathcal{P}(\rho_c^{n+1/2},e_c^{n+1/2})
 	\end{array}
 	\right.
\end{equation}	
\vspace{0.3cm}
 
\begin{displaymath}
\mbox{where} \:\:\:\:\:\: \mathcal{F}\big( (\mathbf{x}_p)_{p \in \{c\}}\big) = \frac{1}{2} \: \left\| \Delta \mathbf{x}_{i-1/2,j} \wedge \Delta \mathbf{x}_{i,j-1/2} \right\|
\end{displaymath}
\begin{displaymath}
 \:\:\:\:\:\:\:\:\:\:\:\:\:\:\:\:\:\: + \frac{1}{2} \:\left\| \Delta \mathbf{x}_{i+1/2,j} \wedge \Delta \mathbf{x}_{i,j+1/2} \right\|
\end{displaymath}

\begin{displaymath}
\mbox{with} \:\:\:\:\:\: \Delta \mathbf{x}_{i\pm1/2,j} = \mathbf{x}_{i\pm1/2,j+1/2} - \mathbf{x}_{i\pm1/2,j-1/2}
\end{displaymath}
\begin{displaymath}
\: \mbox{and} \:\:\:\:\:\:\: \Delta \mathbf{x}_{i,j\pm1/2} = \mathbf{x}_{i+1/2,j\pm1/2} - \mathbf{x}_{i-1/2,j\pm1/2}
\end{displaymath}$\:$\\
Moreover, we have, at node $p = (i-1/2,j-1/2)$:
\begin{displaymath}
\begin{array}{l}
\big( \nabla (P + Q) \big)_p^{n+1/2}\cdot\mathbf{e}_x = \dfrac{1}{2 \: \Delta x} \left[\: \big(\:P_{i,j-1}^{n+1/2} + Q_{i,j-1}^n + P_{i,j}^{n+1/2} + Q_{i,j}^n\:\big) \right. \\
\: \\
\:\:\:\:\:\:\:\:\:\:\:\:\:\:\:\:\:\:\:\:\:\:\:\:\:\:\:\:\:\:\: \left. - \:\big(\:P_{i-1,j-1}^{n+1/2} + Q_{i-1,j-1}^n + P_{i-1,j}^{n+1/2} + Q_{i-1,j}^n\:\big)\:\right]\\
\: \\
\big( \nabla (P + Q) \big)_p^{n+1/2}\cdot\mathbf{e}_y = \dfrac{1}{2 \: \Delta y} \left[\: \big(\:P_{i-1,j}^{n+1/2} + Q_{i-1,j}^n + P_{i,j}^{n+1/2} + Q_{i,j}^n\:\big) \right. \\
\: \\
\:\:\:\:\:\:\:\:\:\:\:\:\:\:\:\:\:\:\:\:\:\:\:\:\:\:\:\:\:\:\: \left. - \:\big(\:P_{i-1,j-1}^{n+1/2} + Q_{i-1,j-1}^n + P_{i,j-1}^{n+1/2} + Q_{i,j-1}^n\:\big)\:\right]\\
\end{array}
\end{displaymath}

\noindent where nodal quantities (mass $m_p^n$ and density $\rho_p^n$) are defined as follows
\begin{displaymath}
\rho_p^n = \dfrac{m_p^n}{\Delta x \Delta y} \:, \:\:\:\:\:\: \mbox{with} \:\:\:\:\: m_p^n = \dfrac{1}{4} \displaystyle \sum_{c \in \{p\}} m_c^n \:.
\end{displaymath}

 In order to compute the discrete pseudo-viscosity, the value of $\Delta \mathbf{u}$ is approximated by $L \:\mathrm{div}\:\mathbf{u}$, $L=\sqrt{\Delta x \Delta y}$ being the characteristic length of a cell, which leads to:
\begin{displaymath}
Q_c^n = \left\lbrace 
	\begin{array}{l r}
	 \rho_c^n a_1 L \:c_c^n \:\lvert (\mathrm{div}\:\mathbf{u})_c^n \rvert + a_2 L^2 \:{\lvert (\mathrm{div}\:\mathbf{u})_c^n \rvert}^2 & \:\:\:\:\:\: \mbox{if} \:\: (\mathrm{div}\:\mathbf{u})_c^n < 0, \\
	 0 & \mbox{else.}
	 \end{array}
	 \right. 
\end{displaymath}
 
\begin{figure}[h]
\begin{center}
\includegraphics[width=0.44\textwidth, height=0.36\textwidth]{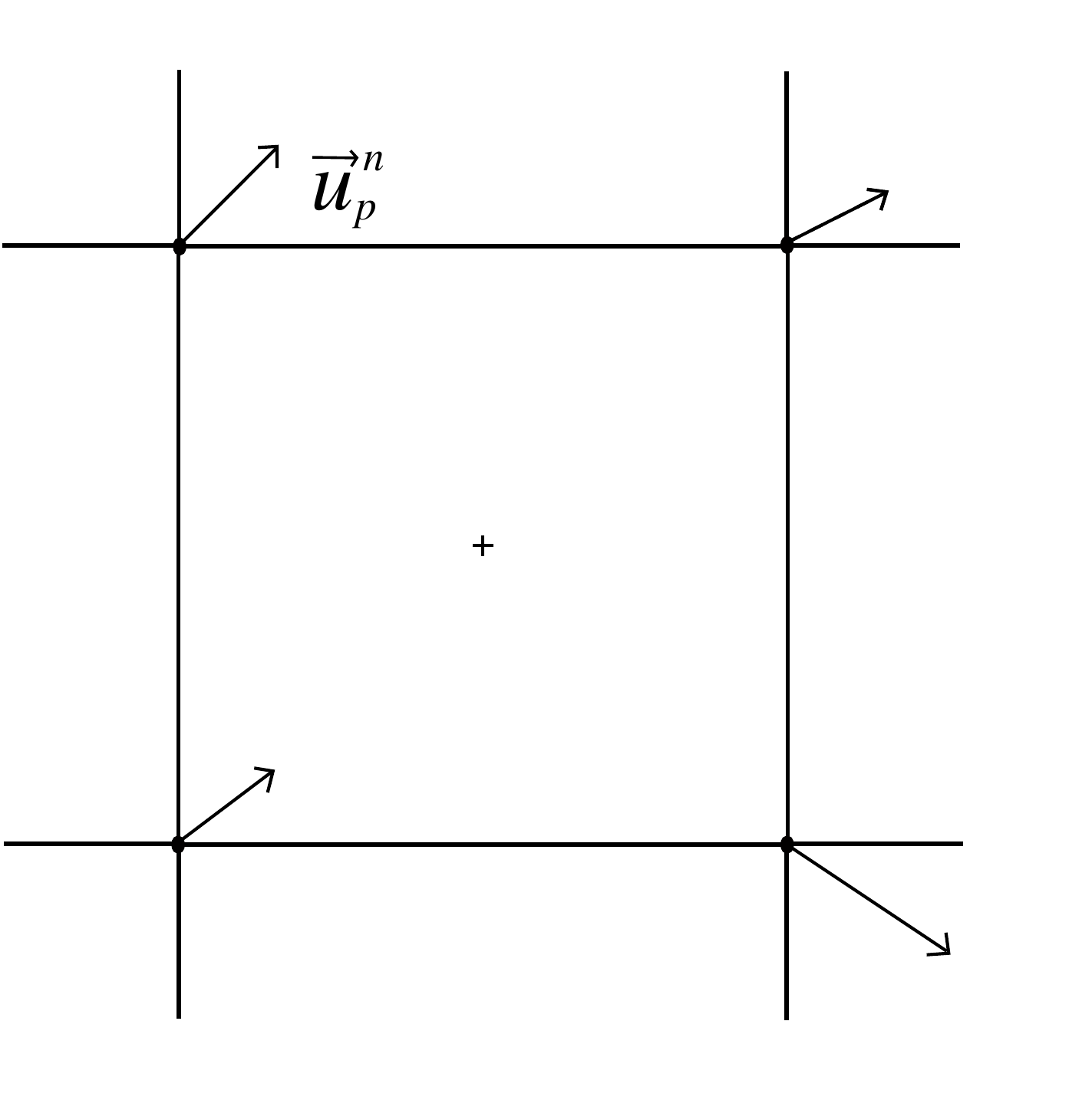}\hspace{0.05\textwidth}
\includegraphics[width=0.44\textwidth, height=0.36\textwidth]{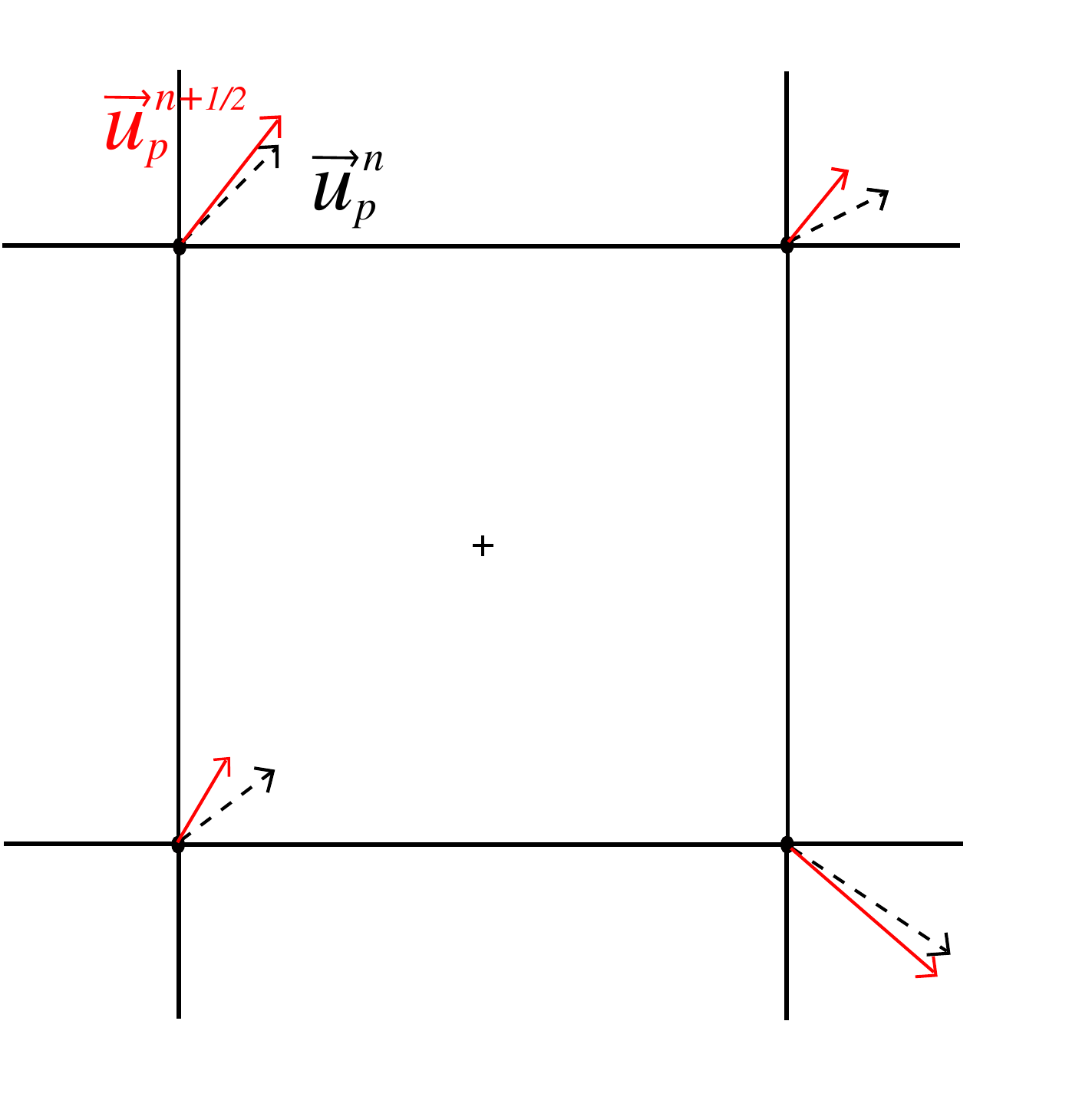}
\caption{Beginning (left) and end (right) of prediction phase.}
\label{fig:prediction_phase}
\end{center}
\end{figure}

\begin{figure}[h]
\begin{center}
\includegraphics[width=0.44\textwidth, height=0.36\textwidth]{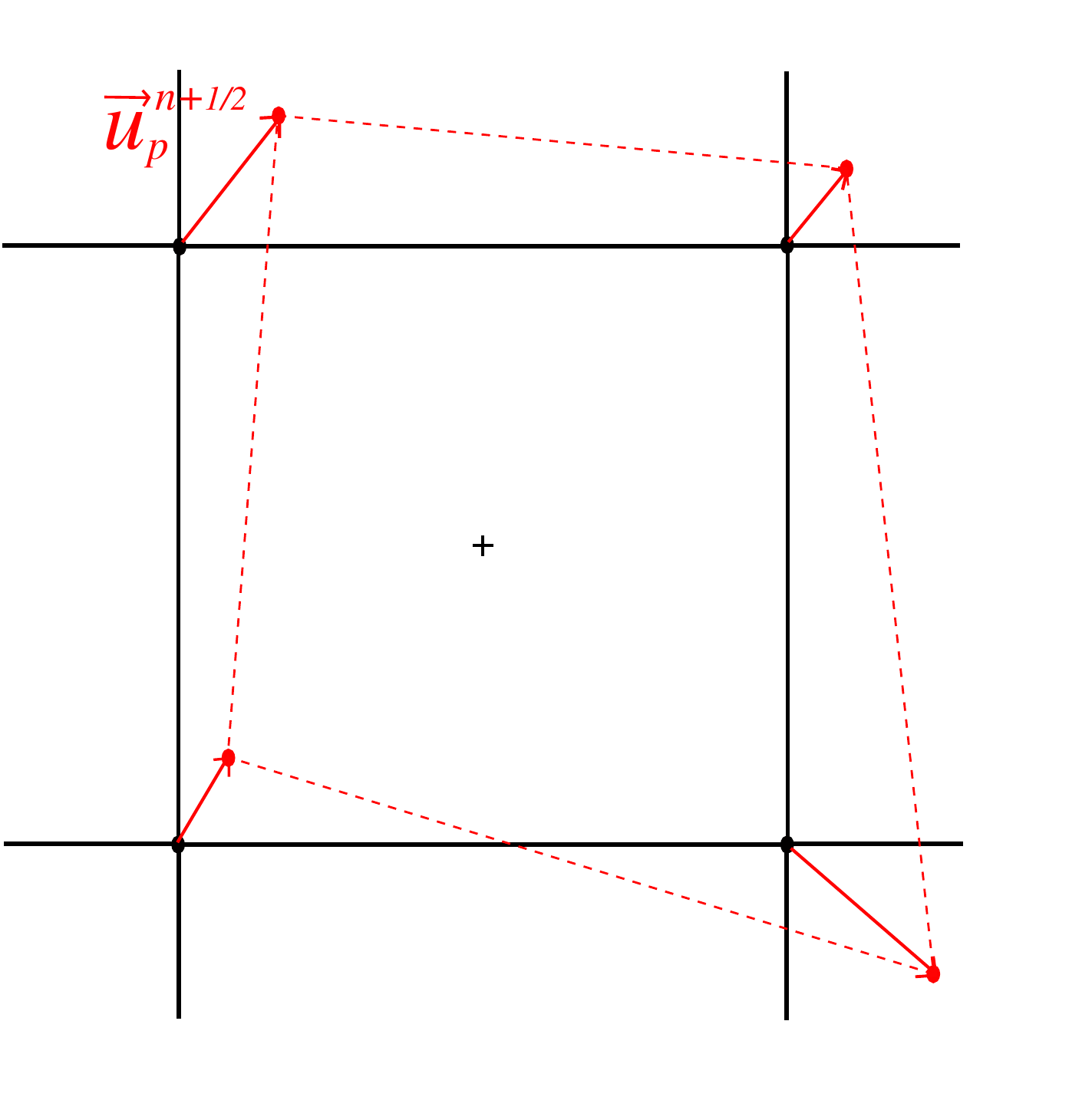}\hspace{0.05\textwidth}
\includegraphics[width=0.44\textwidth, height=0.36\textwidth]{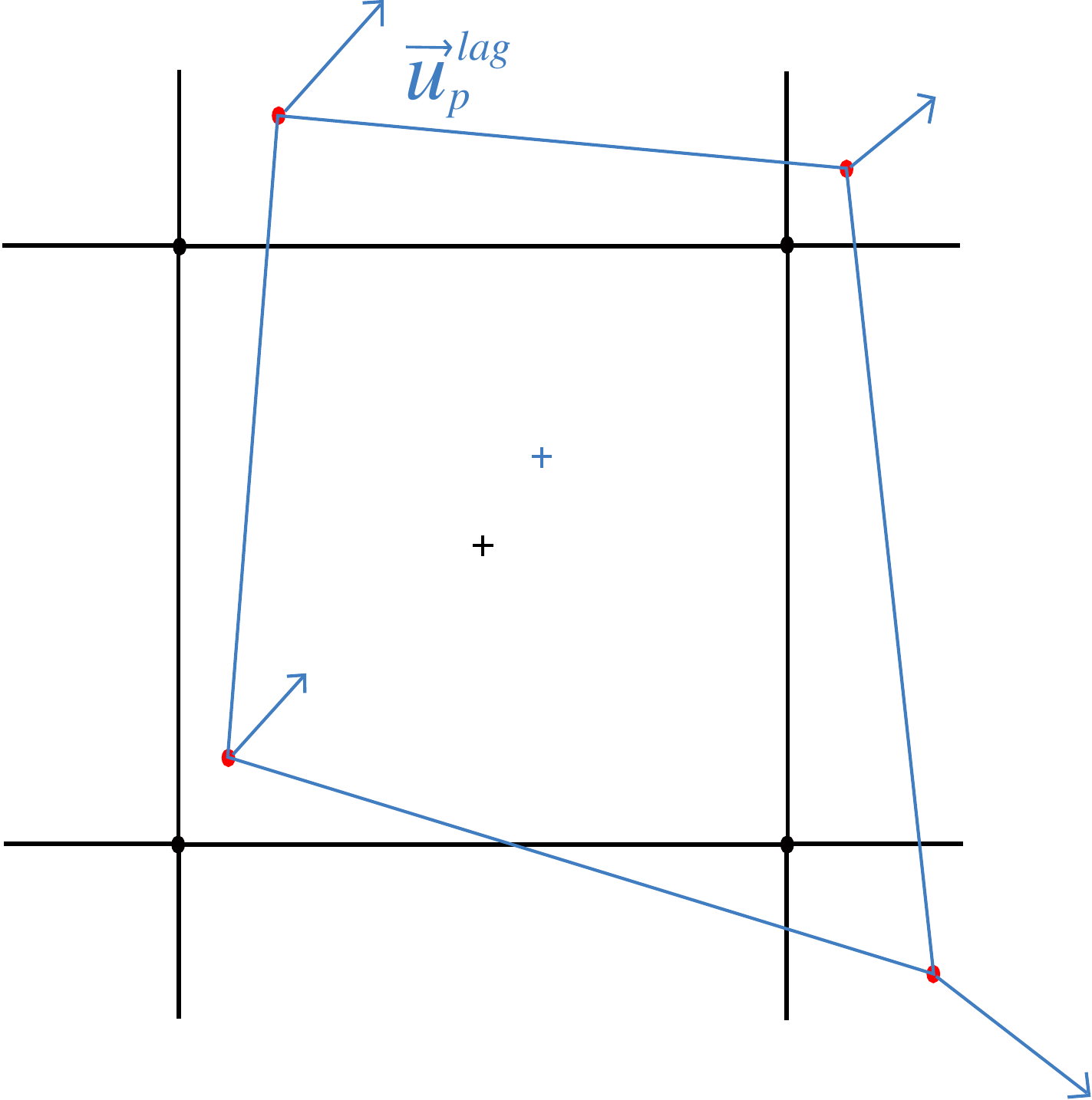}
\caption{Beginning (left) and end (right) of correction phase.}
\label{fig:correction_phase}
\end{center}
\end{figure}

  The values of $\mathbf{u}_p^{n+1/2}$ enable to perform the Lagrangian displacement of the mesh, see Figure~\ref{fig:prediction_phase}. Then the correction phase uses the values at $t^{n+1/2}$ to express Lagrangian variables, i.e.\ variables at the centres and nodes of the Lagrangian deformed mesh, see Figure~\ref{fig:correction_phase}. The end of the Lagrangian phase is denoted $t^{n+1, \:lag} = t^n + \Delta t^n$.
  
	\paragraph{Correction phase}
	
\begin{equation}
\label{correction}
	\left\lbrace \hspace{0.05\textwidth}
	\begin{array}{l}
 		m_c^{n+1, \:lag} = m_c^n\\
 		\: \\
 		\mathbf{x}_p^{n+1, \:lag} = \mathbf{x}_p^n + \Delta t^n \mathbf{u}_p^{n+1/2}\\
 		\: \\ 		
 		Vol_c^{n+1, \:lag} = \mathcal{F}\big( (\mathbf{x}_p^{n+1})_{p \in \{c\}}\big)\\
 		\: \\
 		\mathbf{u}_p^{n+1, \:lag} = 2 \: \mathbf{u}_p^{n+1/2} - \mathbf{u}_p^{n}\\
 		\: \\
 		e_c^{n+1, \:lag} = e_c^n - (P_c^{n+1/2} + Q_c^n) \left( \dfrac{1}{\rho_c^{n+1, \:lag}} - \dfrac{1}{\rho_c^{n}} \right)\\
 		\: \\
 		P_c^{n+1, \:lag} = \mathcal{P}(\rho_c^{n+1, \:lag},e_c^{n+1, \:lag})
 	\end{array}
 	\right.
\end{equation}

\vspace{0.1cm}

The two successive phases enable to get the second order in time. Plus, the scheme is entropic (under a CFL type condition) thanks to the pseudo-viscosity. It is also conservative in mass (locally and globally). But regarding total energy or momentum, it is not so clear since one cannot define rigorously those quantities because of the node-centred velocities. However, there exist a way to define both quantities such that the scheme is conservative in momentum, and in total energy as well, but with an error in $O({\Delta t}^2)$. For further details about the computations, see \cite{courtesinglard2012schemas}.

 \clearpage
 \newpage
 \section{Remap phase}

\begin{figure}[h]
\begin{center}
\includegraphics[width=0.8\textwidth, height=0.35\textwidth]{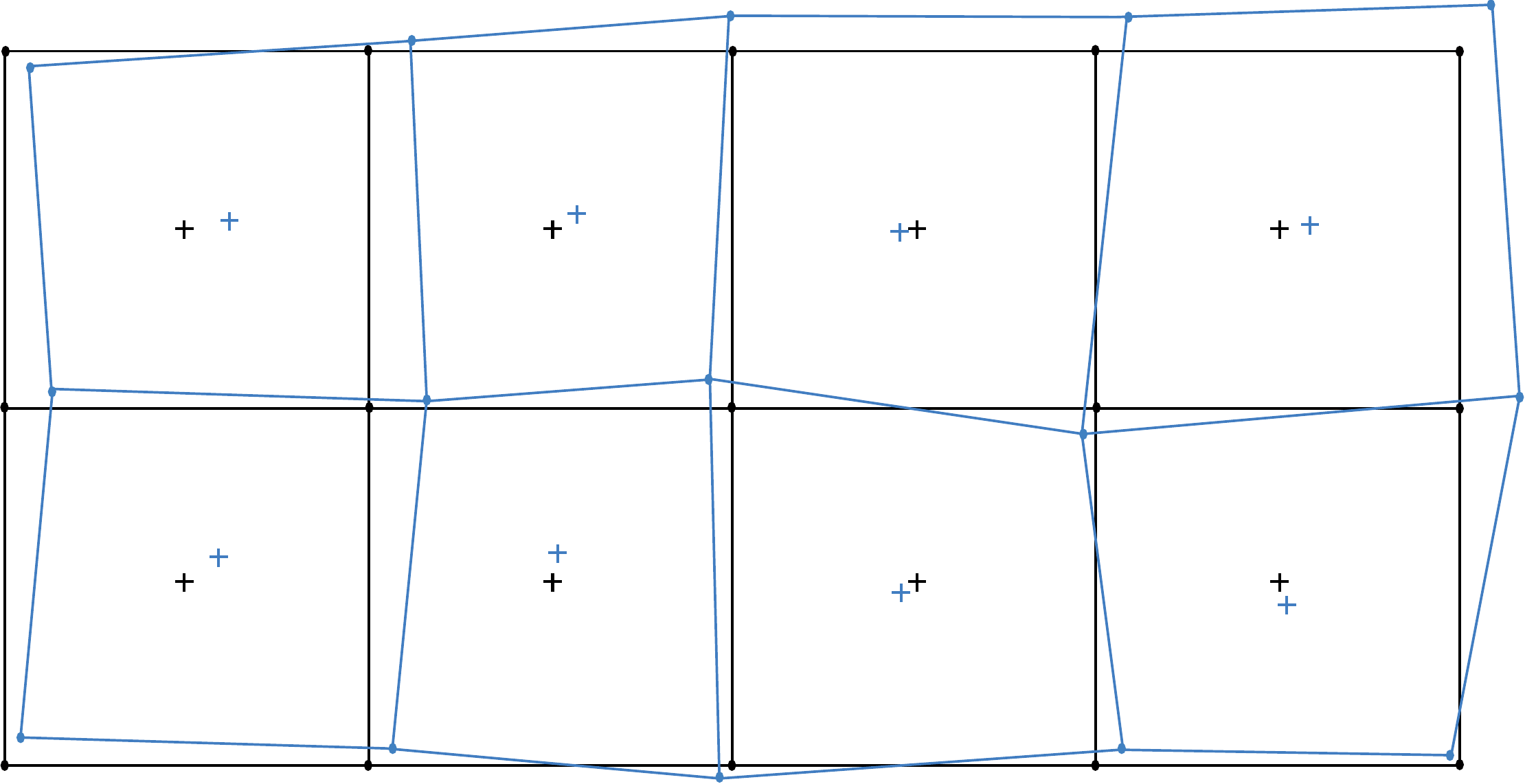}
\caption{Lagrangian and Eulerian meshes.}
\label{fig:interpolation_phase}
\end{center}
\end{figure}

\hspace{1.5em}The principle is the following: the mesh moves during the Lagrangian phase and the values of Lagrangian variables are computed, then new values of $\rho$, $e$ and $\mathbf{u}$ on the Eulerian mesh are interpolated from values on the Lagrangian mesh. In Figure~\ref{fig:interpolation_phase}, it means that one has to interpolate values at black crosses and dots from values at blue crosses and dots. This phase is only a geometric interpolation, all equations ruling physical properties have been solved numerically during the Lagrangian phase. 

\vspace{0.5cm}

\noindent Variables are remapped (interpolated) thanks to a Finite-Volume type remap scheme. That is to say, volume fluxes at the faces of the Eulerian mesh are computed from the Lagrangian displacement of the nodes, which enables to get variable fluxes and finally values after the remap phase in a conservative formalism. \\
Somehow, this second phase can be seen as a pure advection of Lagrangian variables, which explains the Finite-Volume type scheme. The remap phase is often called advection phase. 

\begin{figure}[h]
\begin{center}
\includegraphics[width=0.45\textwidth, height=0.35\textwidth]{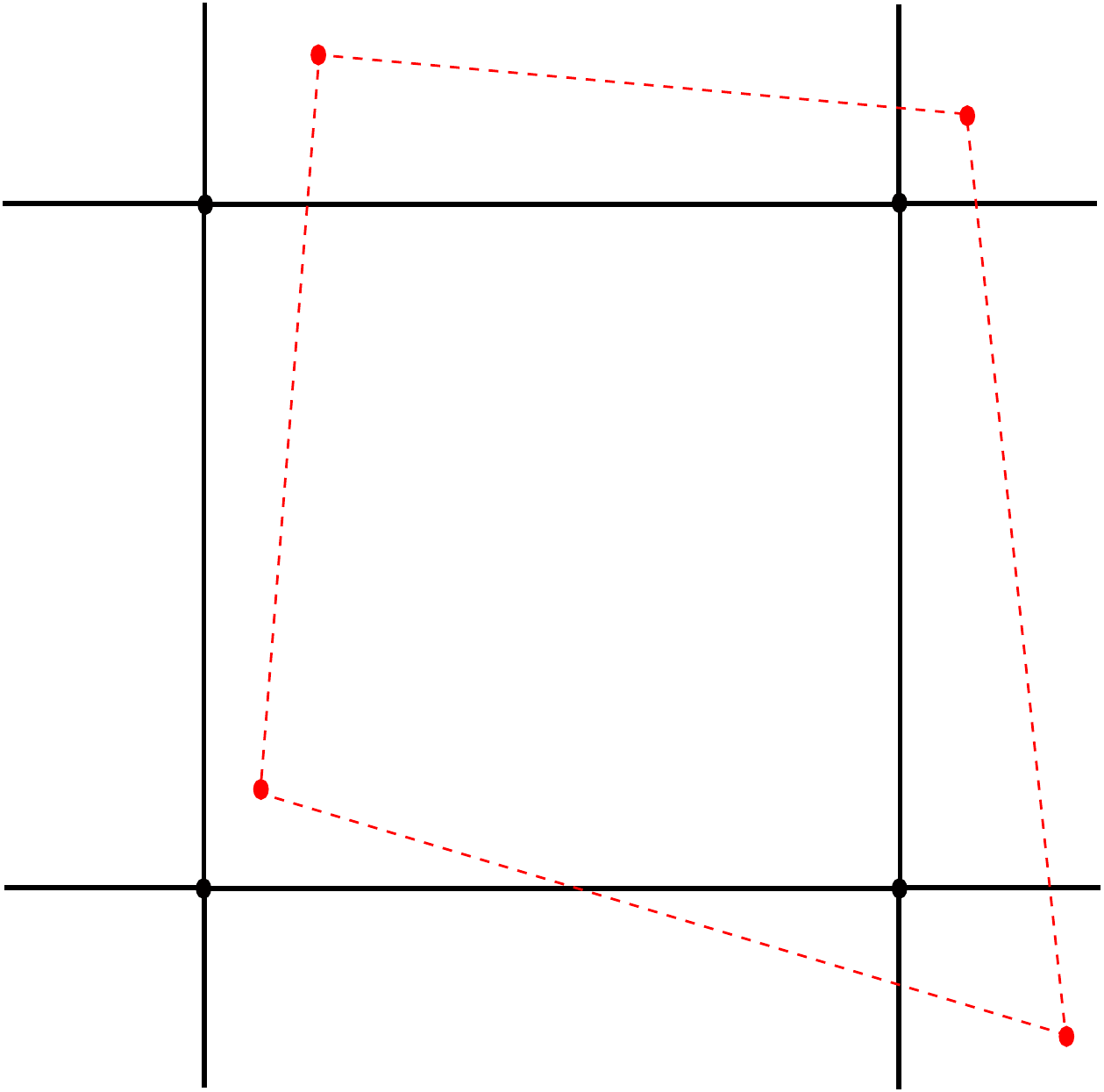}
\caption{Lagrangian displacement of cell $c = (i,j)$.}
\label{fig:lagrangian_displacement}
\end{center}
\end{figure}

In the following sections, several types of remap will be presented. The Lagrangian displacement of one reference cell $c = (i,j)$, knowing the velocities at each node, is represented Figure~\ref{fig:lagrangian_displacement}. At this step, the Lagrangian phase is over, i.e.\ the values of variables at $t^{n+1, \:lag}$ have been computed from those at $t^n$. From now, the exponent $^{n+1, \:lag}$ will be replaced by $^{lag}$ in order to simplify notations.

\vspace{0.5cm}

\noindent Black lines define the fixed Eulerian mesh, whereas red dotted lines define the faces of the Lagrangian cell.

\subsection{Alternate Directions Remap (AD)}

\hspace{1.5em}As mentioned before, two successive steps (remaps) are required to process the full AD remap. First, only $x$ components of the velocities are taken into account to compute fluxes at all $X$-faces (i.e.\ vertical faces), and remap all variables. Second, only $y$ components of the velocities are taken into account to compute fluxes at all $Y$-faces (i.e.\ horizontal faces), and remap all variables once again. The name "Alternate Directions" comes from this $X$-$Y$ alternating.

\subsubsection{$X$-Remap phase}

\begin{figure}[h]
\begin{center}
\includegraphics[width=0.5\textwidth, height=0.4\textwidth]{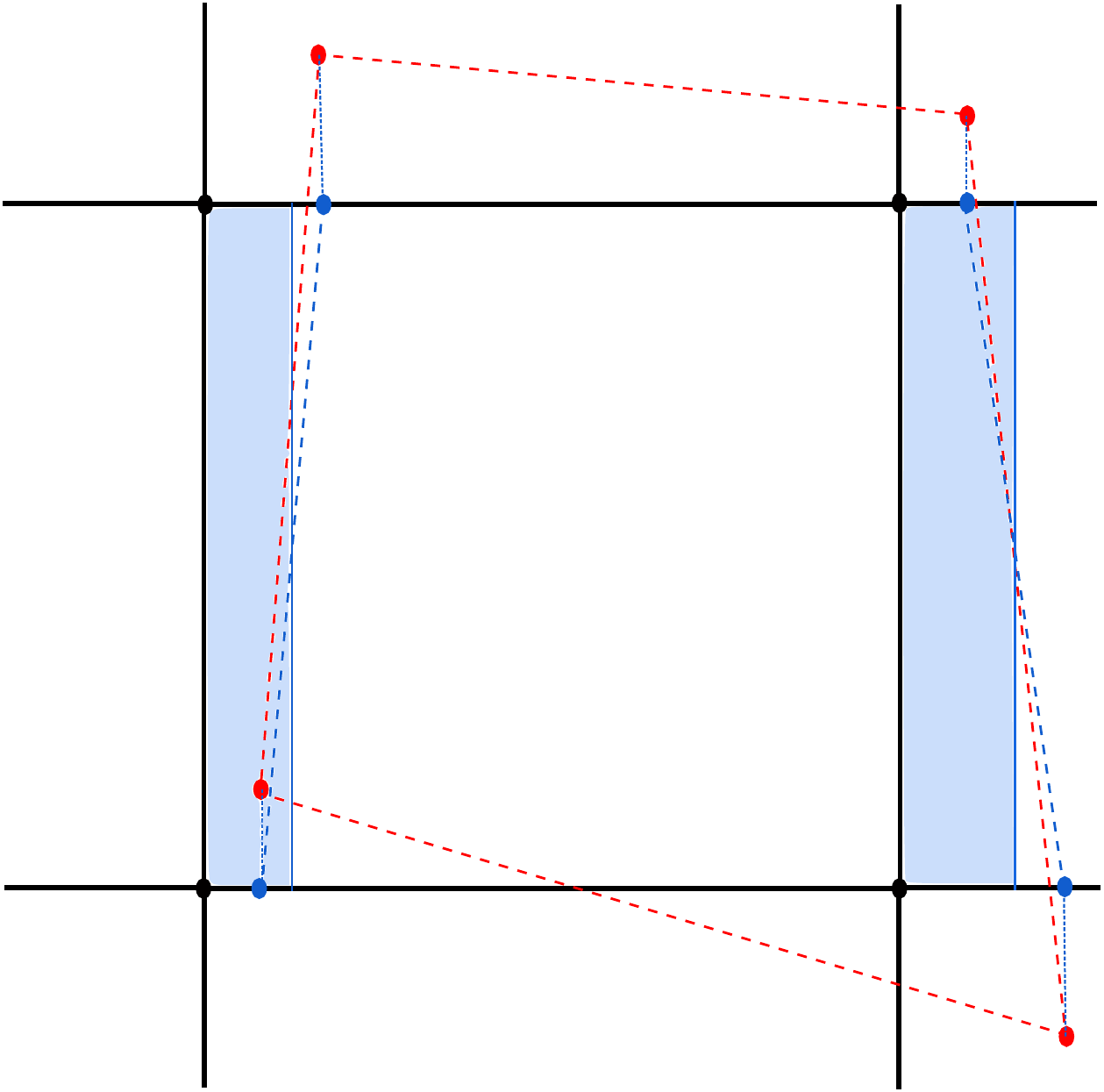}
\caption{Computation of volume fluxes at $X$-faces.}
\label{fig:adi_projection_x}
\end{center}
\end{figure}

\hspace{1.5em}Let us first consider the Lagrangian displacement along $X$, computed from horizontal node predicted velocities $\mathbf{u}_p^{n+1/2}\cdot\mathbf{e}_x$, according to \eqref{correction}, see Figure~\ref{fig:adi_projection_x}. Blue dots define the positions of $X$-Lagrangiannodes. Thus, the volume fluxes considered by this $X$-remap are the ones delimited by blue dotted lines. They are approximated by the blue rectangles, drawn from the mid-points of those blue dotted segments.

Formally, in cell $(i,j)$, the respective volume fluxes at left ($l = i-1/2,j$) and right ($r = i+1/2,j$) faces write:
\begin{equation*}
dVol_{i-1/2,j} = \mathbf{u}_{i-1/2,j}^{n+1/2}\cdot\mathbf{e}_x \: \Delta t^n \Delta y := \frac{1}{2}\big(\mathbf{u}_{i-1/2,j-1/2}^{n+1/2} + \mathbf{u}_{i-1/2,j+1/2}^{n+1/2} \big)\cdot\mathbf{e}_x\: \Delta t^n \Delta y
\end{equation*}
\begin{equation*}
dVol_{i+1/2,j} = \mathbf{u}_{i+1/2,j}^{n+1/2}\cdot\mathbf{e}_x \: \Delta t^n \Delta y:= \frac{1}{2}\big(\mathbf{u}_{i+1/2,j-1/2}^{n+1/2} + \mathbf{u}_{i+1/2,j+1/2}^{n+1/2} \big)\cdot\mathbf{e}_x\: \Delta t^n \Delta y
\end{equation*}

\noindent Those values enable to compute a 1D representation of the Lagrangian volume and density at this step (denoted $t^{lagx}$) as follows, which will be the approximation of the $X$-remap: 

\begin{equation*}
	Vol^{lagx} = Vol^n + dVol_l - dVol_r \:,
\end{equation*}
\begin{equation*}
	\rho^{lagx} = \rho^n  \frac{Vol^n}{Vol^{lagx}} \:\:\:\: \mbox{and} \:\:\:\: e^{lagx} = e^{lag} \:.
\end{equation*}

Let us now compute the variable fluxes at each face. Given the volume flux at a face, one only has to define a value of the variable at this face to get the variable flux. Such a definition can be based on first or second order reconstruction (or even higher order but this will not be discussed in the present work).

\paragraph{First order: upwind reconstruction}
$\:$

Let $a$ be the cell-centred variable on the Lagrangian mesh ($\rho$, $e$, or $\mathbf{u}$ on the Lagrangian dual mesh) to reconstruct on the Eulerian one. At $X$-face $l = (i-1/2,j)$, the upwind reconstruction writes:
\begin{equation*}
a_{i-1/2,j}^{upw} = a_c \:\:\: \mbox{where} \:\:\: c = \left\lbrace
\begin{array}{l r}
(i-1,j) & \mbox{if} \:\:\: dVol_{i-1/2,j} > 0 \\
(i,j) & \mbox{if} \:\:\: dVol_{i-1/2,j} < 0
\end{array}
\right.
\end{equation*}

\paragraph{Second order: linear reconstruction}
$\:$

In order to build a second order scheme, the remap phase needs to be second order as well. But second order spatial reconstruction implies the computation of a gradient. Since the solution can be discontinous, the use of slope limiters when computing discrete gradient is necessary, and avoids to get spurious oscillations. In mathematical terms, taking $c$ defined as in the upwind reconstruction, the value of $a$ at $X$-face  $l = (i-1/2,j)$ is:
\begin{equation*}
a_{i-1/2,j}^{o2} = a_c + \frac{{\delta}^{o2}_x a_c}{2}\:\big(\:\mbox{sgn}(dVol_{i-1/2,j}) \: {\Delta x}_c^{lag} - \Delta t^n \: \mathbf{u}_{i-1/2,j}^{n+1/2}\cdot\mathbf{e}_x\: \big)
\end{equation*}
where ${\Delta x}_c^{lag} = \Delta x + \Delta t^n \:(\mathbf{u}_l^{n+1/2} - \mathbf{u}_r^{n+1/2})\cdot\mathbf{e}_x$ is the width of the Lagrangian cell, and ${\delta}^{o2}_x a_c$ represents the limited value of the gradient in cell $c$. For example, if $c = (i,j)$, and $x_c^{lag}$ denotes the $x$ coordinate of the centre of Lagrangian cell $c$:

\begin{equation*}
{\delta}^{o2}_x a_{i,j} = \frac{1}{2} \: \big( \: \Phi_{i,j}^+\:{\delta}_x a_{i+1/2,j} + \Phi_{i,j}^- \:{\delta}_x a_{i-1/2,j} \: \big)
\end{equation*}

\begin{equation*}
\mbox{with} \:\:\:\: {\delta}_x a_{i\pm1/2,j} = \frac{a_{i\pm1,j} - a_{i,j}}{x_{i\pm1,j}^{lag} - x_{i,j}^{lag}}
\end{equation*}

\begin{equation*}
\mbox{and} \:\:\:\:\: \Phi_{i,j}^+ = \varphi (r)\:, \:\:\: \Phi_{i,j}^- = \varphi (1/r)\:, \:\:\: r = \frac{{\delta}_x a_{i-1/2,j}}{{\delta}_x a_{i+1/2,j}}
\end{equation*}

\noindent The limiter function $\varphi$ makes the scheme degenerate to order 1 when strong gradients are detected, see \cite{courtesinglard2012schemas} for further details. The choice made here is Van Leer limiter:
\begin{equation*}
\varphi(r) = \frac{r + \lvert r \rvert}{1 + r}
\end{equation*}

Now the reconstructed values at $X$-faces have been determined, variable fluxes can be computed easily, and variables can be remapped. This is the $X$-remap phase, denoted $t^{projx}$. According to previous notations, in cell $c$, the conservation of cell-centred variables leads to:

\begin{equation*}
	m_c^{projx} = m_c^{lagx} + \rho_l^{lagx \:o2} dVol_l - \rho_r^{lagx \:o2} dVol_r 
\end{equation*}
\vspace{0.05cm}
\begin{equation*}
	m_c^{projx}e_c^{projx} = m_c^{lagx}e_c^{lagx} + \rho_l^{lagx \:o2}e_l^{lagx \:o2} dVol_l - \rho_r^{lagx \:o2}e_r^{lagx \:o2} dVol_r
\end{equation*}
\vspace{0.05cm}

\textit{Remark:} Instead of reconstructing the quantity $\rho e$ at the faces, the variables $\rho$ and $e$ are independently reconstructed, for robustness issues.

\begin{figure}[h]
\begin{center}
\includegraphics[width=0.65\textwidth, height=0.5\textwidth]{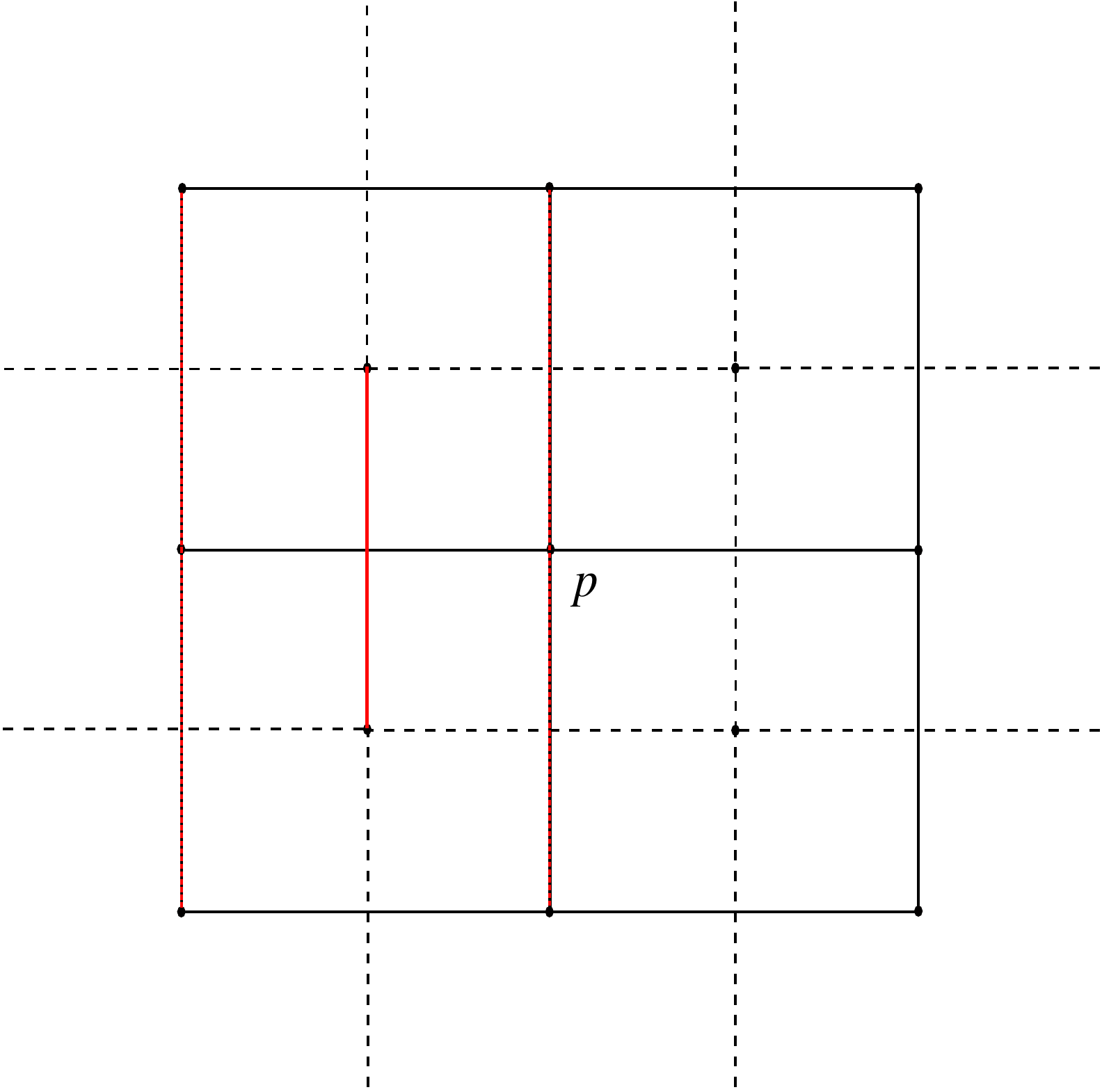}
\caption{Primary and dual meshes around node $p$.}
\label{fig:dual_mesh}
\end{center}
\end{figure}

Regarding node-centred variables, i.e.\ $x$ and $y$ components of the velocities, the principle is exactly the same, but on the dual mesh. That is to say, variables are reconstructed (at first or second order) at the faces of the dual mesh, and then remapped. The nodes of the dual mesh are defined as the centres of the primary mesh, see Figure~\ref{fig:dual_mesh}. The black lines and dots respectively draw the primary and dual meshes. Here, in the regular orthogonal case, the dual mesh remains regular orthogonal.

\vspace{0.5cm}

Let us use the index $_p$ to denote objects related to the dual cell (also called nodal cell) containing node $p$. For example, the faces of this dual cell $p$ are the $f_p$. The nodal mass and mass fluxes at its faces write:

\begin{equation*}
	m_p = \frac{1}{4} \sum_{cell \in \{p\}} \!\!\! m_{cell}
\end{equation*}		
\begin{equation*}
	dm_{f_p} = \frac{1}{4} \sum_{cell \in \{p\}} \!\!\! (dm_f)_{cell}
\end{equation*}

\noindent For example, the mass flux at the left face of the dual mesh, corresponding to $f_p = l_p$,  which is in red in Figure~\ref{fig:dual_mesh}, is computed by taking the average algebraic value of the mass fluxes at the left faces of each cell $c$ containing node $p$, that are all in red as well on the figure. \\
Thus, when writing momentum conservation at node $p$, one obtains:

\begin{equation*}
	m^{projx}_p \mathbf{u}^{projx}_p = m^{lagx}_p \mathbf{u}^{lagx}_p + dm_{l_p} \mathbf{u}^{lagx \:o2}_{l_p} - dm_{r_p} \mathbf{u}^{lagx \:o2}_{r_p}
\end{equation*}
		
\noindent which gives the velocities at $t^{projx}$, and ends the $X$-remap phase.

\subsubsection{$Y$-Remap phase}

\begin{figure}[h]
\begin{center}
\includegraphics[width=0.5\textwidth, height=0.4\textwidth]{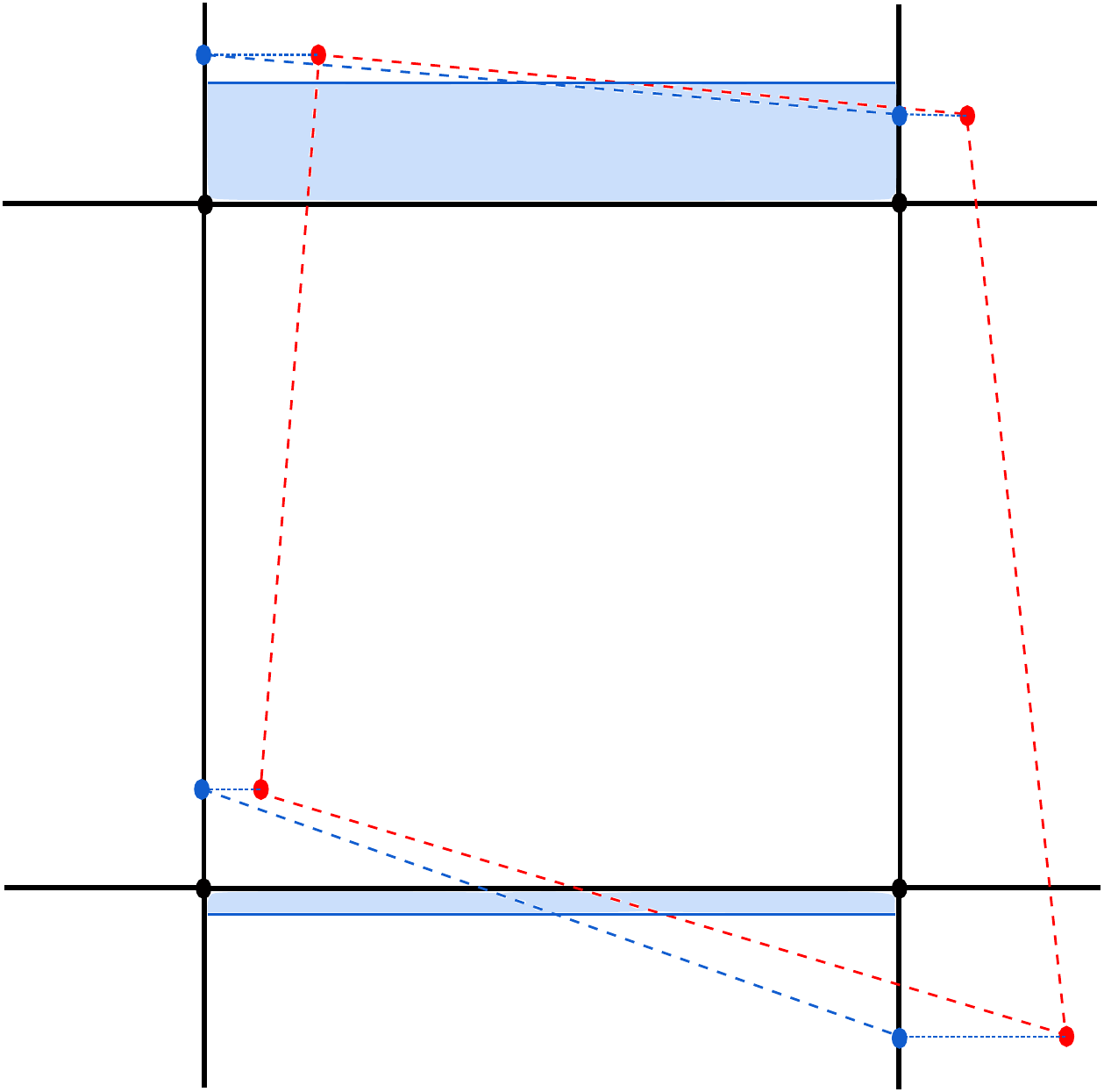}
\caption{Computation of volume fluxes at $Y$-faces.}
\label{fig:adi_projection_y}
\end{center}	
\end{figure}

\hspace{1.5em}Details of the computations will not be written because they are strictly equivalent to those performed in the $X$-remap phase. The Lagrangian displacement along Y is only computed from vertical node predicted velocities $\mathbf{u}_p^{n+1/2}\cdot\mathbf{e}_y$, see Figure~\ref{fig:adi_projection_y}. The only difference is that, instead of starting computations with values of the variables at $t^{lag}$ like for the $X$-remap, this phase starts with values at $t^{projx}$, i.e.\ values of the variables after the $X$-remap has been performed.

\vspace{0.5cm}	

The computation of volume fluxes at the top $t = (i,j+1/2)$ and bottom $b = (i,j-1/2)$ faces enable to write at $t^{lagy}$, taking into account that $Vol^{projx} = Vol^n$:

\begin{equation*}
	Vol^{lagy} = Vol^n + dVol_b - dVol_t \:,
\end{equation*}
\begin{equation*}
\rho^{lagy} = \rho^n  \frac{Vol^n}{Vol^{lagy}} \:\:\:\: \mbox{and} \:\:\:\: e^{lagy} = e^{projx} \:.
\end{equation*}

\noindent Then $Y$-remap (after reconstruction of $Y$-faces values) on both primary and dual meshes provides the values of all variables at $t^{projy} = t^{n+1}$.

Let us recap, at each iteration, once the Lagrangian phase is done, the ADI remap splits in four successive steps: 
\begin{enumerate}
\item Lagrangian displacement of the mesh along $X$ ($lagx$)
\item Remap of Lagrangian variables ($projx$)
\item Lagrangian displacement of the mesh along $Y$ ($lagy$)
\item Remap of Lagrangian variables ($projy =$ iteration $n+1$)
\end{enumerate}

\begin{figure}[h]
\begin{center}
\includegraphics[width=0.2\textwidth, height=0.17\textwidth]{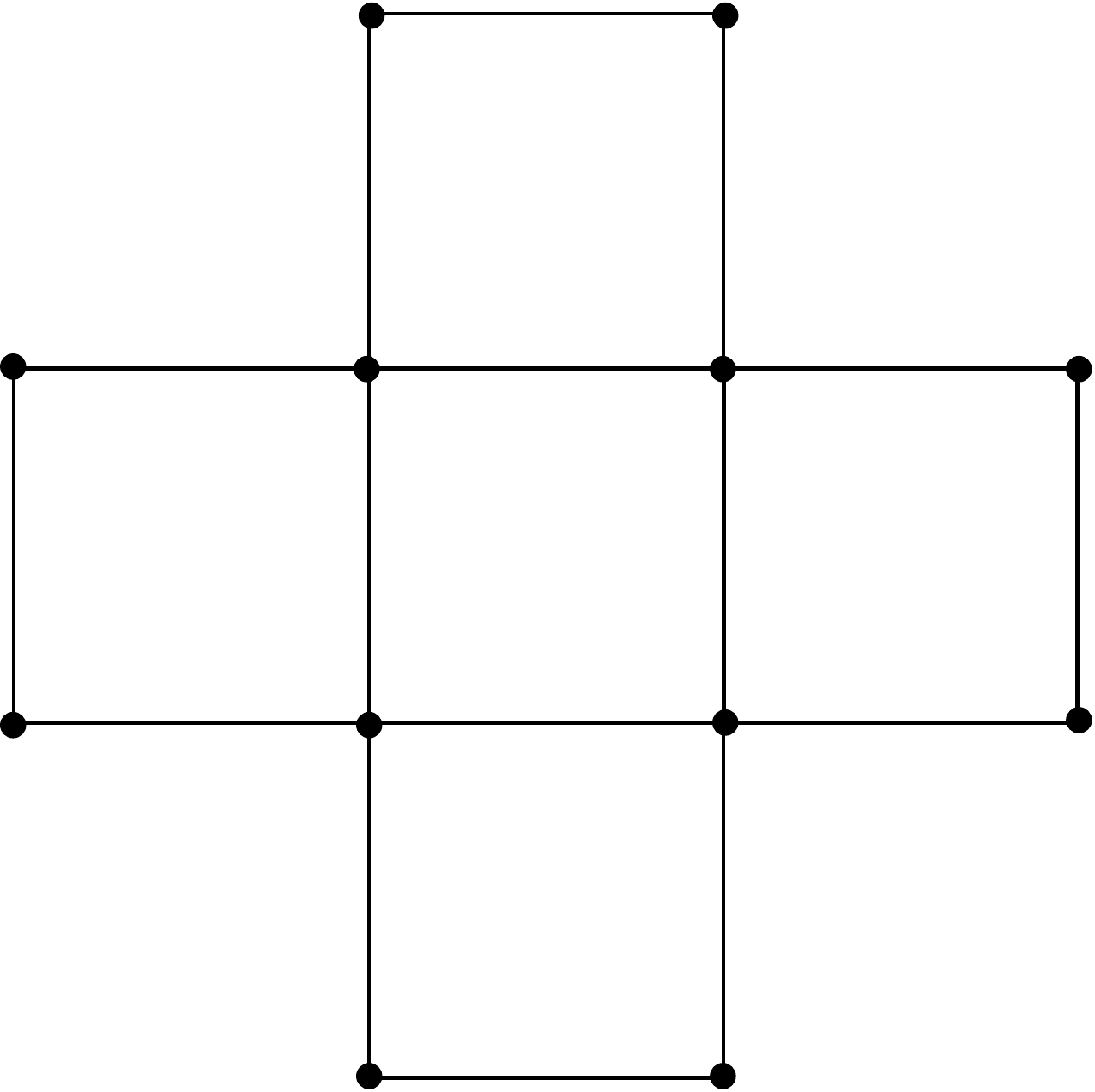}\hspace{3cm}
\includegraphics[width=0.2\textwidth, height=0.17\textwidth]{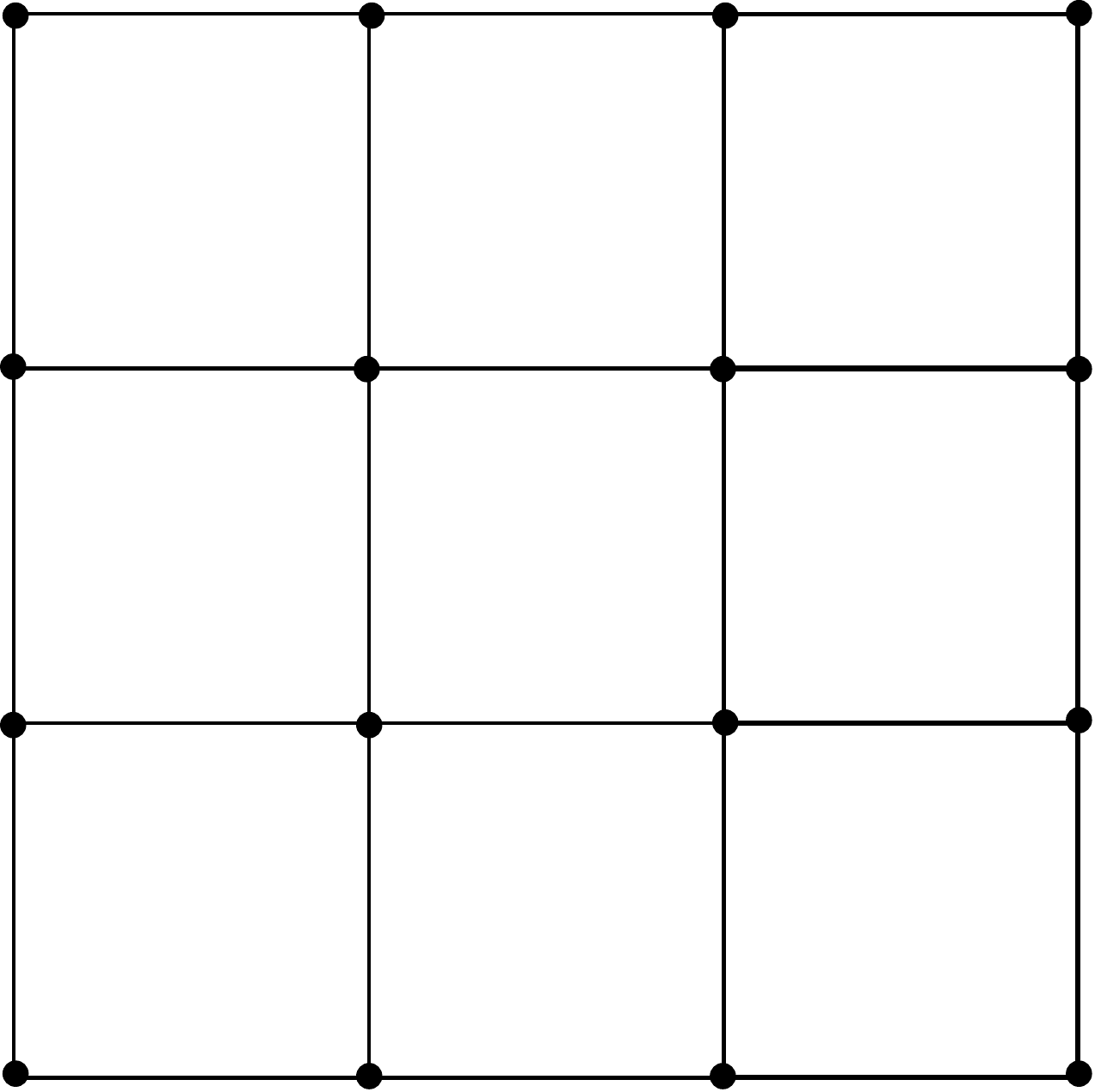}
\caption{Stencils of 5 points remaps (left) and 9 points remaps (right).}
\label{fig:remap_stencil}
\end{center}	
\end{figure}

\textit{Remark 1:} Because of the unidirectional way fluxes are computed and remapped, one can think at first that the AD remap is a 5 points remap. That is to say, the stencil is the cross-shaped group of cells containing the cell at the centre and its four first neighbours, see Figure~\ref{fig:remap_stencil}. Actually, thanks to the directional splitting, matter can be exchanged through the corners. For example, from cell $(i,j)$, it is possible to transfer matter to cell $(i+1,j)$ through face $(i+1/2,j)$ during the $X$-remap phase, and then transfer a part of this matter to cell $(i+1,j+1)$ through face $(i+1,j+1/2)$ during the $Y$-remap phase. Somehow, the directional splitting computes implicit corner fluxes. That is why the AD remap must be thought as a 9 points remap (star-shaped group of 5 plus the cells in the 4 diagonal directions, as shown Figure~\ref{fig:remap_stencil}). 

\vspace{0.5cm}	

\textit{Remark 2:} Let us count the number of MPI communications required during one time step: one synchronization for the Lagrangian phase, one for the $X$-remap phase, and another one for the $Y$-remap phase. That is a total of \textbf{3 in 2D}. 

\textit{Remark 3:} A fixed order of treating $X$ and $Y$ directions necessarily introduces a loss of symmetry in the geometry, especially in cases where the solution is supposed to be symmetric with respect to $x=y$. In order to diminish this effect, another alternating is processed: $X$-$Y$ remap at odd iterations and $Y$-$X$ remap at even iterations. It also produces an almost second order in time remap scheme because it is similar to a Strang splitting when regarding two successive iterations.

\subsection{Direct Remap}

\begin{figure}[h]
\begin{center}
\includegraphics[width=0.5\textwidth, height=0.4\textwidth]{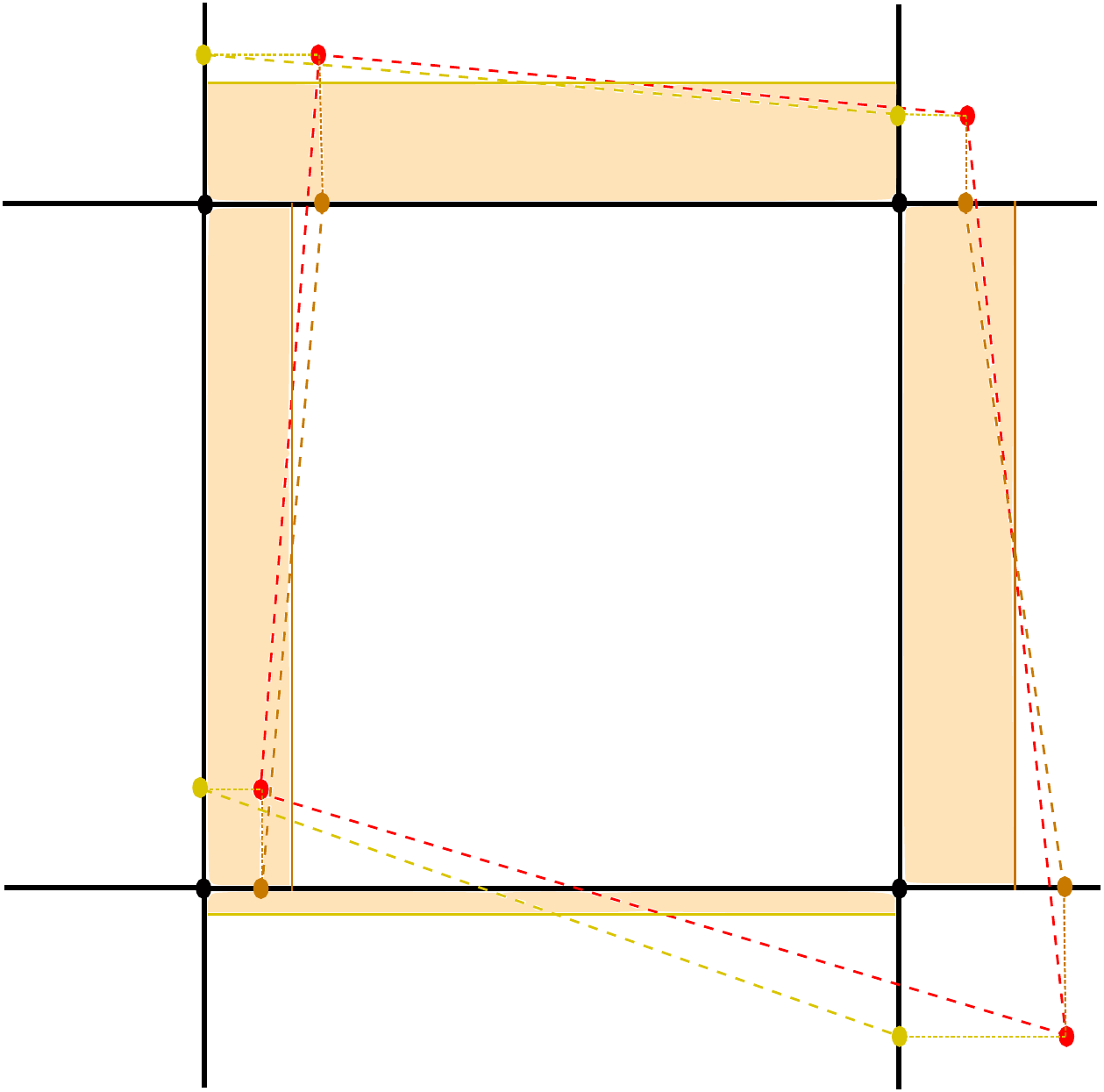}
\caption{Computation of volume fluxes at both $X$ and $Y$-faces.}
\label{fig:direct_projection}
\end{center}	
\end{figure}

\hspace{1.5em}The name "Direct" means that fluxes at the $X$ and $Y$-faces are independently computed and remapped, but at the same time. This is a one-step remap, see \cite{debray2013projection}. The Lagrangian displacement along both directions $X$ and $Y$ and the fluxes are computed respectively from $\mathbf{u}_p^{n+1/2}\cdot\mathbf{e}_x$ and $\mathbf{u}_p^{n+1/2}\cdot\mathbf{e}_y$, exactly the same way as for the ADI remap, see Figure~\ref{fig:direct_projection}.

\vspace{0.5cm}	

Orange and Yellow dots refer to $X$ and $Y$ displacements, respectively. The values of volume fluxes at each face $l$, $r$, $t$, $b$ lead to the following Lagrangian quantities after the grid has been displaced, at $t^{lagxy}$:		

\begin{equation*}
	Vol^{lagxy} = Vol^n + dVol_l - dVol_r + dVol_b - dVol_t ,
\end{equation*}
\begin{equation*}
	\rho^{lagxy} = \rho^n  \frac{Vol^n}{Vol^{lagxy}} \:\:\:\: \mbox{and} \:\:\:\: e^{lagxy} = e^{lag} \:.
\end{equation*}

\noindent The reconstruction of face variables is also the same as for the AD remap, i.e.\ variables at $X$ and $Y$-faces are interpolated using 1D second order limited linear reconstructions in both $X$ and $Y$ directions. When writing mass and energy conservations, one gets: 
		
\begin{equation*}
	m_c^{proj} = m_c^{lagxy} + \sum_{faces} \!\rho_f^{lagxy \:o2} dVol_f
\end{equation*}
\begin{equation*}
	m_c^{proj} e_c^{proj} = m_c^{lagxy} e_c^{lagxy} + \sum_{faces} \!\rho_f^{lagxy \:o2} e_f^{lagxy \:o2} dVol_f 
\end{equation*}

\noindent Similarly, momentum conservation on the dual mesh writes:

\begin{equation*}
	m^{projx}_p \mathbf{u}^{projx}_p = m^{lagx}_p \mathbf{u}^{lagx}_p + \sum_{faces} \! dm_{f_p} \mathbf{u}^{lagx \:o2}_{f_p}
\end{equation*}
	
As for the previous type of remap, let us recap the successive steps of the Direct remap, once the Lagrangian phase is over:	
\begin{enumerate}
\item Lagrangian displacement of the mesh along $X$ and $Y$ ($lagxy$)
\item Remap of Lagrangian variables ($proj =$ iteration $n+1$)
\end{enumerate}

\textit{Remark 1:} Unlike the AD remap, since the remap is in one step, there is no other flux than the explicit fluxes at $X$ and $Y$-faces. This remap is a strict 5 points remap with a star-shaped stencil, just like a classical Finite-Volume scheme. Since there are no corner fluxes, corner effects will not be taken into account by the scheme. This might lead to an accuracy not quite as good as the AD remap, even in simple cases like the linear advection in the $x=y$ direction. 

\textit{Remark 2:} Thanks to one-step remap, the number of MPI communications reduces to \textbf{2 in 2D} (one at each phase), and in fact 2 whatever the dimension.

\subsection{Direct Remap with Corner Fluxes}

\hspace{1.5em}The initial idea is to build a type of remap that combines the advantages of both remaps presented before. That is to say, a one-step remap that takes into account corner effects. Those two requirements together imply necessarily the definition of corner fluxes, since one-step face fluxes cannot report corner effects. 

\begin{figure}[h]
\begin{center}
\includegraphics[width=0.5\textwidth, height=0.4\textwidth]{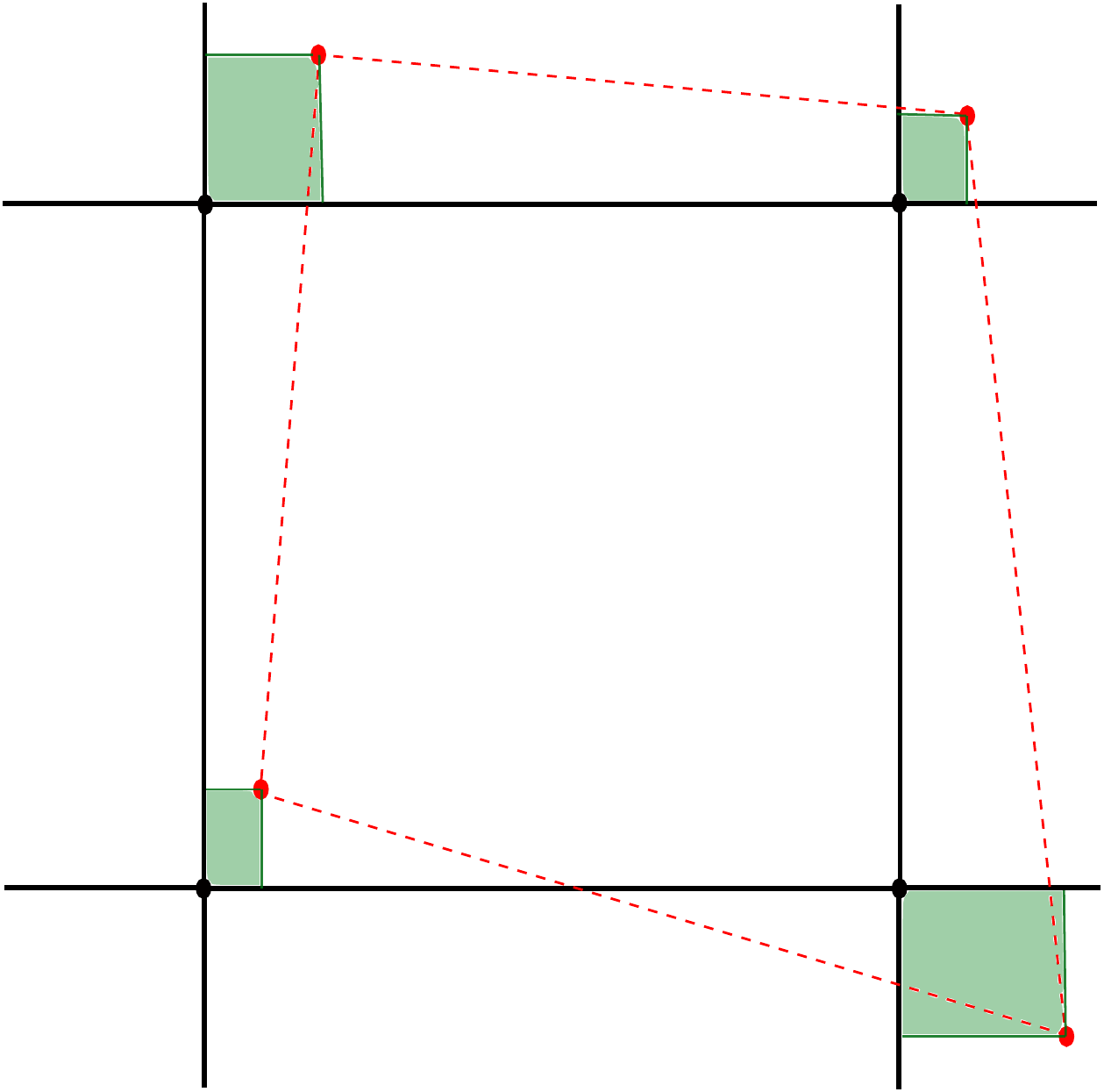}
\caption{Definition of corner volume fluxes.}
\label{fig:corner_fluxes}
\end{center}	
\end{figure}

For the corner flux at a given node $p$ to be defined unequivocally, it has to depend only on node quantities. The displacement of the node seems to be the most natural way to define the volume flux at this node. In the same way as what has been done for face fluxes, a rectangular approximation of the corner flux is computed, setting that the vector $\mathbf{u}_p^{n+1/2} \Delta t^n$ represents the diagonal of the rectangular corner flux, see Figure~\ref{fig:corner_fluxes}. 

\vspace{0.5cm}	

\noindent Thus, the volume flux at each corner $p$ can be expressed as follows:
\begin{equation*}
dVol_p = \lvert \mathbf{u}_p^{n+1/2} \cdot \mathbf{e}_x \rvert \Delta t^n \times \lvert \mathbf{u}_p^{n+1/2} \cdot \mathbf{e}_y\rvert \Delta t^n \:\:\: > 0
\end{equation*}

\noindent Setting  $\:dx_p^{n+1/2} = ( \mathbf{u}_p^{n+1/2} \cdot \mathbf{e}_{x} ) \Delta t^n\:\:$ and $\:\:dy_p^{n+1/2} = (\mathbf{u}_p^{n+1/2} \cdot \mathbf{e}_{y} ) \Delta t^n\:$, one obtains
\begin{equation*}
dVol_p = \lvert \: dx_p^{n+1/2} \times dy_p^{n+1/2}  \:\rvert
\end{equation*}

\textit{Remark:} Corner fluxes are completely defined from a positive value $dVol_p$,  and both signs of $dx_p^{n+1/2}$ and $dy_p^{n+1/2}$. These signs enable to identify the donor cell and the receiving cell.

\begin{figure}[h]
\begin{center}
\includegraphics[width=0.45\textwidth, height=0.36\textwidth]{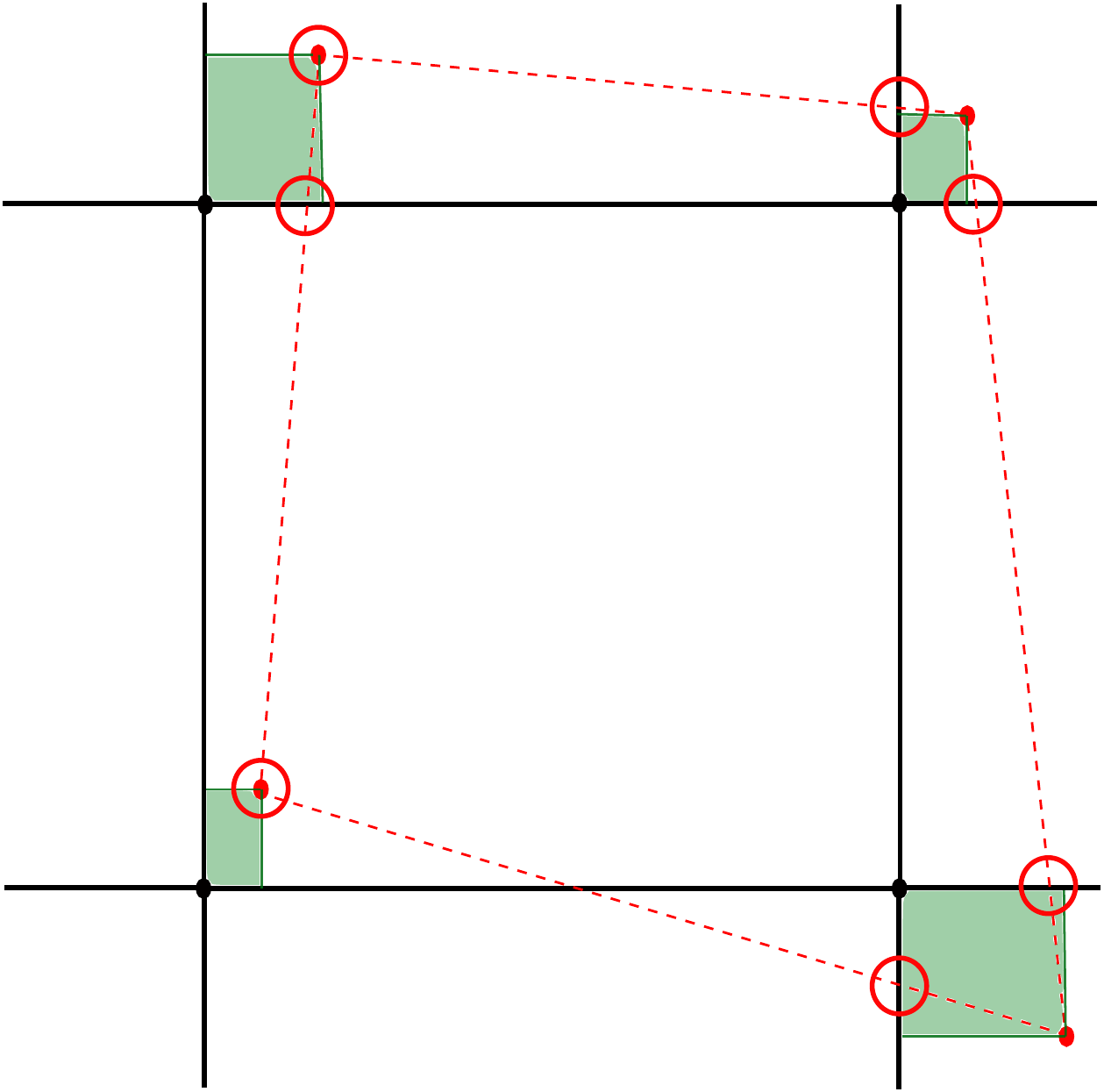}\hspace{0.05\textwidth}
\includegraphics[width=0.45\textwidth, height=0.36\textwidth]{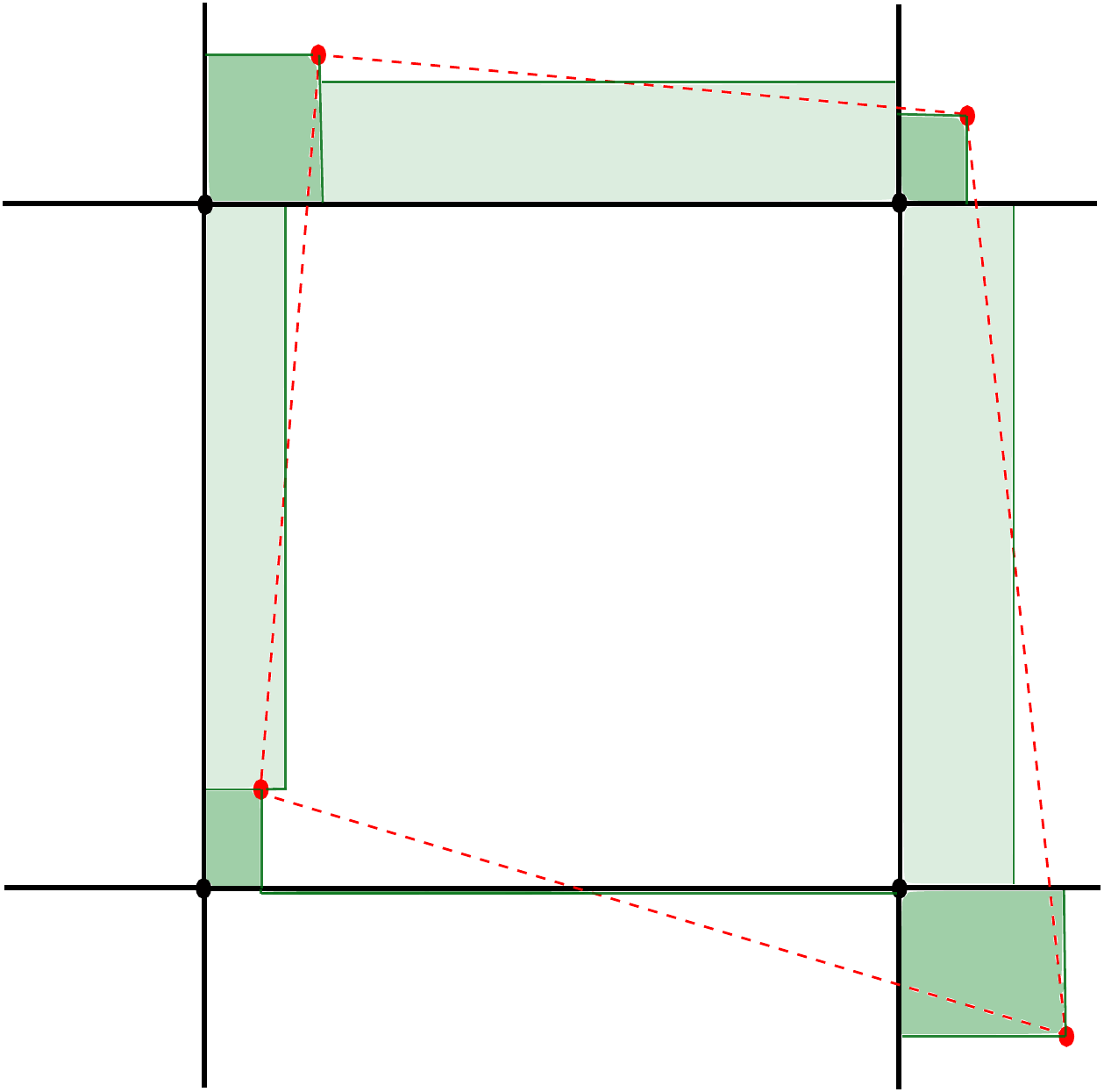}
\caption{Definition of volume fluxes at both $X$ and $Y$-faces.}
\label{fig:face_fluxes_cf}
\end{center}	
\end{figure}

Now corner fluxes have been well defined, one needs to compute fluxes at the faces of the cell. Here again, according to what is done in other remaps, a natural choice is the rectangular approximation of the trapezoid, see Figure~\ref{fig:face_fluxes_cf}. Note that at each face, the volume flux is well defined since it only depends on the velocities at its nodes, and not on the donor cell. This is mandatory for the sake of conservation.

\vspace{0.5cm}	

\noindent A general methodology can be given for face volume fluxes computation. At each face, one has to determine the coordinates $(x_f^-,y_f^-)$, $(x_f^+,y_f^+)$ of the vertices of the trapezoid, i.e.\ the red circles in Figure~\ref{fig:face_fluxes_cf}. For example, at any $X$-face $f$ which extremities are $p^-$ and $p^+$, we start computing $y_f^-$ and $y_f^+$ using this simple formula:
\begin{equation*}
y_f^- = y_{p^-}^n + \max (0, dy_{p^-}^{n+1/2}) \:, \:\:\:\:\:\: y_f^+ = y_{p^+}^n + \min (0, dy_{p^+}^{n+1/2})
\end{equation*}
Then $x_f^-$ and $x_f^+$ are determined thanks to the equation of the line corresponding to Lagrangian face $f$, computed from $\mathbf{x}_{p^\pm}^{n+1/2}$. This finally leads to the following volume flux at $X$-face $f$:
\begin{equation*}
dVol_f = \frac{1}{2} ( x_f^- + x_f^+ ) \times ( y_f^+ - y_f^- )
\end{equation*}

\textit{Remark:} In this remap, the volume flux at each face depends on both $x$ and $y$ components of the velocities at its extremities.

\begin{figure}[h]
\begin{center}
\includegraphics[width=0.5\textwidth, height=0.4\textwidth]{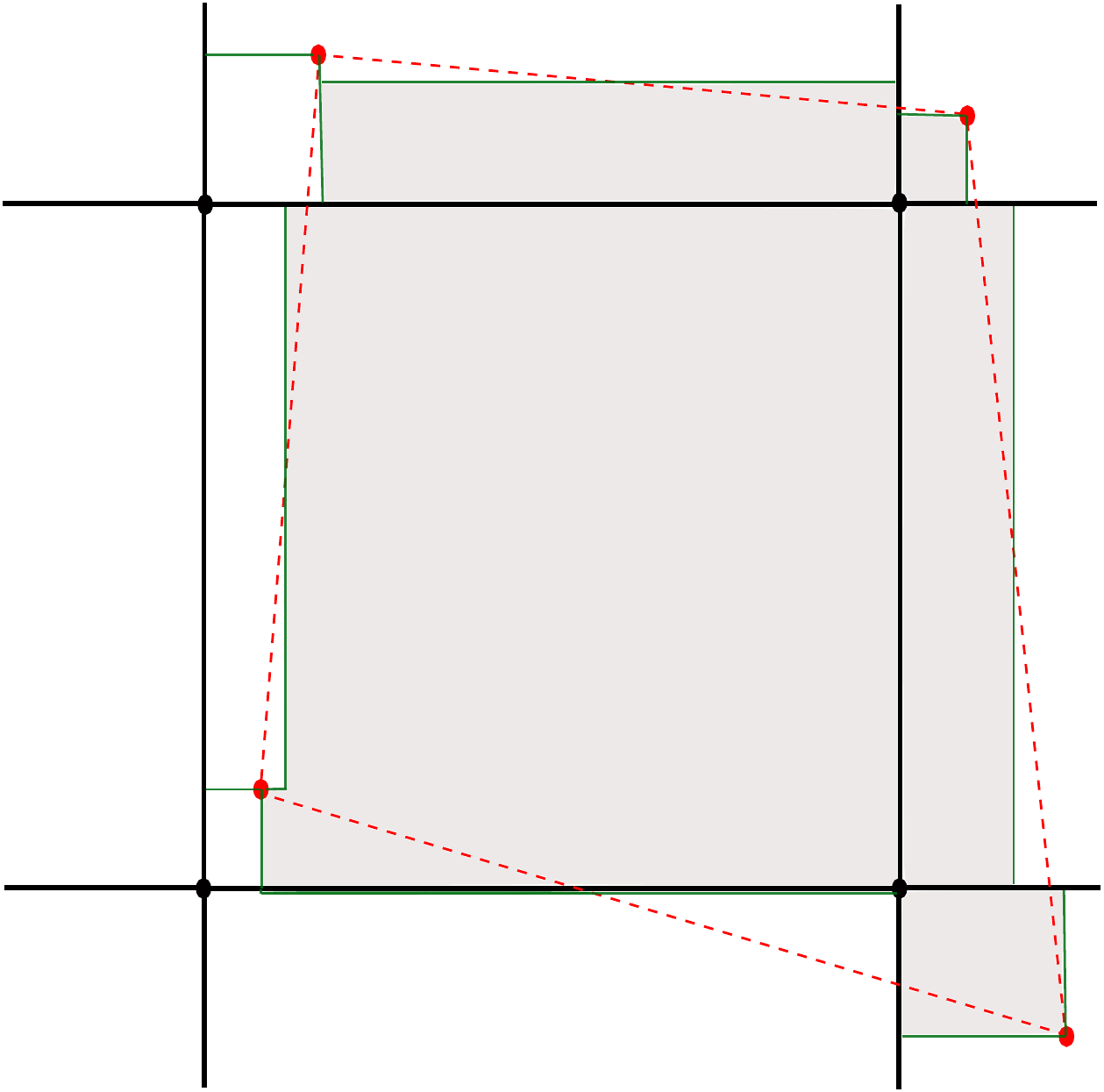}
\caption{Lagrangian representation of the volume.}
\label{fig:vol_lag_cf}
\end{center}	
\end{figure}

Given volume fluxes at the faces and at the corners, one is able to compute Lagrangian geometric quantities such as the volume (see Figure~\ref{fig:vol_lag_cf}) as well as thermodynamical quantities at time $t^{lagCF}$:

\begin{equation*}
	Vol^{lagCF} = Vol^n + \sum_{faces}\!dVol_f + \!\!\sum_{corners}\!\!\!dVol_p \:,
\end{equation*}
\begin{equation*}
	\rho^{lagCF} = \rho^n  \frac{Vol^n}{Vol^{lagCF}}\:\:\:\: \mbox{and} \:\:\:\: e^{lagCF} = e^{lag} \:.
\end{equation*}

\vspace{0.5cm}	

In this case, values have to be defined at faces, but also at corners. As for the faces, a choice has to be made between upwind reconstruction and limited linear reconstruction. The importance of this choice and the way of defining second order values at the corners will be discussed in section 3.2. Here, those reconstructed values at the corners are assumed to be known. The equations of the remap on the primary mesh write in cell $c$:
		
\begin{equation*}
	m_c^{proj} = m_c^{lagCF} + \sum_{faces} \rho_f^{lagCF \:o2} dVol_f + \!\!\sum_{corners} \!\!\rho_p^{lagCF \:o2} dVol_p
\end{equation*}		

\begin{equation*}
\left.
\begin{array}{l c c}
	m_c^{proj}e_c^{proj} & = & m_c^{lagCF}e_c^{lagCF} + \displaystyle\sum_{faces} \rho_f^{lagCF \:o2} e_f^{lagCF \:o2} dVol_f \\
	\: & \: & \: \\
	\: & \: &  + \!\!\displaystyle\sum_{corners} \!\!\rho_p^{lagCF \:o2} e_p^{lagCF \:o2} dVol_p
\end{array}
\right.
\end{equation*}
\vspace{0.1cm}	

For obvious reasons of compatibility between the primary and dual meshes, it is essential to define fluxes at the corners $p'_p$ of the dual mesh, i.e.\ at the centres of the cells of the primary mesh $c = p'_p$.  There are several ways to define those corner fluxes. Note that as it is the case for the faces, only mass corner fluxes need to be defined to remap $x$ and $y$ components of velocities. In order to stay as consistent as possible with what is done at the faces (i.e.\ in the AD remap), one must express nodal mass corner fluxes as follows:

\begin{equation*}
	dm_{p'_p} = \frac{1}{4} \sum_{cell \in \{p\}} \!\!\! (dm_{p'})_{cell}
\end{equation*} 

\vspace{0.5cm}	
\begin{figure}[h]
\begin{center}
\includegraphics[width=0.42\textwidth, height=0.37\textwidth]{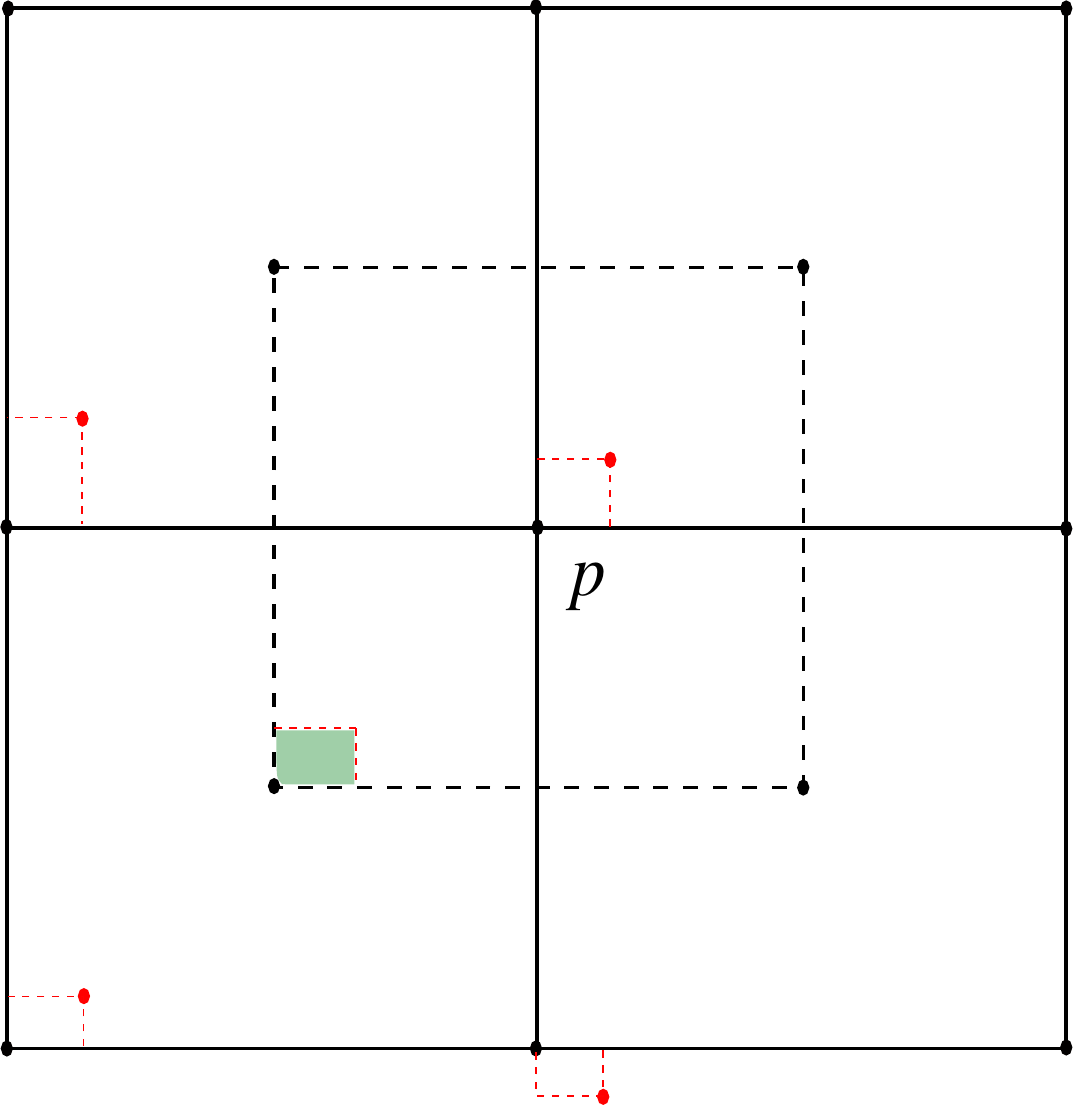}
\caption{Definition of corner fluxes on the dual mesh.}
\label{fig:nodal_corner_fluxes}
\end{center}
\end{figure}

\noindent For instance, if $p'$ stands for left-bottom, like in Figure~\ref{fig:nodal_corner_fluxes}, the mass corner flux at the left-bottom node of the dual cell $p$ is defined as the average of all mass corner fluxes at left-bottom nodes of the cells $cell \in \{p\}$.

\textit{Remark:} In this particular case, the left-bottom mass corner flux in the cell situated at the right-bottom of node $p$, say cell $p_{+-}$, equals zero. Indeed, the direction of this corner flux makes it equal zero for this cell, and non zero for the cell situated at the left-bottom of node $p$, say $p_{--}$. Thus, this flux will be counted as the flux at the right-bottom node of cell $p_{--}$ when computing the right-bottom mass corner flux in the dual cell.

Regarding mass fluxes at the faces, they are computed similarly to AD or Direct remaps. All mass fluxes being known, let us write momentum conservation equation:
\begin{equation*}
	m^{proj}_p \mathbf{u}^{proj}_p = m^{lagCF}_p \mathbf{u}^{lagCF}_p + \sum_{faces}\! dm_{f_p} \mathbf{u}^{lagCF \:o2}_{f_p} + \!\!\sum_{corners}\!\!\! dm_{c} \mathbf{u}^{lagCF \:o2}_{c}
\end{equation*}

Here is a brief recap of the successive steps of the Direct Remap with Corner Fluxes (denoted DirectCF), after the end of the Lagrangian phase:	
\begin{enumerate}
\item Lagrangian displacement of the mesh along $X$, $Y$ and diagonal directions ($lagCF$)
\item Remap of Lagrangian variables ($proj =$ iteration $n+1$)
\end{enumerate}

\textit{Remark 1:} The desired properties seem to be recovered. Indeed, in terms of MPI communications, the one-step aspect of this remap implies \textbf{2 in 2D}. It also implies a 
\textbf{higher arithmetic intensity} than the AD remap. Plus, the corner effects are obviously taken into account, in a way that appears to be more faithful than the AD remap since corner fluxes are computed directly from the displacement of the nodes.

\textit{Remark 2:} Geometric computations are more complex than in the AD case.

 \clearpage
 \newpage
 \section{Specificities of the Direct Remap with Corner Fluxes}

\hspace{1.5em}Before implementing this new remap, it is natural to wonder in what extent could it be beneficial, i.e.\ whether it is worth performing those complex computations, and  how do these corners affect the initial scheme.

\subsection{Specific case: linear advection with first order reconstruction for the remap}

\hspace{1.5em}Let us first compare the DirectCF and AD remaps in the most basic situation: linear advection, which corresponds to the case of a constant velocity field $\mathbf{u} = (u_x, u_y)$. The remaps are performed with a first order, \i.e. upwind, reconstruction of variables at corners and faces. Let $a$ be a scalar volume quantity defined on the primary mesh, the goal is to compare $a^{proj, \:CF}$ and $a^{proj, \:AD}$. Horizontal and vertical displacements of the nodes are respectively defined as $dx = u_x \Delta t$, $dy = u_y \Delta t$. To simplify computations, $dx$ and $dy$ are assumed to be positive. Two dimensionless variables are introduced: $\varepsilon_x = dx / \Delta x$, $\varepsilon_y = dy / \Delta y$.

\paragraph{Direct remap with Corner Fluxes}
$\:$
\vspace{0.3cm}

At each node $p$, $X$-face $f_X$ and $Y$-face $f_Y$, volume fluxes write:
\begin{displaymath}
dVol_p = dx \: dy, \:\:\:\: dVol_{f_X} = dx\:(\Delta y - dy), \:\:\:\: dVol_{f_Y} = dy\:(\Delta x - dx)
\end{displaymath}

\noindent Since the velocity field is constant, the grid is not deformed. In particular, $Vol^{lag} = \Delta x \Delta y$. Moreover, in the pure advection case, $a_{i,j}^{lag} = a_{i,j}$ for all $(i,j)$. Thus, when writing the equation of mass conservation in cell $(i,j)$, one obtains:

\begin{displaymath}
\begin{array}{c c l}
\Delta x \Delta y \:a_{i,j}^{proj} & = & \Delta x \Delta y \:a_{i,j} - dVol_{f_X}(a_{i,j} - a_{i-1,j}) \\
\:  & \: & \: \\
\: & \: & \:\: - dVol_{f_Y}(a_{i,j} - a_{i,j-1}) - dVol_p(a_{i,j} - a_{i-1,j-1})
\end{array}
\end{displaymath}
\begin{displaymath}
\boxed{
\:a_{i,j}^{proj,\:CF}  = a_{i,j}(1 - \varepsilon_x)(1 - \varepsilon_y) + a_{i-1,j}\: \varepsilon_x (1 - \varepsilon_y) +  a_{i-1,j}\: \varepsilon_y (1 - \varepsilon_x) + a_{i-1,j-1}\:\varepsilon_x \varepsilon_y
}
\end{displaymath}

\paragraph{Alternate Directions remap}
$\:$
\vspace{0.3cm}

Let us first consider the $X$-remap. This time, at each $X$-face $f_X$, the volume flux writes:
\begin{displaymath}
dVol_{f_X} = dx\:\Delta y
\end{displaymath}

\noindent  and for all $(i,j)$, mass conservation leads to:

\begin{displaymath}
\Delta x \Delta y \:a_{i,j}^{projx} = \Delta x \Delta y \:a_{i,j} - dVol_{f_X}(a_{i,j} - a_{i-1,j}) 
\end{displaymath}
\begin{displaymath}
\mbox{so that} \:\:\:\:\:\:\:\:\:\: \begin{array}{c c l}
a_{i,j}^{projx} & = & a_{i,j}\: (1 - \varepsilon_x) + a_{i-1,j}\: \varepsilon_x \\
\:  & \: & \: \\
a_{i,j-1}^{projx} & = & a_{i,j-1}\: (1 - \varepsilon_x) + a_{i-1,j-1}\: \varepsilon_x  
\end{array}
\end{displaymath}
  
Then for the $Y$-remap, at each $Y$-face $f_Y$, the volume flux writes:
\begin{displaymath}
dVol_{f_Y} = dy\:\Delta x
\end{displaymath}
\begin{displaymath}
\mbox{thus,} \:\:\:\:\: \Delta x \Delta y \:a_{i,j}^{proj} = \Delta x \Delta y \:a_{i,j}^{projx}  - dVol_{f_Y}(a_{i,j}^{projx} - a_{i,j-1}^{projx}) 
\end{displaymath}
\begin{displaymath}
\boxed{
\:a_{i,j}^{proj,\:AD}  = a_{i,j}(1 - \varepsilon_x)(1 - \varepsilon_y) + a_{i-1,j}\: \varepsilon_x (1 - \varepsilon_y) +  a_{i-1,j}\: \varepsilon_y (1 - \varepsilon_x) + a_{i-1,j-1}\:\varepsilon_x \varepsilon_y
}
\end{displaymath}

The result obtained is the same for both schemes. In other words, the two remaps are strictly equivalent in this specific case. First, it shows that the DirectCF remap degenerates to the right well-known scheme at order 1 in such a simple situation. Second, this feature underlines the importance of second order reconstruction for the DirectCF remap. The accuracy of the results and the differences with AD remap results will highly depend on this reconstruction.

\subsection{Accuracy of the geometric representation for non-linear deformations}

\begin{figure}[h!]
\begin{center}
\includegraphics[width=0.45\textwidth, height=0.35\textwidth]{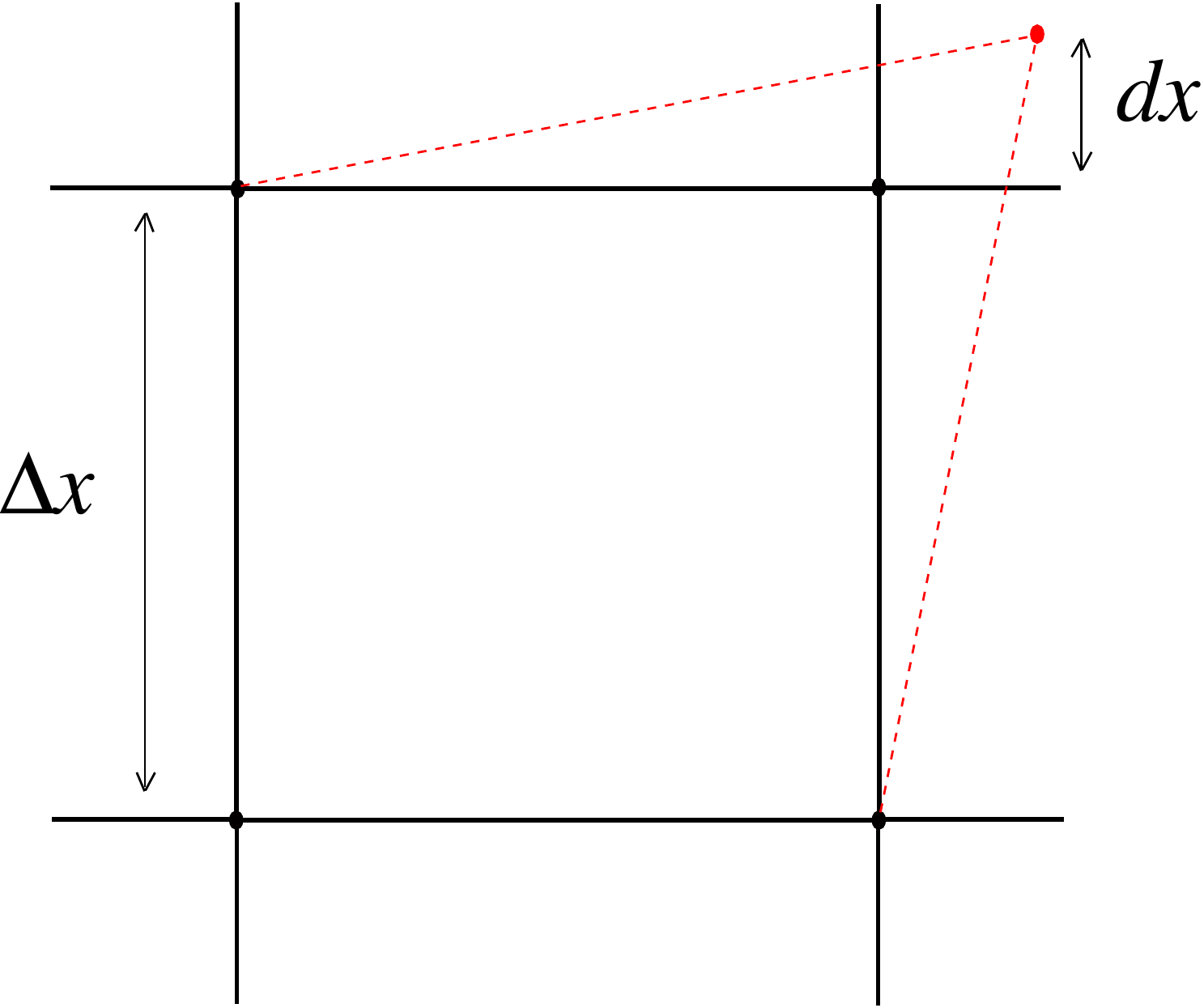}
\caption{Lagrangian displacement of one single node.}
\label{fig:one_node_displacement}
\end{center}
\end{figure}

\hspace{1.5em}Let us now consider a very simple non-linear situation: the displacement of one single node in the direction $x=y$, see Figure~\ref{fig:one_node_displacement}. The initial density equals $1$ in the cell at the centre and $0$ in every other cell. It is assumed that $dx = dy$. In order to estimate the impact of the corner fluxes on both geometric representation and effective mass fluxes, for each remap, the Lagrangian volume will be computed and compared to the exact Lagrangian volume, so will the mass flux at the top right corner.

\begin{figure}[h]
\begin{center}
\includegraphics[width=0.32\textwidth, height=0.26\textwidth]{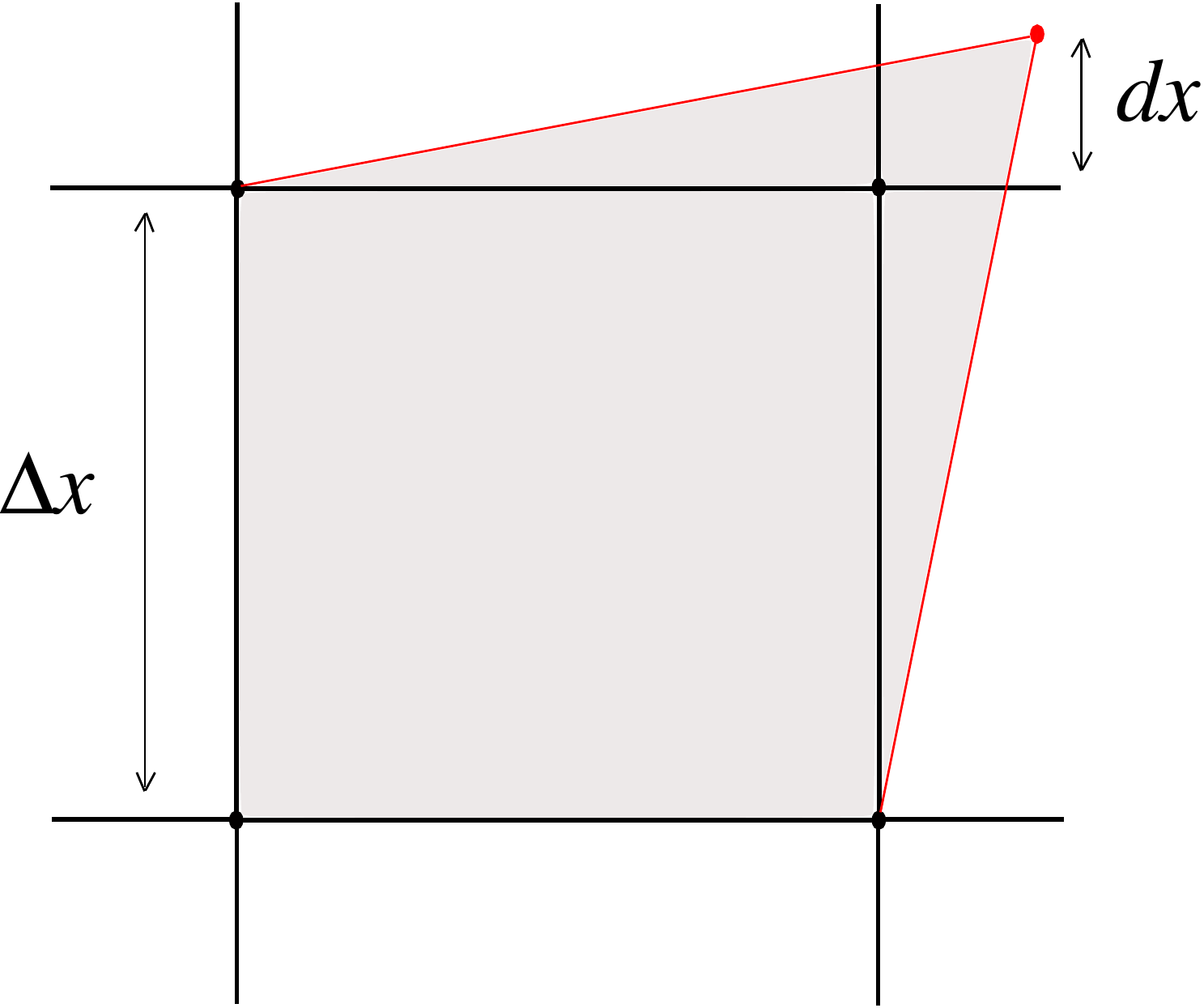}\hspace{0.15\textwidth}
\includegraphics[width=0.32\textwidth, height=0.26\textwidth]{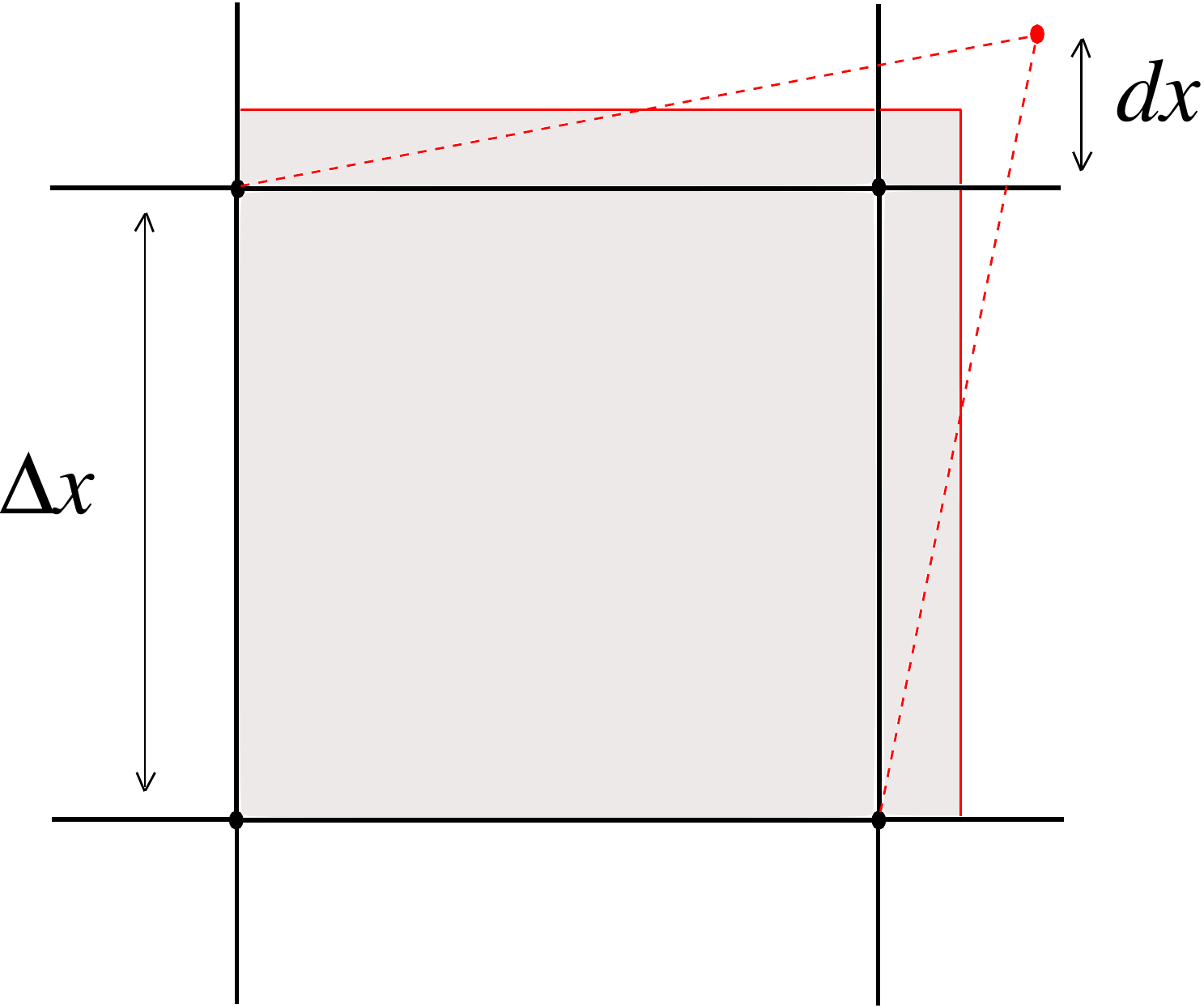}
\end{center}
\begin{center}
\includegraphics[width=0.32\textwidth, height=0.26\textwidth]{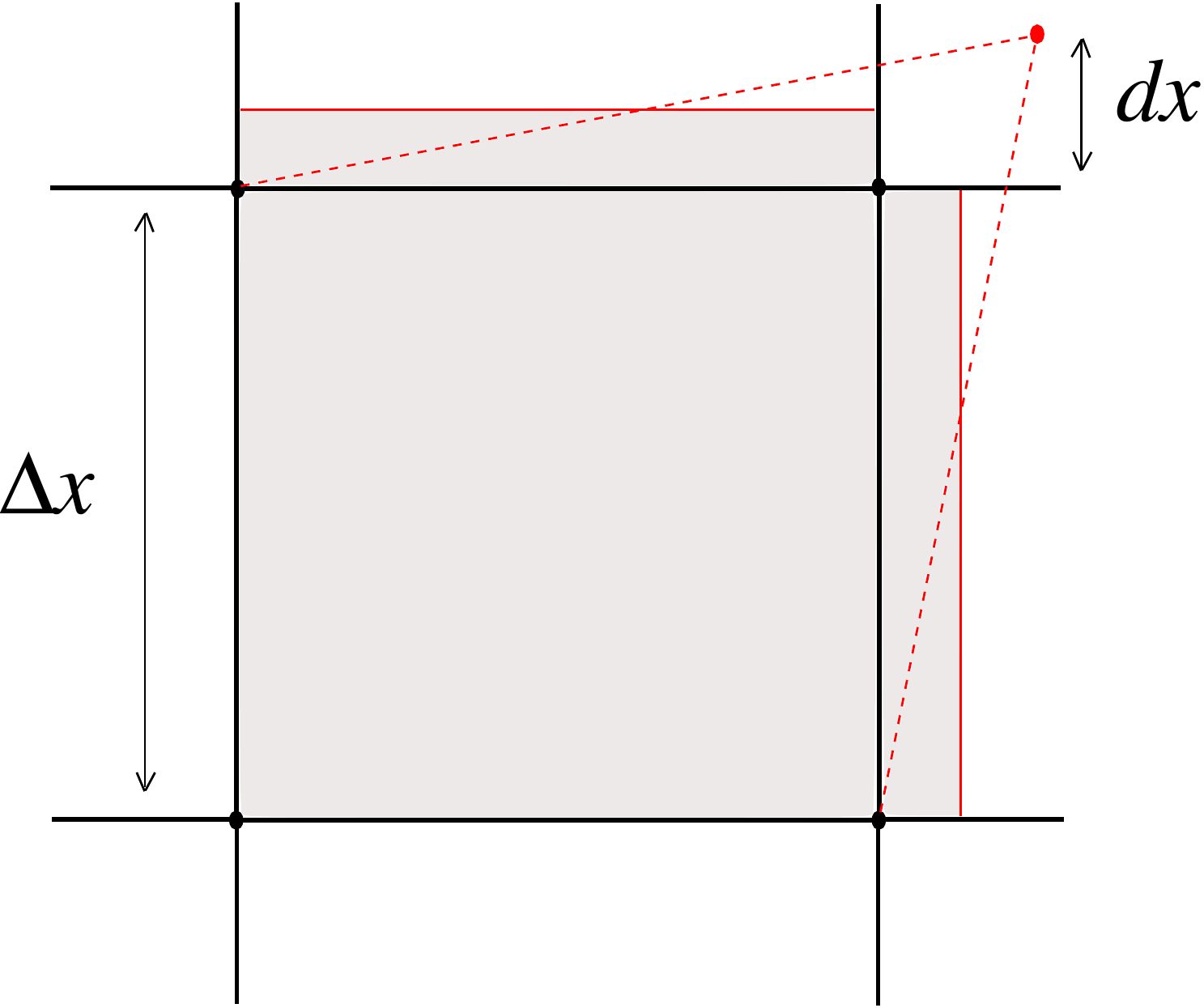}\hspace{0.15\textwidth}
\includegraphics[width=0.32\textwidth, height=0.26\textwidth]{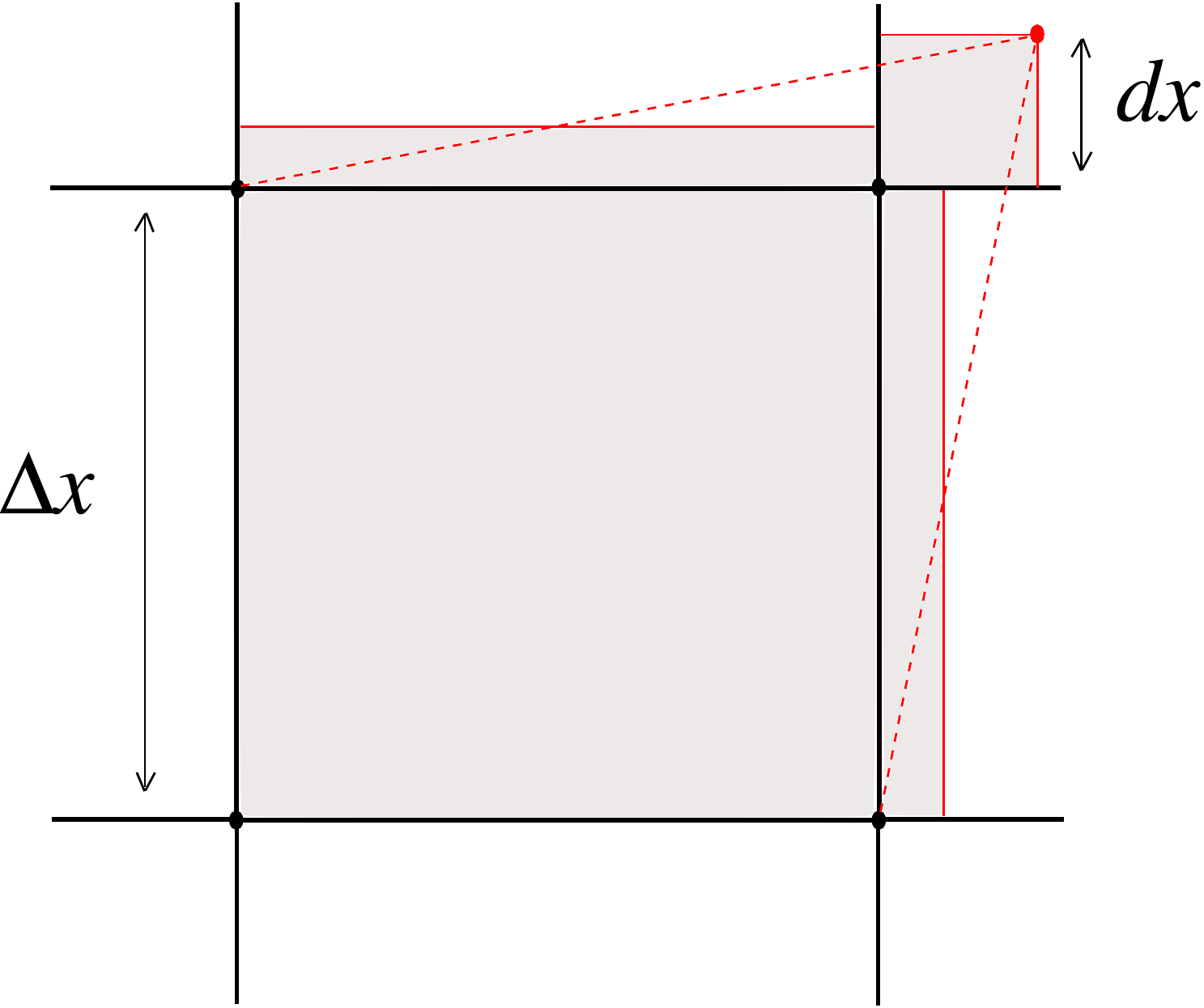}
\end{center}
\caption{Representation of the Lagrangian volume in the exact (top-left), AD (top-right), Direct (bottom-left) and DirectCF (bottom-right) cases.}
\label{fig:vol_lag_comparison}
\end{figure}

To be more specific, the quotient $Vol^{lag} / (\Delta x)^2$ is computed in function of $\varepsilon = dx / \Delta x$. Geometric representations of the different Lagrangian volumes are compared in Figure~\ref{fig:vol_lag_comparison}.

\noindent The results of the computations give, stopping at the third order in $\varepsilon$ for the Direct remap with Corner Fluxes:
\begin{itemize}
    \item[] Exact: $1 + \varepsilon$ 
    \item[] AD: $1 + \varepsilon + \varepsilon^2 / 4$
    \item[] Direct: $1 + \varepsilon$ 
    \item[] Direct with Corner Fluxes: $1 + \varepsilon  + \varepsilon^3  + o(\varepsilon^3)$ 
\end{itemize}

\textit{Remark 1:} In the AD case, since the effective Lagrangian volume $Vol^{lag, \:AD}$ does not appear directly, the value is deduced from $\rho^{projy}$, defining this volume such that $\rho^{projy} Vol^{lag, \:AD} = m^n$. Besides, in this simple case, this representation is strictly equivalent to a rectangular approximation, which consists in taking the middles of Lagrangian faces and taking the rectangle passing through these points. \\
Let us be a little more specific, starting with denoting $\rho^n = m^n / \Delta x^2$. It is clear that:
\begin{displaymath}
Vol^{lagx} = \Delta x^2 + \frac{dx}{2} \Delta x, \:\:\:\: \mbox{thus} \:\:\: m^{projx} = \frac{m^n}{1 + \varepsilon / 2}
\end{displaymath}
Then, since $Vol^{lagy} = Vol^{lagx}$, one gets:
\begin{displaymath}
m^{projy} = \frac{m^{projx}}{1 + \varepsilon / 2} = \frac{m^n}{(1 + \varepsilon / 2)^2}
\end{displaymath}
And finally, according to the definition of $Vol^{lag, \:AD}$ given below:
\begin{displaymath}
\frac{Vol^{lag, \:AD}}{\Delta x^2} = \frac{m^n}{m^{projy}} = (1 + \varepsilon / 2)^2 = 1 + \varepsilon + \varepsilon^2 / 4
\end{displaymath}

\textit{Remark 2:} The Direct gives exactly the right Lagrangian volume. And the DirectCF is one order more accurate than the AD in this particular case.

In order to estimate corner effects, let us determine the mass flux passing through the top-right corner during the remap phase. Mass fluxes are coloured in green in the Figures.

\vspace{0.5cm}	
\begin{figure}[h]
\begin{center}
\includegraphics[width=0.32\textwidth, height=0.26\textwidth]{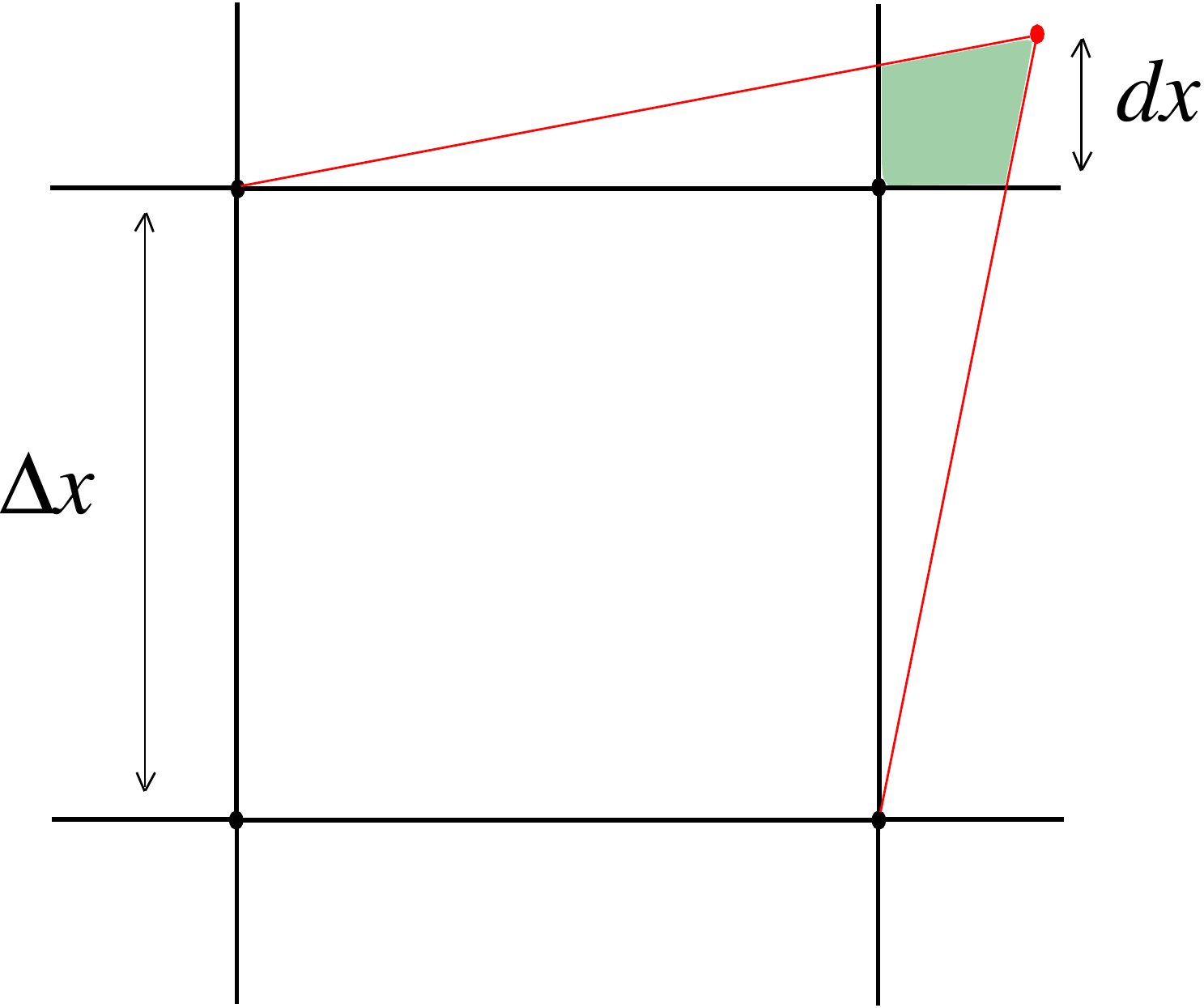}
\includegraphics[width=0.32\textwidth, height=0.26\textwidth]{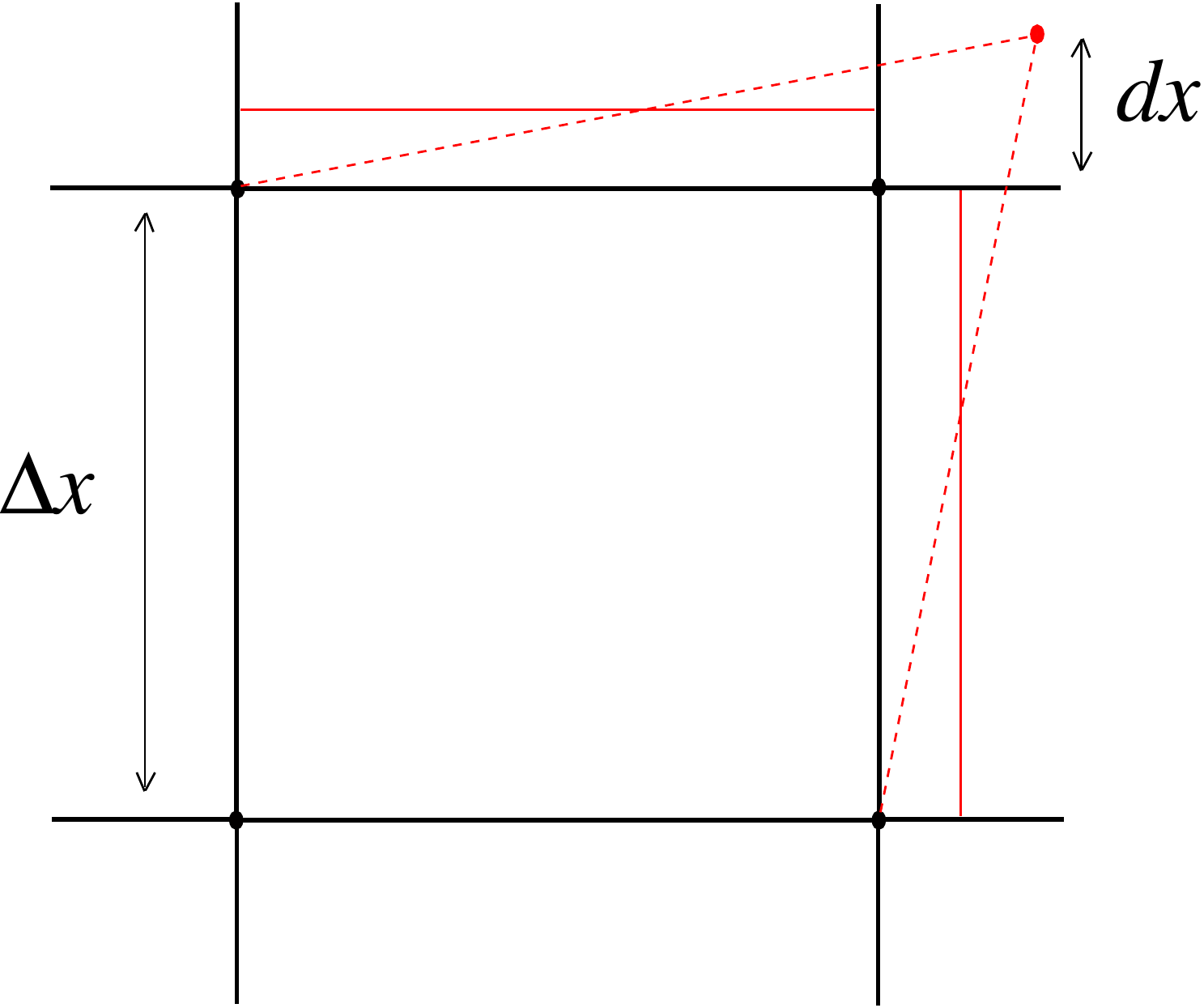}
\includegraphics[width=0.32\textwidth, height=0.26\textwidth]{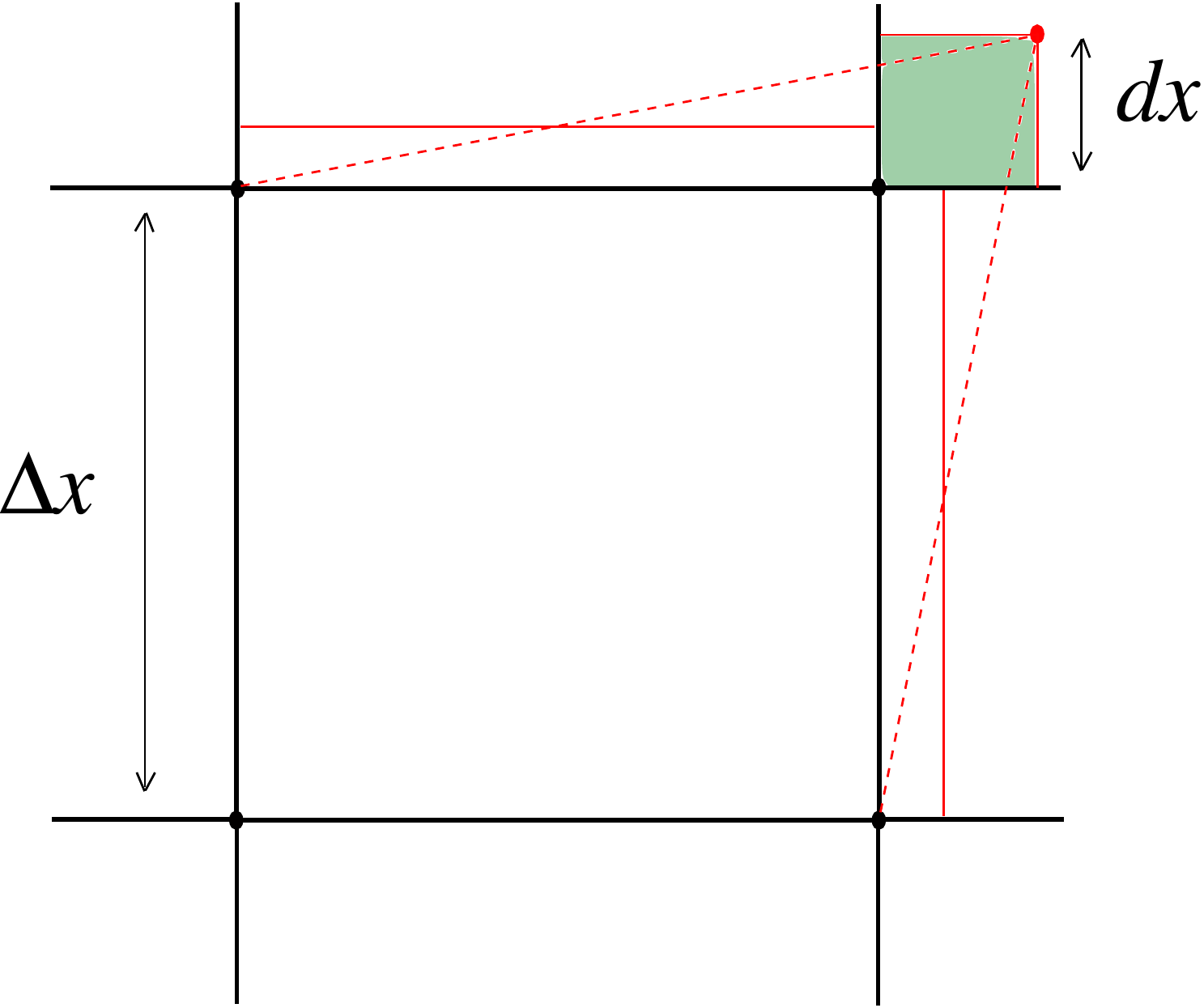}
\caption{Mass corner flux in the  exact (left), Direct (centre) and DirectCF (right) cases.}
\label{fig:mass_cf_direct_directcf}
\end{center}
\end{figure}

\noindent In the Direct and DirectCF cases, this mass corner flux is easily computed knowing $\rho^{lag}$ and $dVol_p$, see Figure~\ref{fig:mass_cf_direct_directcf}. 

\begin{figure}[h]
\begin{center}
\includegraphics[width=0.49\textwidth, height=0.26\textwidth]{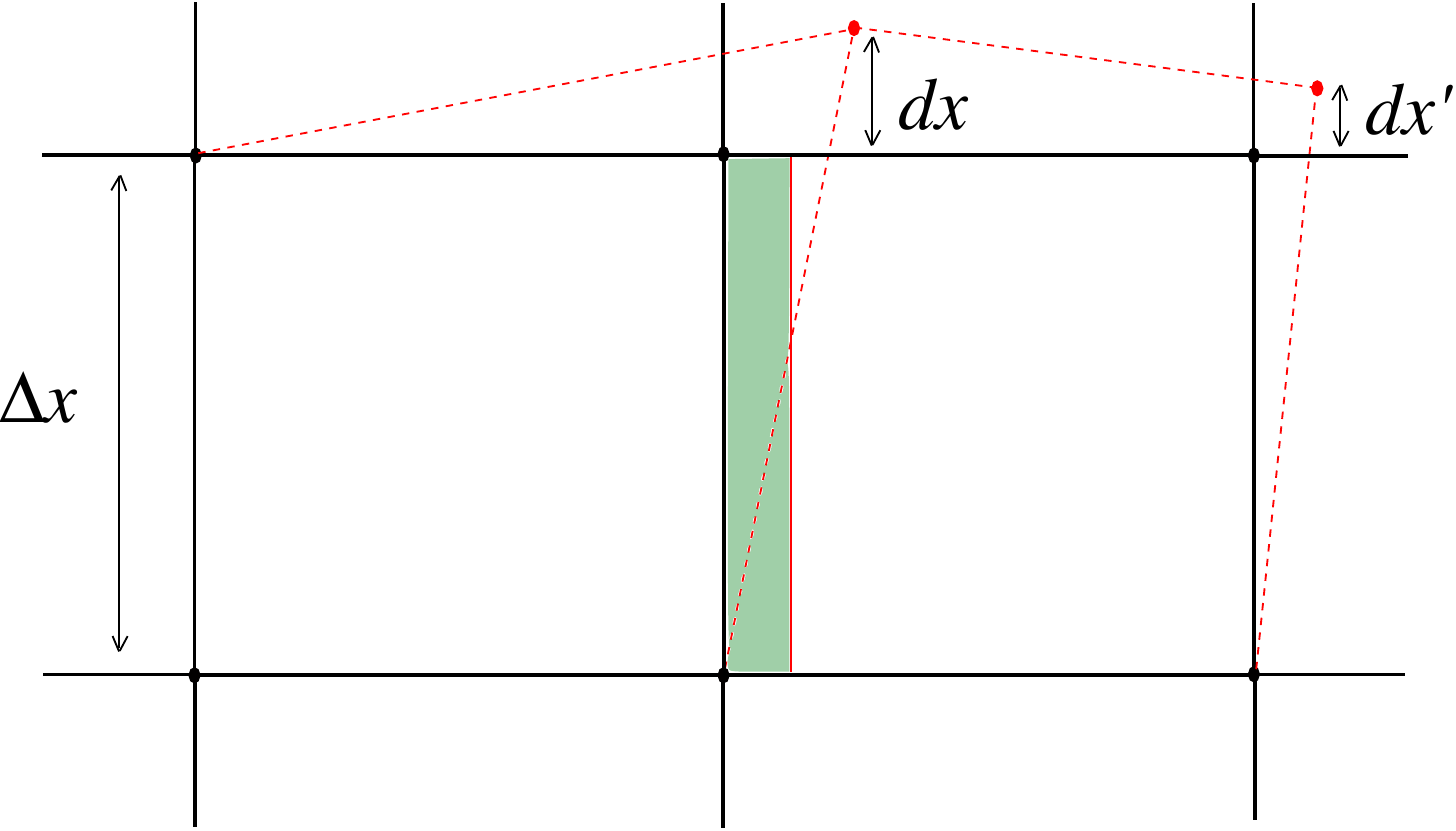}
\includegraphics[width=0.49\textwidth, height=0.26\textwidth]{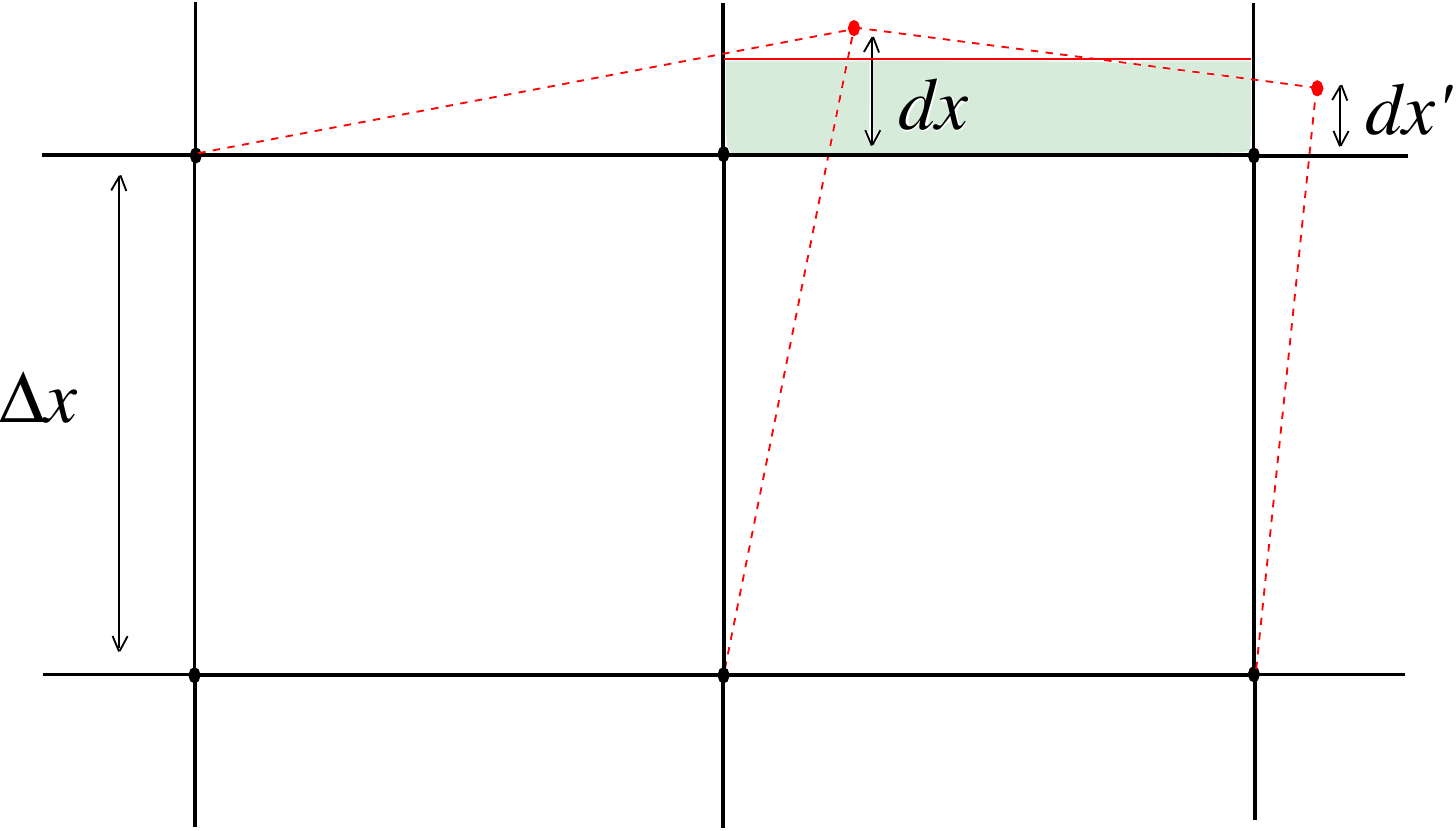}
\end{center}
\caption{Mass face fluxes at $X$ (top) and $Y$ (bottom) remap phases in the AD case.}
\label{fig:mass_cf_adi}
\end{figure}

\noindent But in the AD case, the computation is less straight-forward because of the directional splitting, see Figure~\ref{fig:mass_cf_adi}. Indeed, the result depends on the displacement of the node situated at the right of the top-right node considered, which does not seem really natural. This displacement is denoted $dx'$.

\vspace{0.7cm}	

Similarly to the representation of the volume, the quotient $dm_{p} / m$ is expressed in function of $\varepsilon = dx / \Delta x$ and $\varepsilon' = dx' / \Delta x$ (for the AD remap), which gives at the third order
\begin{itemize}
    \item[] Exact: $\varepsilon^2 - 2 \varepsilon^3$
    \item[] AD: $\varepsilon ^2/ 4 + \varepsilon \varepsilon'/ 4 -
    \frac{1}{4} \big( \varepsilon^3 + \varepsilon \varepsilon'^2 + \varepsilon' \varepsilon^2 \big)$ 
    \item[] Direct: $0$ 
    \item[] Direct with Corner Fluxes: $\varepsilon^2 - \varepsilon^3$
\end{itemize}
    
\textit{Remark 1:} The results for the Direct and AD remaps are not close to the exact value, they are not even accurate at the highest order (second) in this simple case. The fact that the result in the AD case depends on $\varepsilon'$ shows a drawback of this remap, at least formally.  

\textit{Remark 2:} In this example, the DirectCF gives results for Lagrangian volume and mass corner flux that coincide with exact values at the highest order. Thus, it seems that the DirectCF offers a better geometric lagrangian representation than the one in the AD remap.

	\subsection{Impact of the corners on the scheme}

\hspace{1.5em}The first notable impact of adding corners fluxes in the Lagrangian representation is the mandatory addition of corner fluxes on the dual mesh, as seen in section 2.3.

Another consequence, also mentioned in section 2.3, is the definition of reconstructed values of variables at the corners. As for the reconstruction at faces, the order of interpolation is the user's choice, which has to be consistent with the way face values are reconstructed. Contrary to the first order reconstruction that is obvious, it seems that second order reconstruction can be computed in different ways. Three of them are presented in the following paragraphs.

\paragraph{First order: upwind reconstruction}
$\:$

Let $a^{lag}$ be the cell-centred lagrangian variable to reconstruct at the corner. This exponent $^{lag}$ is dropped to simplify notations: $a^{lag} \leftrightarrow a$. At the left-bottom corner $p=(i-1/2,j-1/2)$, the upwind reconstruction writes:
\begin{equation*}
a_p^{upw} = a_c \:\:\: \mbox{where} \:\:\: c = \left\lbrace
\begin{array}{l r}
(i-1,j-1) & \mbox{if} \:\:\: dx^{n+1/2}_p > 0, \:\:dy^{n+1/2}_p  > 0 \\
(i,j) & \mbox{if} \:\:\: dx^{n+1/2}_p < 0, \:\:dy^{n+1/2}_p < 0 \\
(i-1,j) & \mbox{if} \:\:\: dx^{n+1/2}_p > 0, \:\:dy^{n+1/2}_p  < 0 \\
(i,j-1) & \mbox{if} \:\:\: dx^{n+1/2}_p < 0, \:\:dy^{n+1/2}_p  > 0
\end{array}
\right.
\end{equation*}

\paragraph{Second order: average reconstruction in the diagonal direction}
$\:$

If corner fluxes are taken into account exactly as face fluxes are, the value of each variable at corner $p=(i-1/2,j-1/2)$ has to be reconstructed along the diagonal direction $x=y$. The most basic reconstruction writes:
\begin{equation*}
\widetilde{a_{p}}^{o2, \:avg} = \frac{1}{2} \: \left\lbrace
\begin{array}{l r}
(a_{i-1,j-1} + a_{i,j}) & \mbox{if} \:\:\: dx^{n+1/2}_p \: dy^{n+1/2}_p  > 0 \\
(a_{i-1,j} + a_{i,j-1}) & \mbox{if} \:\:\: dx^{n+1/2}_p \: dy^{n+1/2}_p  < 0 
\end{array}
\right.
\end{equation*}

\noindent Actually, in order to prevent the creation of negative masses, which can happen if $\rho_{p}^{o2, \:avg} > \rho_c$, with $c$ the donor cell defined in the upwind reconstruction, the following formula is preferred:
\begin{equation*}
a_{p}^{o2, \:avg} = \min\left( a_c, \widetilde{a_{p}}^{o2, \:avg} \right)
\end{equation*}

\paragraph{Second order: linear limited reconstruction in both $X$ and $Y$ directions}
$\:$

Another way of reconstructing corner values comes from the he 2D grid point of view. Second order expressions for face values rely on a writing related to coordinates. For example, at $X$-face $f$, denoting by $\mathbf{x}_f$ the middle of this face and $\mathbf{x}_c$ the centroid of the donor cell, the 1D $X$-reconstruction is based on the following feature:
\begin{equation*}
\mathbf{x}_f = \mathbf{x}_c \pm \frac{\Delta x}{2} \: \mathbf{e}_x
\end{equation*}
Thus, noting that at corner $p=(i-1/2,j-1/2)$, one can write
\begin{equation*}
\mathbf{x}_p = \mathbf{x}_c \pm \frac{\Delta x}{2} \: \mathbf{e}_x \pm \frac{\Delta y}{2} \: \mathbf{e}_y
\end{equation*}
It seems fair to choose at node $p$ a reconstruction based on the one introduced in section 2.1, that is:
\begin{equation*}
a_{p}^{o2, \:xy} = a_c + \frac{{\delta}^{o2}_x a_c}{2}\:\big(\:\mbox{sgn}(dVol_{f_x}) \: {\Delta x}_c^{lag} - \Delta t^n \: \mathbf{u}_{f_x}^{n+1/2}\cdot\mathbf{e}_x\: \big)
\end{equation*}
\begin{equation*}
\:\:\:\:\:\:\:\:\: + \:\frac{{\delta}^{o2}_y a_c}{2}\:\big(\:\mbox{sgn}(dVol_{f_y}) \: {\Delta y}_c^{lag} - \Delta t^n \: \mathbf{u}_{f_y}^{n+1/2}\cdot\mathbf{e}_y\: \big)
\end{equation*}

\begin{equation*}
\mbox{where} \:\:\: f_x = \left\lbrace
\begin{array}{l r}
(i-1/2,j-1) & \mbox{if} \:\:\: dy^{n+1/2}_p  > 0 \\
(i-1/2,j) & \mbox{if} \:\:\: dy^{n+1/2}_p  < 0 
\end{array}
\right.
\end{equation*}
\begin{equation*}
\:\mbox{and} \:\:\:\:\:\: f_y = \left\lbrace
\begin{array}{l r}
(i-1,j-1/2) & \mbox{if} \:\:\: dx^{n+1/2}_p  > 0 \\
(i,j-1/2) & \mbox{if} \:\:\: dx^{n+1/2}_p  < 0 
\end{array}
\right.
\end{equation*}

\textit{Remark:} This formula is quite convenient since it can be expressed as a simple sum or difference of values of $a$ at centres and faces of the cells, which are already known.

\paragraph{Second order: linear limited reconstruction in the diagonal direction}
$\:$

Let us denote by $\mathbf{v}$ and $\mathbf{v}_\bot$ the unit vectors along the diagonal and antidiagonal directions, such that:
\begin{equation*}
\mathbf{v} = \frac{1}{\sqrt{{\Delta x}^2 + {\Delta y} ^2}} \left(\! \begin{array}{clcr} \Delta x \\ \Delta y \end{array} \! \right), \:\:\:\:\: \mathbf{v}_\bot = \frac{1}{\sqrt{{\Delta x}^2 + {\Delta y} ^2}} \left(\! \begin{array}{clcr} \Delta x \\ -\Delta y \end{array} \! \right)
\end{equation*}

\noindent Note that $\mathbf{v} \cdot \mathbf{v}_\bot = 0$ only if $\Delta x = \Delta y$. The approach is the same as for the $x$-$y$ reconstruction, but using the coordinate system $(\mathbf{v}, \mathbf{v}_\bot)$ instead of $(\mathbf{e}_x, \mathbf{e}_y)$. In other words, at corner $p$, let us write:
\begin{equation*}
\:\:\:\:\:\:\mathbf{x}_p = \mathbf{x}_c \pm \frac{1}{2} \sqrt{{\Delta x}^2 + {\Delta y} ^2} \: \mathbf{v}
\end{equation*}
\begin{equation*}
\mbox{or,} \:\:\: \mathbf{x}_p = \mathbf{x}_c \pm \frac{1}{2} \sqrt{{\Delta x}^2 + {\Delta y} ^2} \: \mathbf{v}_\bot
\end{equation*}

\noindent For a quantity $a_c$ in cell $c$, one is able to compute discrete limited gradients ${\delta}^{o2}_\mathbf{v} a_c$ and ${\delta}^{o2}_{\mathbf{v}_\bot} a_c$ along the diagonal and antidiagonal directions, with respect to what is done along the $X$ and $Y$ directions in section 2.1.1. In cell $c = (i,j)$, this leads:

\begin{equation*}
{\delta}^{o2}_\mathbf{v} a_{i,j} = \frac{1}{2} \: \big( \: \Phi_{i,j}^+\:{\delta}_\mathbf{v} a_{i+1/2,j+1/2} + \Phi_{i,j}^- \:{\delta}_\mathbf{v} a_{i-1/2,j-1/2} \: \big)
\end{equation*}

\begin{equation*}
\mbox{with} \:\:\:\: {\delta}_\mathbf{v} a_{i\pm1/2,j\pm1/2} = \frac{a_{i\pm1,j\pm1} - a_{i,j}}{(\mathbf{x}_{i\pm1,j\pm1}^{lag} - \mathbf{x}_{i,j}^{lag}) \cdot \mathbf{v}}
\end{equation*}

\begin{equation*}
\mbox{and} \:\:\:\:\: \Phi_{i,j}^+ = \varphi (r)\:, \:\:\: \Phi_{i,j}^- = \varphi (1/r)\:, \:\:\: r = \frac{{\delta}_\mathbf{v} a_{i-1/2,j-1/2}}{{\delta}_\mathbf{v} a_{i+1/2,j+1/2}}
\end{equation*}

\noindent The expression of ${\delta}^{o2}_{\mathbf{v}_\bot} a_{i,j}$ is similar. Given those discrete gradients, one can write a second order reconstruction as follows:

\begin{equation*}
a_{p}^{o2, \:diag} = \left\lbrace 
\begin{array}{l}
a_c + \dfrac{{\delta}^{o2}_\mathbf{v} a_c}{2}\:\big(\:\mbox{sgn}(dx_p^{n+1/2}) \: {\Delta \mathbf{x}}_c^{lag} - \Delta t^n \: \mathbf{u}_{p}^{n+1/2} \: \big)\cdot\mathbf{v} \\
\: \\
a_c + \dfrac{{\delta}^{o2}_{\mathbf{v}_\bot} a_c}{2}\:\big(\:\mbox{sgn}(dx_p^{n+1/2}) \: {\Delta \mathbf{x}}_c^{lag} - \Delta t^n \: \mathbf{u}_{p}^{n+1/2} \: \big)\cdot\mathbf{v}_\bot
\end{array}
\right.
\end{equation*}

\noindent The reconstruction is performed along $\mathbf{v}$ if $dx^{n+1/2}_p \: dy^{n+1/2}_p > 0$ and along $\mathbf{v}_\bot$ if $dx^{n+1/2}_p \: dy^{n+1/2}_p < 0$. 

\textit{Remark:} Signs of quantities in the formula depend on the sign of $dx^{n+1/2}_p$ since $\mathbf{v}$ and $\mathbf{v}_\bot$ have been chosen in the $x>0$ direction.

\paragraph{Second order: multidimensional reconstruction}
$\:$

Instead of computing 1D linear reconstructions along all directions, the idea here is to build a bidimensional patch on the 9-points stencil. In order to be consistent with what is done for the second order in the 1D case, this patch has to be a plane. According to the multidimensional reconstruction model presented in \cite{diot2012methode}, with respect to the notations previously introduced, the analytic expression of the parametrization $\check{a}_{i,j}$ of the plane reconstructed at cell $c = (i,j)$ for the quantity $a$ writes:
\begin{displaymath}
\check{a}_{i,j}(\mathbf{x}) = a_{i,j} + \mathbf{G}_{i,j} \cdot (\mathbf{x} - \mathbf{x}_c)
\end{displaymath} 

\noindent Such a writing ensures that
\begin{displaymath}
\frac{1}{\lvert \Omega_{i,j} \rvert} \: \int_{\Omega_{i,j}} \!\!\check{a}_{i,j}(\mathbf{x})\: \mathrm{d}\mathbf{x} = a_{i,j}
\end{displaymath}

The difficulty is to find the components of $\mathbf{G}_{i,j} = \big(G_{i,j}^1,G_{i,j}^2\big)$ that will provide a reconstruction as close to the real solution as possible. This will be done solving a least squares problem. The criterion chosen is the mean value of $\check{a}_{i,j}$ in each cell of the 9-points stencil. Consequently, the functional $\mathrm{E}_{i,j}$ to minimize is the following:

\begin{displaymath}
\mathrm{E}_{i,j} = \frac{1}{2} \sum\limits_{\substack{\lvert k - i \rvert \leq 1 \\ \lvert l - j \rvert \leq 1}} \left( \frac{1}{\lvert \Omega_{k,l} \rvert} \: \int_{\Omega_{k,l}} \!\!\check{a}_{i,j}(\mathbf{x})\: \mathrm{d}\mathbf{x} - a_{k,l} \right)^2
\end{displaymath}

\noindent The minimum is reached if, for all $(k,l)$

\begin{displaymath}
\frac{1}{\lvert \Omega_{k,l} \rvert} \: \int_{\Omega_{k,l}} \!\!\check{a}_{i,j}(\mathbf{x})\: \mathrm{d}\mathbf{x} = a_{k,l}
\end{displaymath}

\noindent Let us denote:
\begin{displaymath}
\left\lbrace \begin{array}{c}
X_{i,j,k,l}^1 = \dfrac{1}{\lvert \Omega_{k,l} \rvert} \: \displaystyle\int_{\Omega_{k,l}} \!\!(x - x_c)\: \mathrm{d}\mathbf{x} \\
\: \\
X_{i,j,k,l}^2 = \dfrac{1}{\lvert \Omega_{k,l} \rvert} \: \displaystyle\int_{\Omega_{k,l}} \!\!(y - y_c)\: \mathrm{d}\mathbf{x}
\end{array}
\right.
\end{displaymath}

\noindent This leads to the following system in the matrix form for the minimum:

\begin{displaymath}
\begin{pmatrix}
X_{i,j,i-1,j-1}^1 & X_{i,j,i-1,j-1}^2 \\
\vdots & \vdots \\
\vdots & \vdots \\
X_{i,j,i-1,j}^1 & X_{i,j,i-1,j}^2
\end{pmatrix} \cdot
\begin{pmatrix}
G_{i,j}^1 \\
\: \\
G_{i,j}^2
\end{pmatrix} =
\begin{pmatrix}
a_{i-1,j-1} - a_{i,j} \\
\vdots \\
\vdots \\
a_{i-1,j} - a_{i,j}
\end{pmatrix}
\end{displaymath}
\vspace{0.2cm}

The overdetermined system to be solved rewrites: $\mathbf{X}_{i,j} \: \mathbf{G}_{i,j} = \mathbf{A}_{i,j}$, where geometric ($\mathbf{X}_{i,j}$) and physical ($\mathbf{A}_{i,j}$) quantities are dissociated. If the rank of the matrix $\mathbf{X}_{i,j}$ is maximal, then the symmetric matrix $^t\mathbf{X}_{i,j} \: \mathbf{X}_{i,j}$ is invertible and the solution of the least squares problem writes:
\begin{displaymath}
\mathbf{G}_{i,j} = (^t\mathbf{X}_{i,j} \: \mathbf{X}_{i,j})^{-1} \: ^t\mathbf{X}_{i,j} \mathbf{A}_{i,j}
\end{displaymath}

\noindent The matrix $(^t\mathbf{X}_{i,j} \: \mathbf{X}_{i,j})^{-1} \: ^t\mathbf{X}_{i,j}$ is called the Moore-Penrose pseudo-inverse of $\mathbf{X}_{i,j}$, see \cite{gw1998matrix}. In this case, note that $^t\mathbf{X}_{i,j} \: \mathbf{X}_{i,j}$ is a $2\!\times\!2$ matrix, thus its inversion is not very costly in terms of calculation time. It can be done directly, using the commatrix formula, or in two steps: first performing a $QR$ decomposition $^t\mathbf{X}_{i,j} \: \mathbf{X}_{i,j} = \mathbf{Q}_{i,j} \mathbf{R}_{i,j}$, then inverting the upper triangular matrix $\mathbf{R}_{i,j}$ and writing
\begin{displaymath}
\mathbf{G}_{i,j} = \mathbf{R}_{i,j}^{-1} \,^t\mathbf{Q}_{i,j}\, ^t\mathbf{X}_{i,j} \, \mathbf{A}_{i,j}
\end{displaymath}
In order to facilitate the possible extension to higher order reconstructions, the second method has been implemented.

The last formula gives the values of $G_{i,j}^1$ and $G_{i,j}^2$, which enable to get the expression of $\check{a}_{i,j}$, the non-limited reconstruction of $a_{i,j}$ on the stencil. Since the main objective is to remain robust, a limitation has to be performed. Formally, the limiter writes as a vector function $\mathbf{\Phi}_{i,j}$ directly applied to $\mathbf{G}_{i,j}$, such that the new limited reconstruction gives:
\begin{displaymath}
\check{a}_{i,j}(\mathbf{x}) = a_{i,j} + \mathbf{\Phi}_{i,j} \left(\mathbf{G}_{i,j}\right) \cdot (\mathbf{x} - \mathbf{x}_c)
\end{displaymath} 
The goal of the limitation is to diminish the norm of $\mathbf{G}_{i,j}$ when high gradients appear, but its direction needs to remain the same. Consequently, the form of the limiter function is the following:
\begin{displaymath}
\mathbf{\Phi}_{i,j} \left(\mathbf{G}_{i,j}\right) = \phi_{i,j}\: \mathbf{G}_{i,j}
\end{displaymath} 

\noindent Once again, giving priority to robustness, the choice made for the scalar $\phi_{i,j}$ is:
\begin{displaymath}
\phi_{i,j} = \frac{1}{\parallel\! \mathbf{G}_{i,j} \!\parallel} \,\min \left(\, \lvert\delta_x^{o2} a_{i,j}\rvert\,, \:\lvert\delta_y^{o2} a_{i,j}\rvert\,, \:\lvert\delta_{\mathbf{v}}^{o2} a_{i,j}\rvert\,, \:\lvert\delta_{\mathbf{v}_\bot}^{o2} a_{i,j}\rvert \,\right)
\end{displaymath}
\noindent All $\delta_\cdot^{o2} a_{i,j}$ are computed using Van Leer slope limitation, as presented before (see section 2.1.1 for $x$, $y$, and previous paragraph for $\mathbf{v}$, $\mathbf{v}_\bot$). Note that the limitation coefficient takes into account the variations of $a$ in all directions. Besides, the norm of multidimensional limited gradient $\mathbf{\Phi}_{i,j} \left(\mathbf{G}_{i,j}\right)$ equals the minimum of the norms of limited 1D gradients in all directions. This ensures that the scheme remains stable.  


Now the analytic expression of limited $\check{a}_{i,j}$ is known, the reconstructed value of $a$ at any corner or face is obtained by evaluating the function $\check{a}_{i,j}$ at this point. For example, at corner $p$, if the centroid of $dVol_p$ is denoted $\mathbf{x}_{c_p}$, one simply gets:
\begin{displaymath}
a_{p}^{o2, \:multid} = \check{a}_{i,j}(\mathbf{x}_{c_p}) = a_{i,j} + \mathbf{\Phi}_{i,j} \left(\mathbf{G}_{i,j}\right) \cdot (\mathbf{x}_{c_p} - \mathbf{x}_c)
\end{displaymath}  

\textit{Remark:} Note that in the multidimensional case, unlike all other second order reconstructions presented before, values at faces are also reconstructed using $\check{a}_{i,j}$. Thus, corners and faces are treated with the same limitation coefficient.

\begin{figure}[h]
\begin{center}
\includegraphics[width=0.22\textwidth, height=0.22\textwidth]{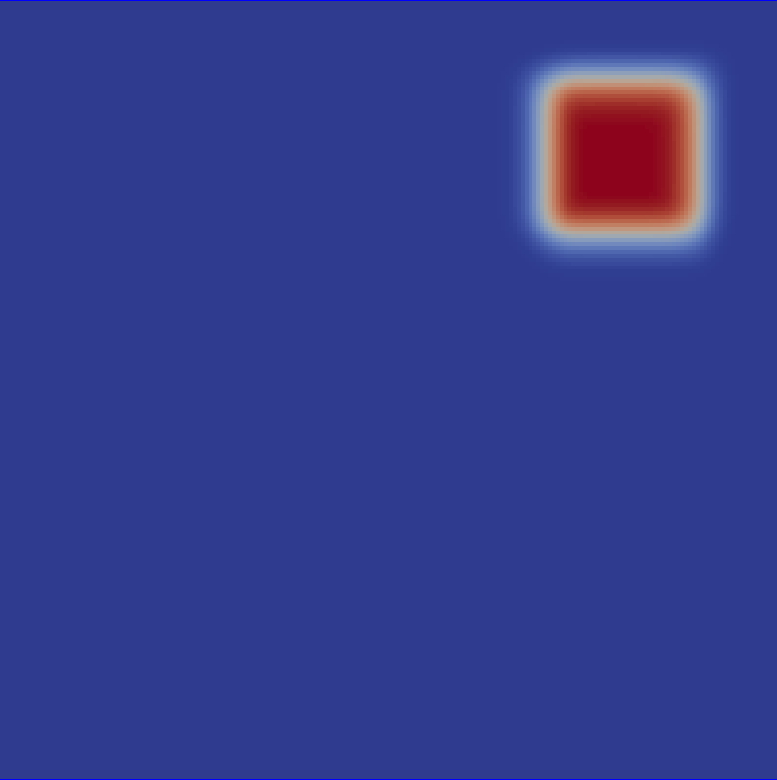}\hspace{0.03\textwidth}
\includegraphics[width=0.22\textwidth, height=0.22\textwidth]{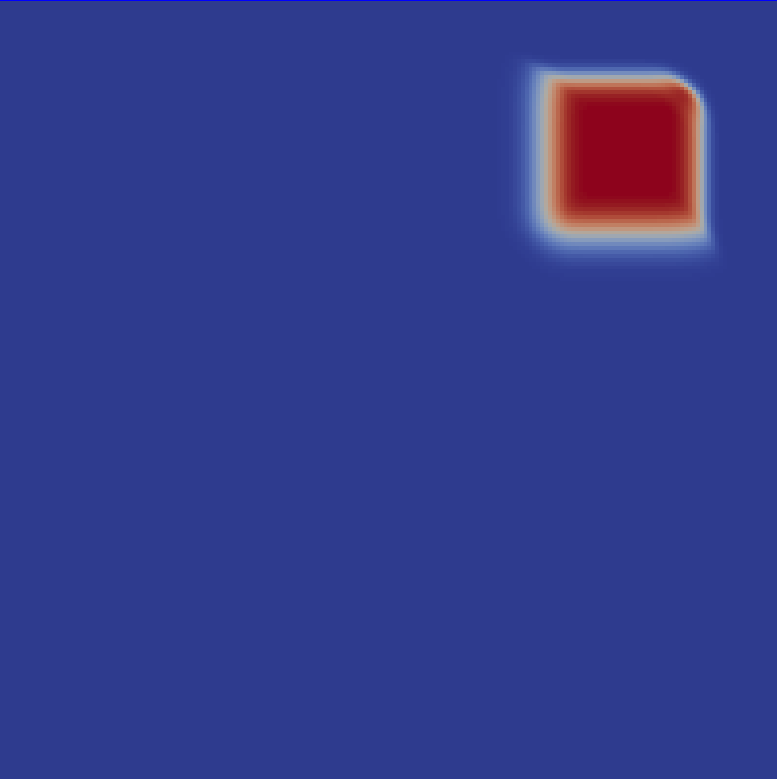}\hspace{0.03\textwidth}
\includegraphics[width=0.22\textwidth, height=0.22\textwidth]{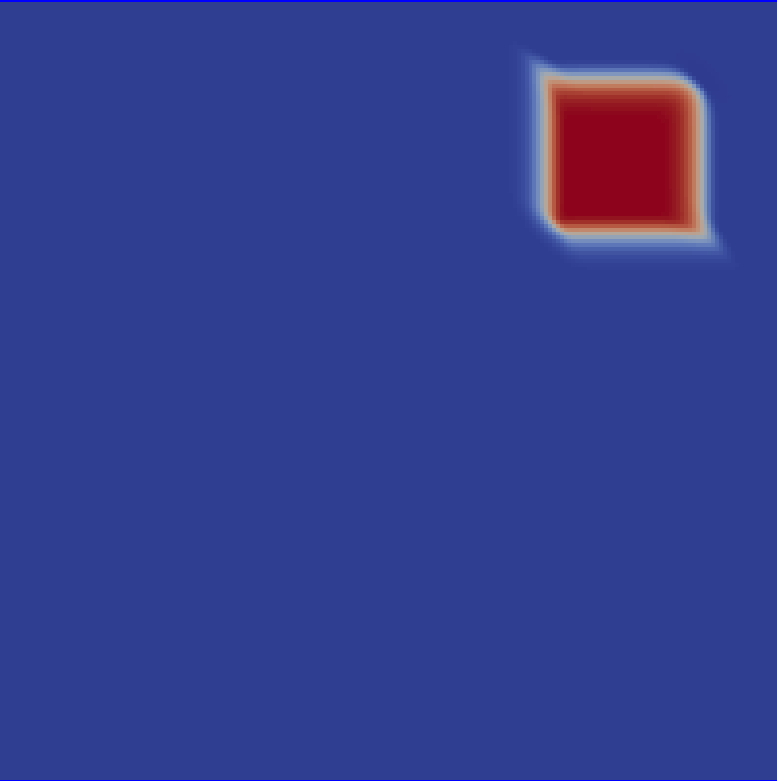}\hspace{0.03\textwidth}
\includegraphics[width=0.22\textwidth, height=0.22\textwidth]{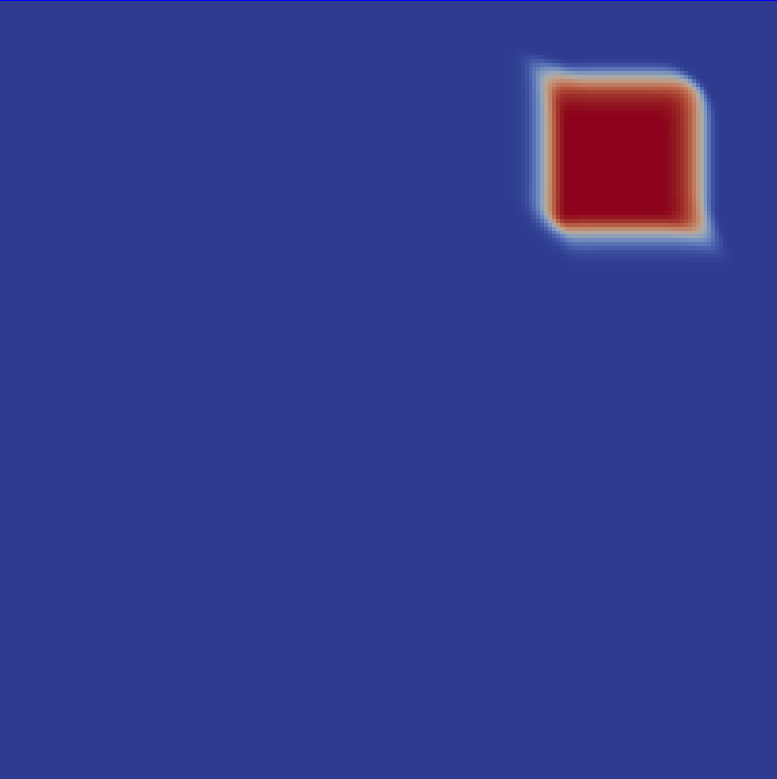}
\caption{Comparison between upwind, average, linear $x$-$y$ and linear diagonal reconstructions at the corners (from left to right).}
\label{fig:adv_monomat_cf}
\end{center}
\end{figure}

\vspace{1cm}
To emphasize the importance of the corners on the scheme, and on the results obtained, a simple test-case has been ran with each one of the four different reconstructions. Let us consider the linear advection of a square ($\rho = 10 \: \mbox{kg.m}^{-3}$) into a field ($\rho = 0.1 \: \mbox{kg.m}^{-3}$) with velocity $\mathbf{u} = (5\:; 5) \: \mbox{m.s}^{-1}$. The mesh contains $200\!\times\!200$ square cells. Fluxes at the faces are computed with a second order reconstruction. Results are gathered in Figure~\ref{fig:adv_monomat_cf}.

\vspace{0.5cm}

\textit{Remark 1:} As expected, the upwind reconstruction is the most diffusive, and the non-limited average reconstruction creates oscillations. In the upwind case, the shape of the square is preserved but the diffusion at the edges makes it not very accurate. 

\begin{figure}[h]
\begin{center}
\includegraphics[width=0.33\textwidth, height=0.28\textwidth]{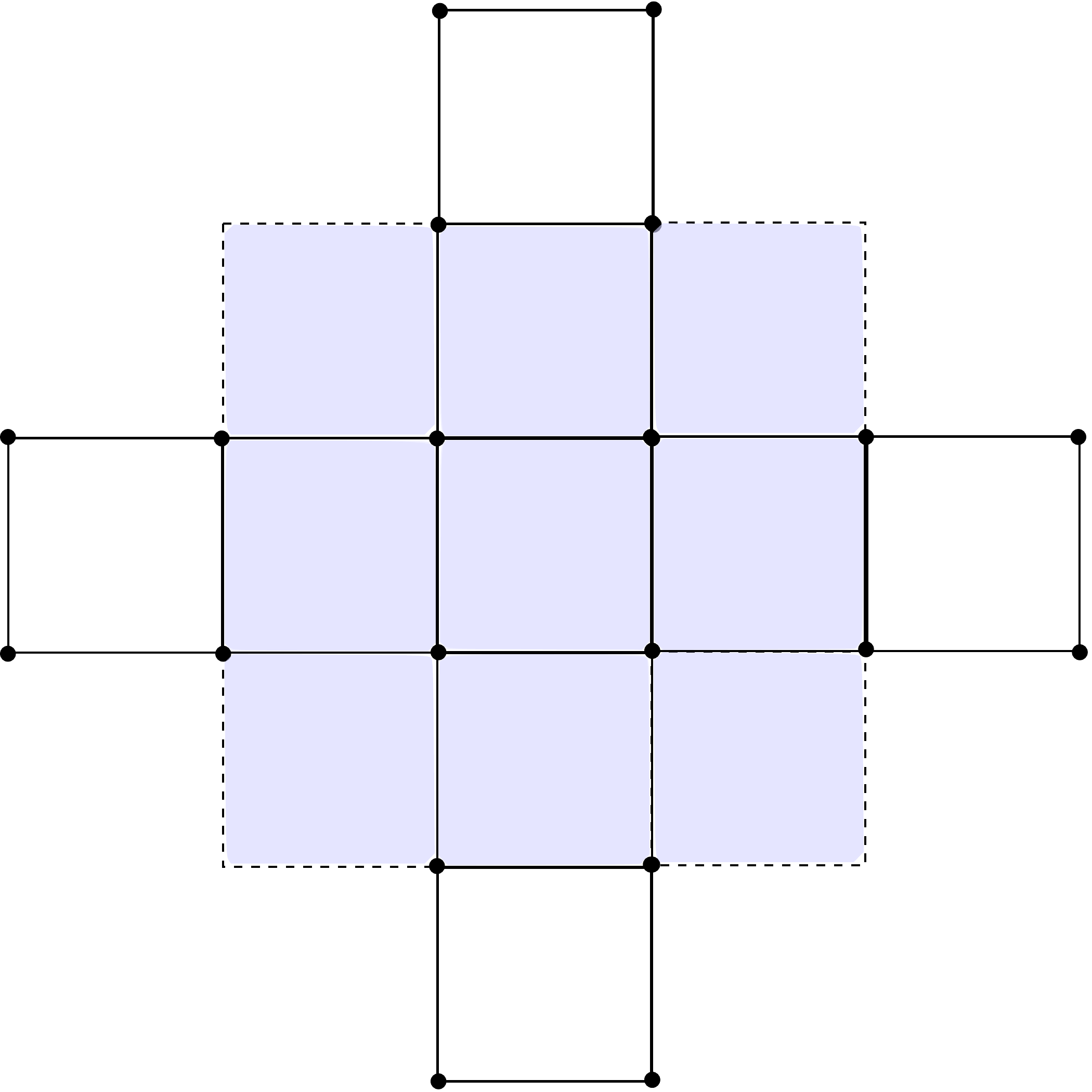}\hspace{2cm}
\includegraphics[width=0.33\textwidth, height=0.28\textwidth]{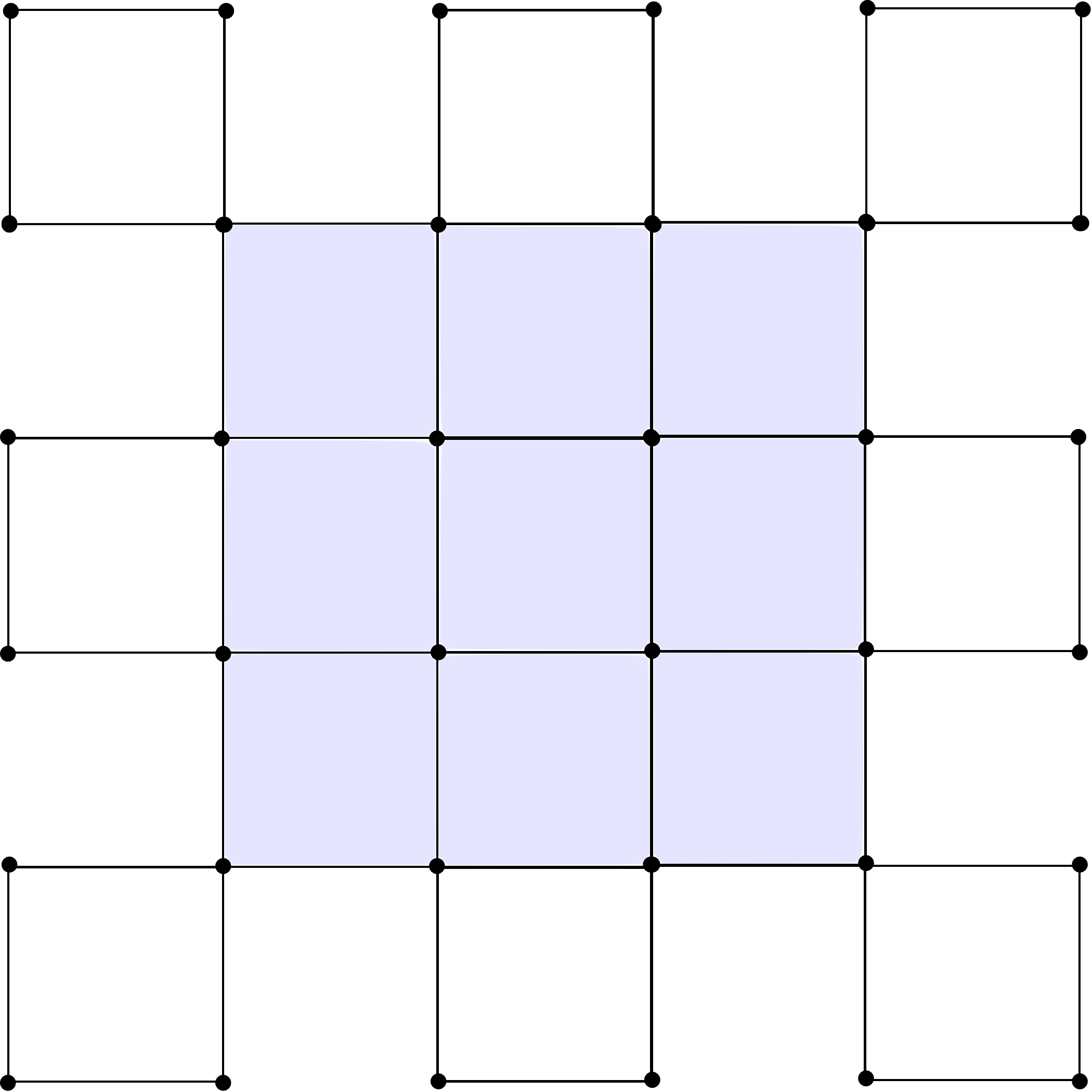}
\caption{Stencils of both $x$-$y$ (left) and diagonal (right) linear reconstructions.}
\label{fig:stencil_reconstruction}
\end{center}
\end{figure}

\textit{Remark 2:} The linear $x$-$y$ reconstruction tends to spread the corners along the $x=-y$ direction (i.e.\ the direction orthogonal to propagation), like with the Direct remap, see section 5. Here, this could be due to a wrong second order reconstruction at the corners. Indeed, the same reconstruction as the AD remap is used in each direction, and then the results are combined. This means that in each direction, the stencil used to reconstruct values has a width of 5 cells, which gives the 9 cells cross-shaped stencil $(i,[j-2, j+2]) \cup ([i-2, i+2],j)$. There is an inconsistency between this stencil (used for reconstruction) and the one of the remap, which is 9 cells square-shaped, coloured in blue on Figure~\ref{fig:stencil_reconstruction}. For example, it seems wrong to take into account the value in cell $(i+2,j)$ and not the one in cell $(i+1,j+1)$ for the reconstruction of a value at the corner $(i+1/2,j+1/2)$.

\textit{Remark 3:} The linear diagonal reconstruction seems more adapted since it shows the same default, but in a lesser extent. The stencil for this reconstruction is the union of the stencil for the $x$-$y$ reconstruction and the $X$-shaped stencil containing the two diagonals, which are 5 cells-long each. The total stencil is 17 cells-wide, see Figure~\ref{fig:stencil_reconstruction}, and it includes the 9 points stencil of the remap.

\vspace{0.5cm}

\textit{Remark 4:} The stencil of reconstruction in the multidimensional case is exactly the one of the remap since the patch is reconstructed over these 9 cells.

\textit{Remark 5:} Even if corner fluxes are very small volumes ($O({\Delta x}^2 {\Delta t}^2$)), their impact on the results is significant.

 \clearpage
 \newpage
 \section{Multi-material extension}

\hspace{1.5em}The numerical multi-material model displayed in this section is adapted to non-miscible fluids. It draws a sharp interface between two distinct materials. This method is called VOF (\textbf{V}olume \textbf{O}f \textbf{F}luids), see \cite{hirt1981volume}. Its implementation in SHY code will be limited to two fluids.

	\subsection{Lagrangian phase}

\hspace{1.5em}The multi-material Lagrangian phase is almost the same as in the mono-material case. The partial quantity of material $\alpha$ in cell $c$ is denoted $\alpha c$. \\
Let us define the volume fraction of material $\alpha$ in cell $c$ at $t^n$ as:
\begin{equation*}
	k_{\alpha c}^{n} = \frac{Vol_{\alpha c}^n}{Vol_{c}^n}
\end{equation*}
The hypothesis made in this model is called "iso-deformation", that is to say, all materials are assumed to have the same compressibility during Lagrangian phase. In other words, the following equality holds at each iteration $n$ 
\begin{equation*}
	k_{\alpha c}^{n+1,\:lag} = k_{\alpha c}^{n+1/2}  = k_{\alpha c}^{n} 
\end{equation*}
Such an assumption could seem quite rough, but it actually gives good results.

	\paragraph{Prediction phase}
	
\begin{equation}
\label{prediction_multi}
	\left\lbrace \hspace{0.05\textwidth}
	\begin{array}{l}
 		m_c^{n+1/2} = m_c^n\\
 		\: \\
 		m_{\alpha c}^{n+1/2} = m_{\alpha c}^n\\
 		\: \\
 		\mathbf{x}_p^{n+1/2} = \mathbf{x}_p^n + \dfrac{\Delta t^n}{2} \mathbf{u}_p^n\\
 		\: \\ 		
 		Vol_c^{n+1/2} = \mathcal{F}\big( (\mathbf{x}_p^{n+1/2})_{p \in \{c\}}\big)\\
 		\: \\
		Vol_{\alpha c}^{n+1/2} = Vol_{\alpha c}^{n} \dfrac{Vol_c^{n+1/2}}{Vol_c^{n}} = k_{\alpha c}^{n} Vol_c^{n+1/2} \\
		\: \\
 		\mathbf{u}_p^{n+1/2} = \mathbf{u}_p^{n} - \dfrac{\Delta t^n}{2} \big( \nabla (P + Q) \big)_p^{n+1/2} / \rho_p^n\\
 		\: \\
 		e_{\alpha c}^{n+1/2} = e_{\alpha c}^n - (P_{\alpha c}^n + Q_c^n) \left( \dfrac{1}{\rho_{\alpha c}^{n+1/2}} - \dfrac{1}{\rho_{\alpha c}^{n}} \right)\\
 		\: \\
		P_{\alpha c}^{n+1/2} = \mathcal{P} (\rho_{\alpha c}^{n+1/2},e_{\alpha c}^{n+1/2})
 	\end{array}
 	\right.
\end{equation}	

Moreover, the average pressure in a mixed-cell writes:
\begin{equation*}
	P_c^{n+1/2} = \sum_{\alpha c} \: k_{\alpha c}^{n+1/2} P_{\alpha c}^{n+1/2}
\end{equation*}	

 In the same way, total mass and internal specific energy are naturally defined as:
\begin{equation*}
	m_c^{n+1/2} = \sum_{\alpha c} \: m_{\alpha c}^{n+1/2}
\end{equation*}
\begin{equation*}
	e_c^{n+1/2} = \left( \sum_{\alpha c} \: m_{\alpha c}^{n+1/2} e_{\alpha c}^{n+1/2} \right) / m_c^{n+1/2} 
\end{equation*}

\textit{Remark:} The equation verified by the velocity is strictly the same as in the mono-material case. Indeed, the deformation of the cell is determined by the average pressure, and not partial pressures. Note that the definition of the pseudo-viscosity $Q$ given in section 1 only depends on $\mathbf{u}$, so there is only one $Q$ in each cell.

	\paragraph{Correction phase}
	
\begin{equation}
\label{correction_multi}
	\left\lbrace \hspace{0.05\textwidth}
	\begin{array}{l}
 		m_c^{n+1, \:lag} = m_c^n\\
 		\: \\
 		m_{\alpha c}^{n+1, \:lag} = m_{\alpha c}^n\\
 		\: \\
		\mathbf{x}_p^{n+1, \:lag} = \mathbf{x}_p^n + \Delta t^n \mathbf{u}_p^{n+1/2}\\
 		\: \\ 		
 		Vol_c^{n+1, \:lag} = \mathcal{F}\big( (\mathbf{x}_p^{n+1})_{p \in \{c\}}\big)\\
 		\: \\ 		
 		Vol_{\alpha c}^{n+1, \:lag} = k_{\alpha c}^{n} Vol_c^{n+1, \:lag} \\
		\: \\
 		\mathbf{u}_p^{n+1, \:lag} = 2 \: \mathbf{u}_p^{n+1/2} - \mathbf{u}_p^{n}\\
 		\: \\
 		e_{\alpha c}^{n+1, \:lag} = e_{\alpha c}^n - (P_{\alpha c}^{n+1/2} + Q_c^n) \left( \dfrac{1}{\rho_{\alpha c}^{n+1, \:lag}} - \dfrac{1}{\rho_{\alpha c}^{n}} \right)\\
 		\: \\
 		P_{\alpha c}^{n+1, \:lag} = \mathcal{P} (\rho_{\alpha c}^{n+1, \:lag},e_{\alpha c}^{n+1, \:lag})
 	\end{array}
 	\right.
\end{equation}

Thanks to system \eqref{correction_multi}, all partial and average variables can be recovered at $t^{n+1, \:lag}$.

	\paragraph{Interface positioning}

$\:$ \\

\begin{figure}[h]
\begin{center}
\includegraphics[width=0.45\textwidth, height=0.35\textwidth]{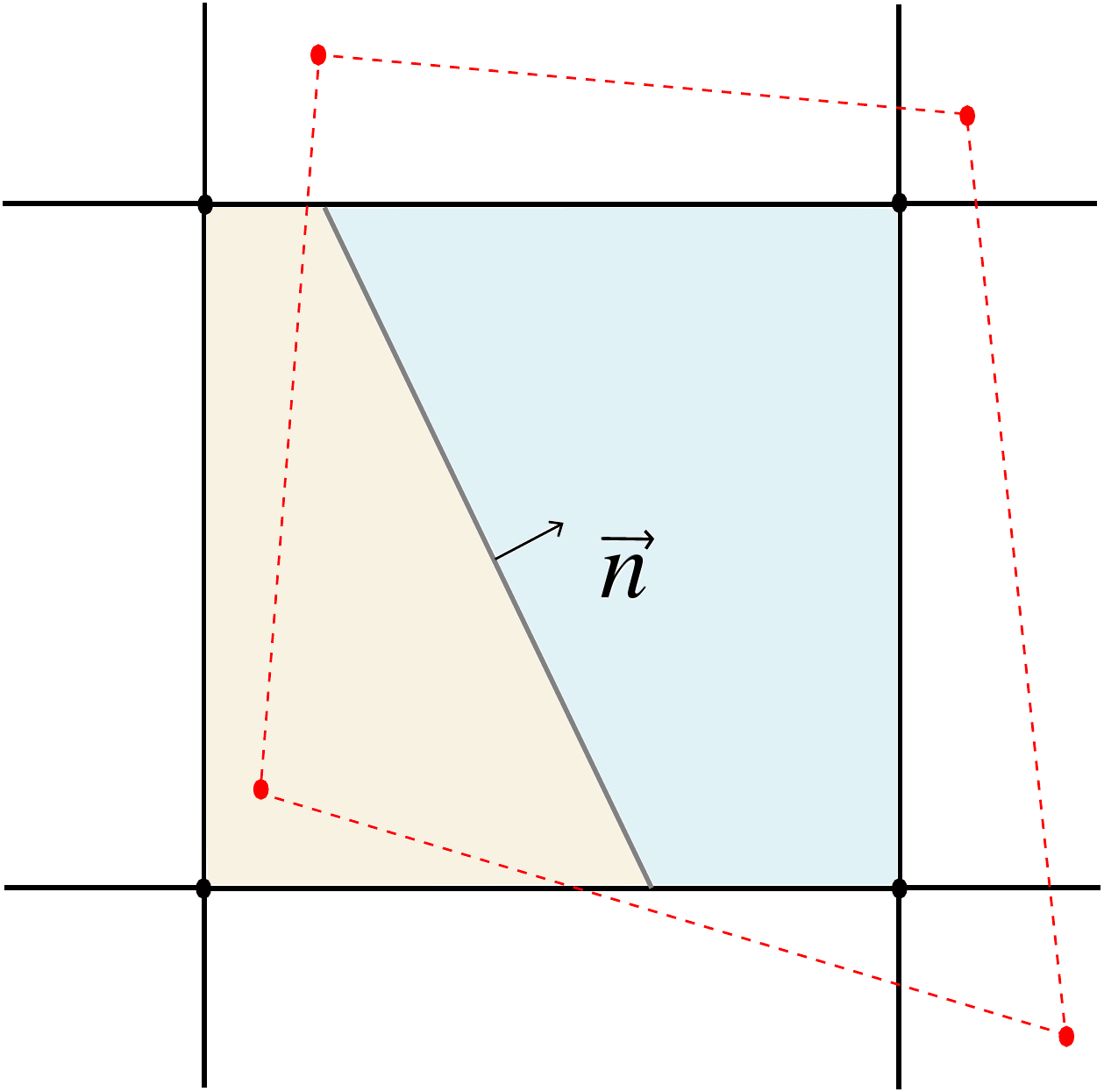}
\caption{Multi-material sharp interface model.}
\label{fig:multimat_interface}
\end{center}
\end{figure}

As explained in the introduction, the model of interface in one cell is a segment unequivocally defined from volume fractions and the normal vector to this interface (in the case of two materials). A method to determine the normal vector to the interface is due to Youngs in \cite{youngs1982time}. If $\mathbf{n}_\alpha$ and $k_\alpha$ respectively denote the outwards normal vector (in the sense of material $\alpha$) and the volume fraction:
\begin{equation*}
	\mathbf{n}_\alpha = \frac{\nabla k_\alpha}{\parallel \nabla k _\alpha \parallel}
\end{equation*}
Then, given the volume fraction in each cell of the domain, by computing a discrete gradient (based on a 9 points square-shaped stencil), the normal vector is determined in every mixed-cell of the domain. \\
The knowledge of normal vectors and volume fractions enables to place the interfaces on each mixed-cell, see Figure~\ref{fig:multimat_interface}. Let us  now perform the remap.

\vspace{0.5cm}

\textit{Remark 1:} The interface is not continuous on the domain.

\textit{Remark 2:}  Since $k_\alpha$ remains constant during the whole Lagrangian phase, $\mathbf{n}_\alpha$ is also assumed to be constant.

	\subsection{Alternate Directions Remap}

	\paragraph{$X$-Remap phase}

$\:$ \\

\begin{figure}[h]
\begin{center}
\includegraphics[width=0.45\textwidth, height=0.35\textwidth]{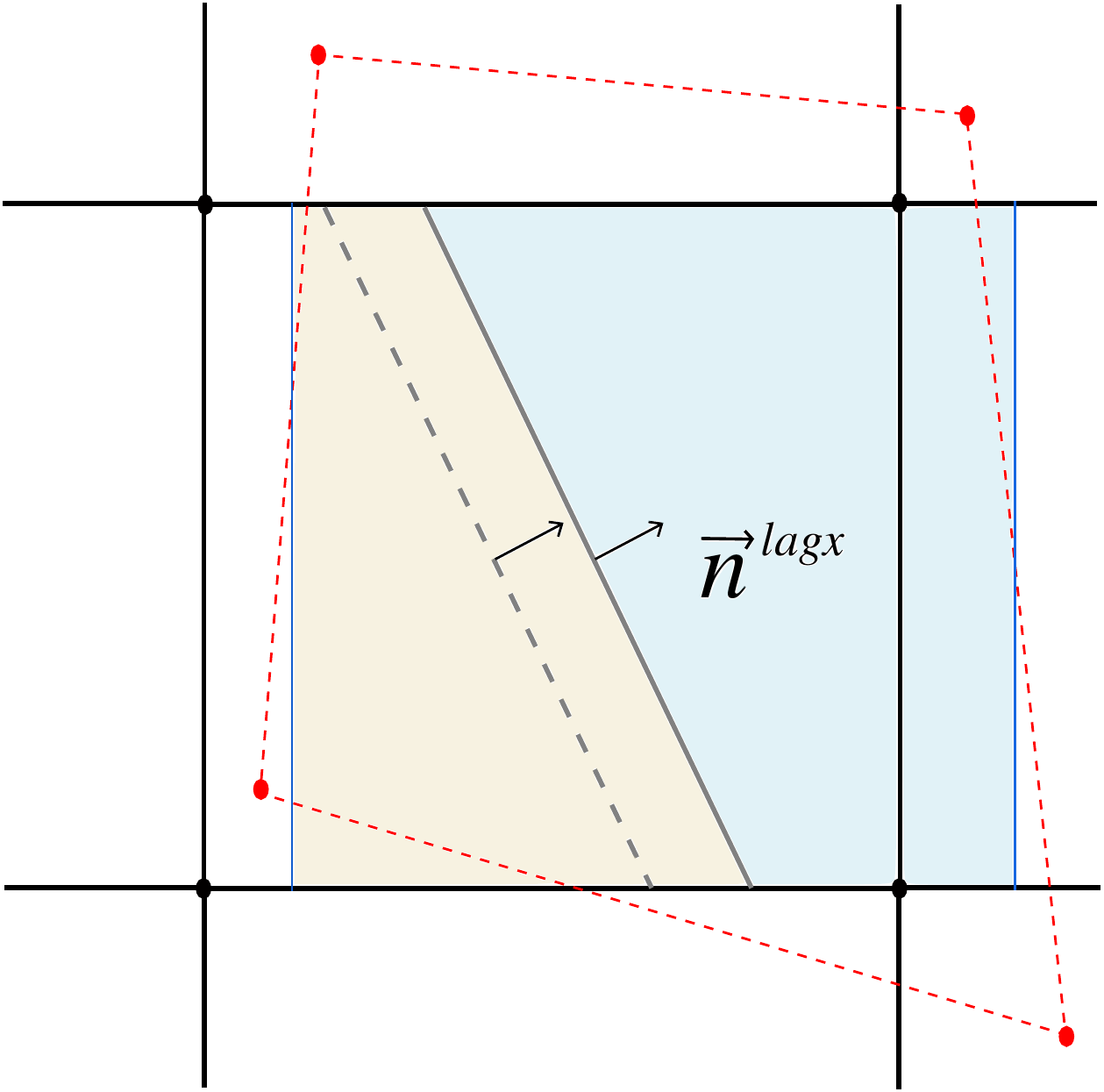}
\caption{$X$-Remap. Intersection between the interface and the volume fluxes.}
\label{fig:adi_x_multimat}
\end{center}
\end{figure}

Similarly to what is done in the mono-material case, the goal here is to compute volume fluxes of each material, i.e.\ the $dVol_{\alpha f}$, at all $X$-faces. The principle remains the same: Lagrangian displacement of the mesh along $X$, then remap. Node velocities give the value of $Vol^{lagx}$, a geometric representation of this volume, which is a rectangle thanks to the rectangular approximation of volume fluxes. \\
Taking into consideration the hypothesis exposed in the previous section, $k_{\alpha c}^{lagx} = k_{\alpha c}^n$ and $\mathbf{n}_{\alpha c}^{lagx} = \mathbf{n}_{\alpha c}^n$ are known. Thus, the interface can be placed on the cell at $t^{lagx}$, see Figure~\ref{fig:adi_x_multimat}.  

\vspace{0.5cm}

\begin{figure}[h]
\begin{center}
\includegraphics[width=0.45\textwidth, height=0.35\textwidth]{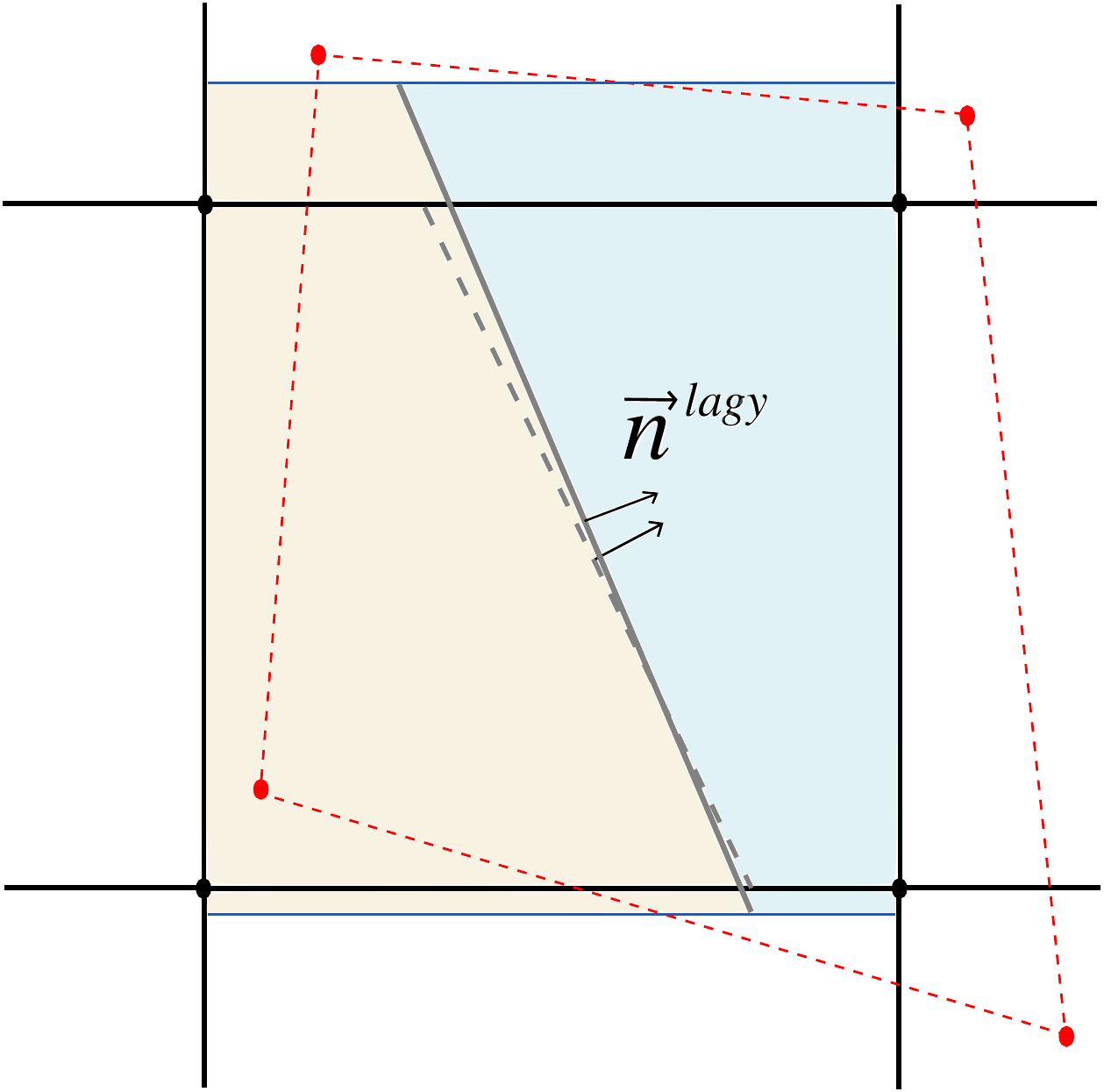}
\caption{$Y$-Remap. Intersection between the interface and the volume fluxes.}
\label{fig:adi_y_multimat}
\end{center}
\end{figure}

Now the interface is positioned on the mesh, one is able to compute its intersection with the volume fluxes at $X$-faces, i.e.\ the $dVol_f$. This simple intersection of a rectangle by a line enables to determine the $dVol_{\alpha f}$, see Figure~\ref{fig:adi_y_multimat}. \\
If the intersection is empty, as it is the case in Figure~\ref{fig:adi_x_multimat}, then volume fluxes are pure. At $X$-face $f$, in order to determine the material for the pure volume flux, one computes the scalar product:
\begin{equation*}
	(\mathbf{x}_f^{lagx} - \mathbf{x}_i^{lagx}) \cdot \mathbf{n}_{\alpha c}^{lagx}
\end{equation*}
where $\mathbf{x}_f^{lagx}$ and $\mathbf{x}_i^{lagx}$ respectively denote the middles of $X$-face $f$ and the interface $i$. The sign of this quantity enables to conclude.

Once the $dVol_{\alpha f}$ have been determined, variable can be reconstructed at first or second order at the faces for each material. Unlike the mono-material case, the volume of each material has to be remapped, this gives $k_{\alpha c}^{projx}$. The rest is strictly equivalent as what is done in the mono-material case. Let us write down the equations:

\begin{equation*}
	Vol_{\alpha c}^{projx} = Vol_{\alpha c}^{lagx} +  dVol_{\alpha l} - dVol_{\alpha r}
\end{equation*}
\vspace{0.05cm}
\begin{equation*}
	k_{\alpha c}^{projx} = \frac{Vol_{\alpha c}^{projx} }{Vol_c^n} = \frac{Vol_{\alpha c}^{projx} }{\Delta x \Delta y} 
\end{equation*}
\vspace{0.05cm}
\begin{equation*}
	m_{\alpha c}^{projx} = m_{\alpha c}^{lagx} + \rho_{\alpha l}^{lagx \:o2} dVol_{\alpha l} - \rho_{\alpha r}^{lagx \:o2} dVol_{\alpha r}
\end{equation*}
\vspace{0.05cm}
\begin{equation*}
	m_{\alpha c}^{projx}e_{\alpha c}^{projx} = m_{\alpha c}^{lagx}e_{\alpha c}^{lagx} + \rho_{\alpha l}^{lagx \:o2}e_{\alpha l}^{lagx \:o2} dVol_{\alpha l} - \rho_{\alpha r}^{lagx \:o2}e_{\alpha r}^{lagx \:o2} dVol_{\alpha r}
\end{equation*}

\textit{Remark:} Compared to the mono-material case, another condition is added for the limitation of second order linear reconstructions at faces. If the face is in the neighbourhood of a contact discontinuity, i.e.\ if it is not surrounded only by pure cells, the reconstruction degenerates to order 1.
\vspace{1cm}

	\paragraph{$Y$-Remap phase}

$\:$ \\
 
Since volume fractions have been remapped, the hypothesis of the model at this time gives: $k_{\alpha c}^{lagy} = k_{\alpha c}^{projx}$ and $\mathbf{n}_{\alpha c}^{lagy} = \mathbf{n}_{\alpha c}^{projx}$, that are in general different from $k_{\alpha c}^n$ and $\mathbf{n}_{\alpha c}^n$. In Figure~\ref{fig:adi_y_multimat}, the interface has rotated.

\vspace{0.5cm}

Apart from that point, all the calculations are similar to those of the $X$-remap.

	\subsection{Direct Remap}

\hspace{1.5em}Let us recall that in the mono-material case, volume fluxes at the faces are exactly the same as those in AD remap. But in the multi-material case, because of the interface positioning, fluxes are different. Indeed, what is performed in this remap can be summarized as follows. Let us define $i^n$ as the interface at $t^n$. 

The computation of the $dVol_{\alpha f}$ along both $X$ and $Y$ directions are performed at the same time. Along $X$, volume fluxes at $X$-faces are computed, the interface $i^n$ is displaced along $X$ (this new interface is denoted $i^n_X$) such that the volume fraction remains equal to $k_{\alpha c}^n$ in the $X$-Lagrangian cell (which is the same as in the $X$-remap in the AD case), and partial volume fluxes are deduced at $X$-faces. Along $Y$, volume fluxes at $Y$-faces are computed, the interface $i^n$ is displaced along $Y$ (this new interface is denoted $i^n_Y$) such that the volume fraction remains equal to $k_{\alpha c}^n$ in the $Y$-Lagrangian cell (which is the same as in the $Y$-remap in the  AD case), and partial volume fluxes are deduced at $Y$-faces.

\textit{Remark:} Two different interfaces $i^n_X$ and $i^n_Y$ are used respectively to determine partial volume fluxes at $X$ and $Y$-faces. They both have the same normal vector as $i^n$,  unlike the AD remap.  

Finally, the formula for the remap of the volume gives:
\begin{equation*}
	Vol_{\alpha c}^{projx} = Vol_{\alpha c}^{lagx} +  dVol_{\alpha l} - dVol_{\alpha r} +  dVol_{\alpha b} - dVol_{\alpha t}
\end{equation*}

The other equations are logically deduced from this one.

	\subsection{Direct Remap with Corner Fluxes}

\begin{figure}[h]
\begin{center}
\includegraphics[width=0.45\textwidth, height=0.35\textwidth]{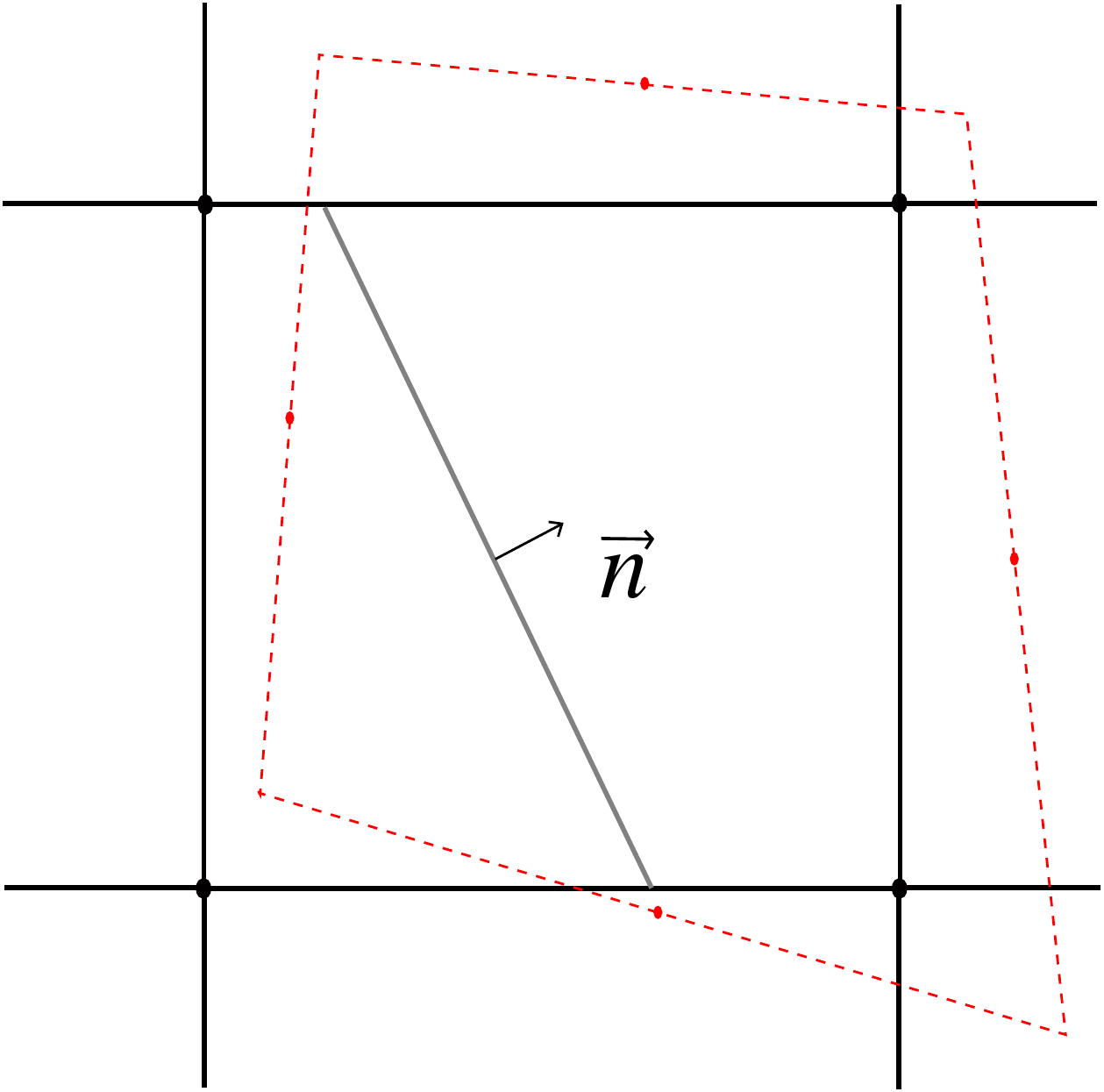}\hspace{0.05\textwidth}
\includegraphics[width=0.45\textwidth, height=0.35\textwidth]{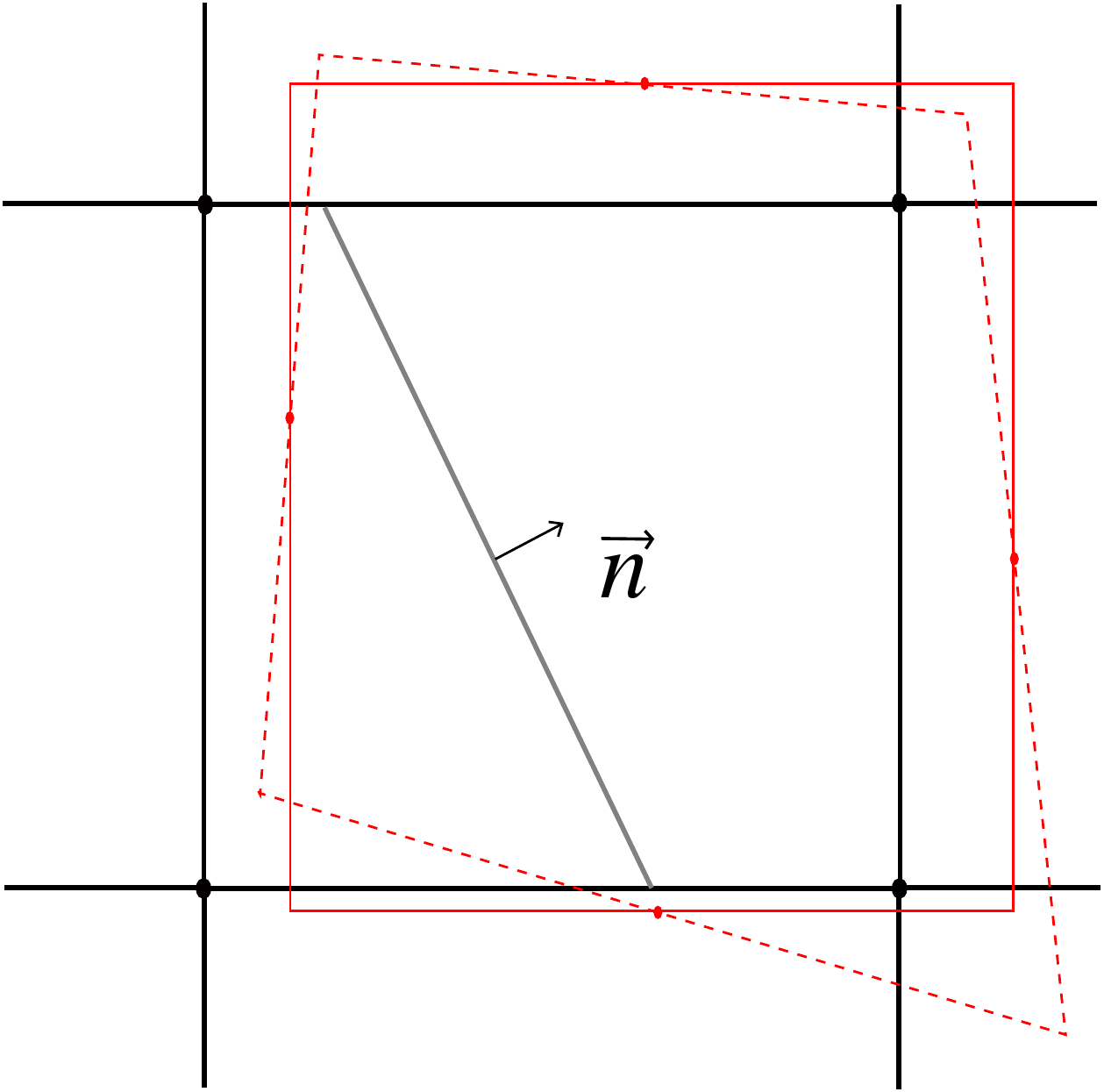}
\caption{Rectangular approximation of Lagrangian cell.}
\label{fig:rectangular_approximation}
\end{center}
\end{figure}

\hspace{1.5em}Since the main idea behind the DirectCF remap is the accuracy of geometric representation, one single Lagrangian interface will be used when intersecting all volume fluxes. In order to place this interface thanks to $k_{\alpha c}^n$, $\mathbf{n}_{\alpha c}^n$ and $Vol^{lagCF}$, an approximation has to be made. Indeed, the Lagrangian representation of the volume is quite complex, see Figure~\ref{fig:vol_lag_cf}. To avoid expensive calculations, the choice is to reconstruct the interface on a rectangular approximation of the Lagrangian cell. \\
This approximation is very simple: it is the only rectangle whom faces pass through the middles of the four Lagrangian faces, see Figure~\ref{fig:rectangular_approximation}.
	
\vspace{0.5cm}

\textit{Remark:} The error committed by this approximation on the lagrangian geometry is assumed to be reasonable because this rectangular volume is exactly the Lagrangian volume taken into consideration in the AD remap, see section 3.1.

\begin{figure}[h]
\begin{center}
\includegraphics[width=0.45\textwidth, height=0.35\textwidth]{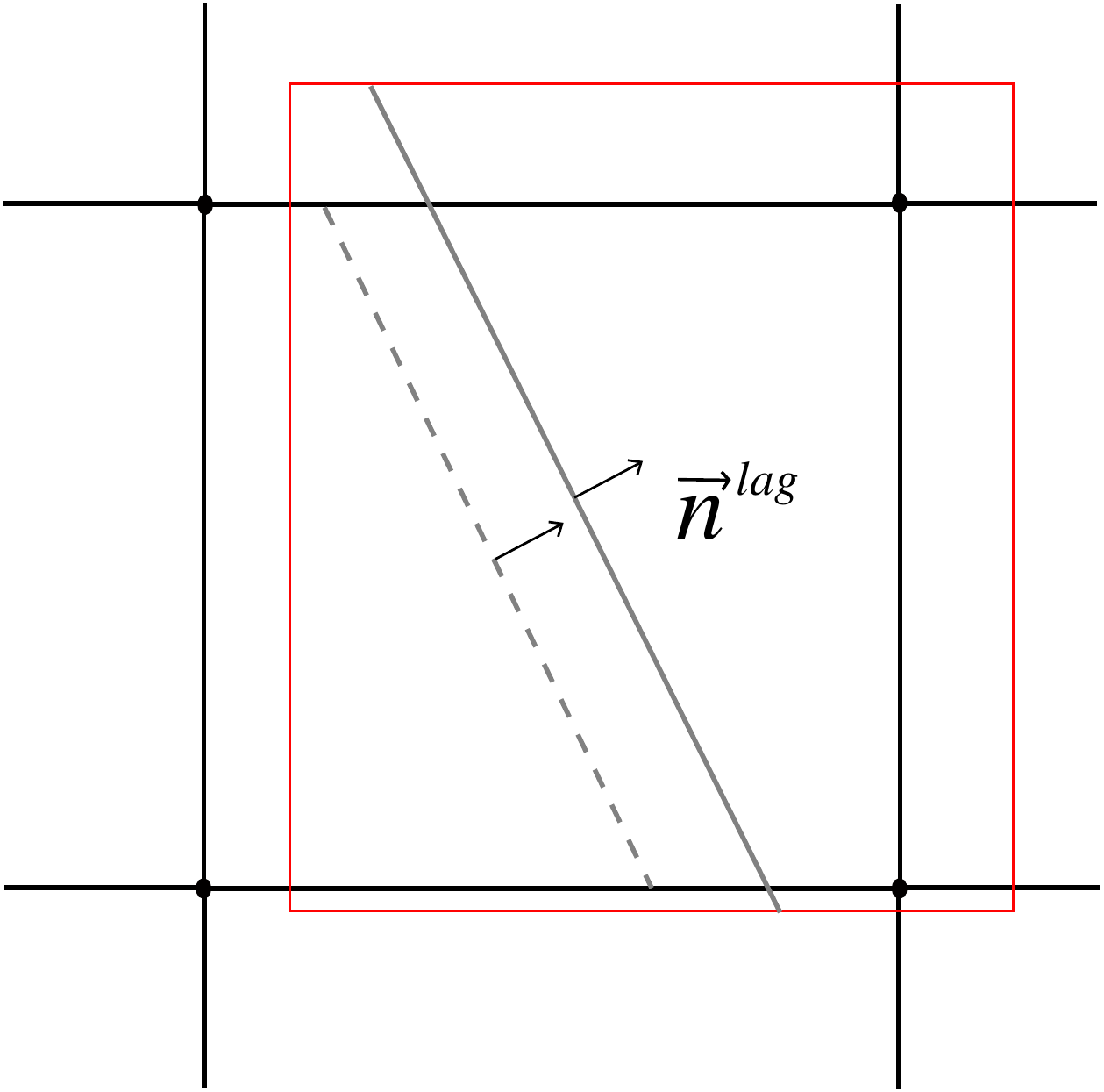}
\caption{Placement of Lagrangian interface.}
\label{fig:placement_interface}
\end{center}
\end{figure}

Now the approximated Lagrangian cell has been defined, and chosen rectangular, it is easy to place the interface using $k_{\alpha c}^n$, $\mathbf{n}_{\alpha c}^n$, see Figure~\ref{fig:placement_interface}.

\vspace{0.5cm}

\noindent Then, the equation of the interface enables to perform the intersection with all volume fluxes, at faces and corners. And the values of partial volume fluxes at faces and corners come out.

\vspace{0.5cm}
\begin{figure}[h]
\begin{center}
\includegraphics[width=0.45\textwidth, height=0.38\textwidth]{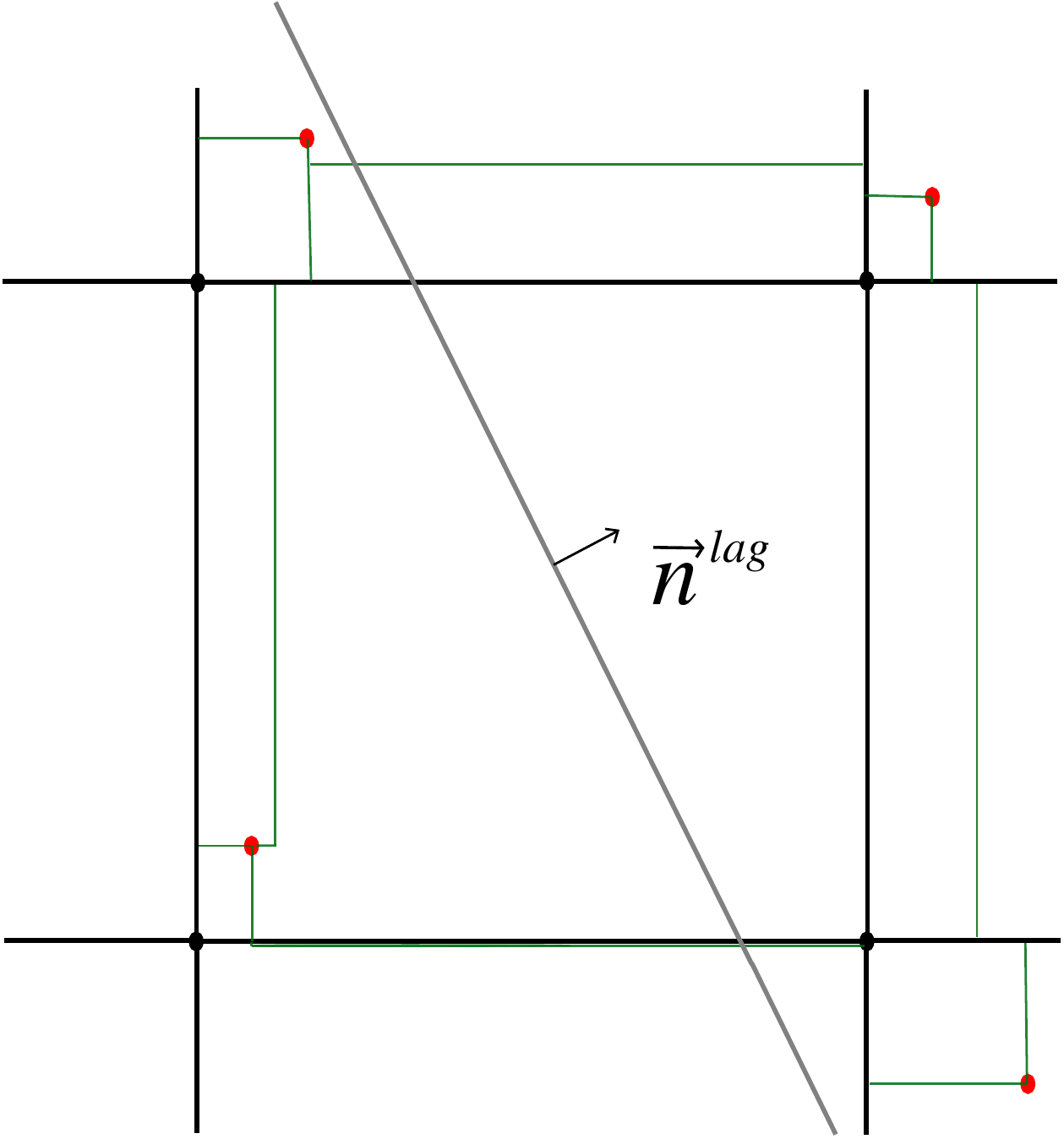}\hspace{0.05\textwidth}
\includegraphics[width=0.45\textwidth, height=0.38\textwidth]{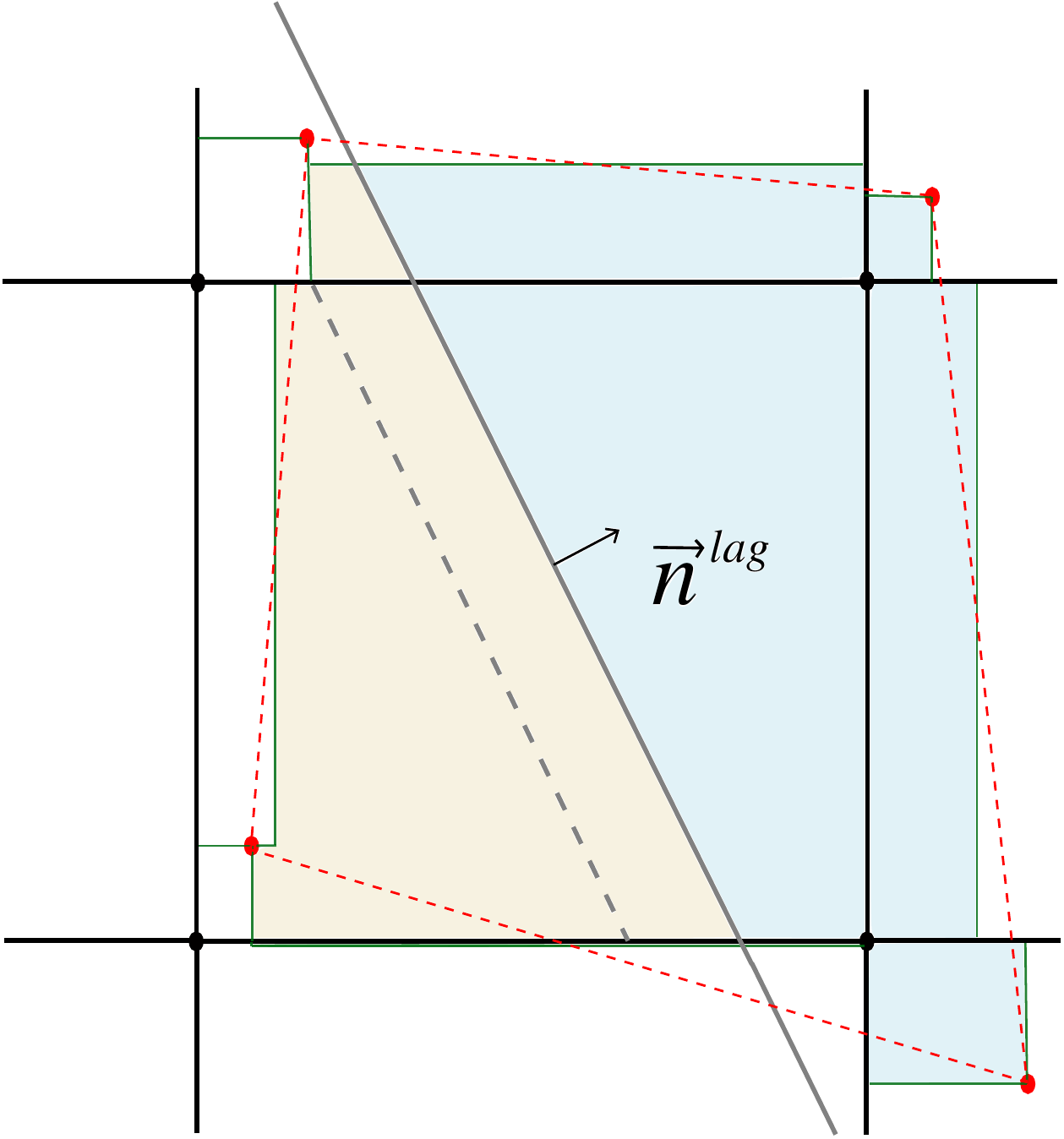}
\caption{DirectCF remap. Intersection between the interface and volume fluxes.}
\label{fig:dvol_fraction}
\end{center}
\end{figure}

\textit{Remark 1:} In the same manner as for the faces, in the case of a pure corner flux at node $p$, the sign of $(\mathbf{x}_p^{n+1/2} - \mathbf{x}_i^{lagCF}) \cdot \mathbf{n}_{\alpha c}^{lagCF}$ indicates the material.

\textit{Remark 2:} Reconstructing the interface on a volume which is different from the one that enables to compute volume fluxes does not ensure the constancy of volume fractions during the Lagrangian phase. Such a little deviation from the isocompressibility model could create artificial compressions or expansions of partial volumes. In order to fix this, the values of lagrangian partial volumes actually considered by the remap (\i.e. yellow and blue volumes in Figure~\ref{fig:dvol_fraction}) need to be known. The too many possible situations make this computation extremely costly. Thus, this error will be neglected. In practice, even on severe benchmarks, it actually seems to be negligible.

\textit{Remark 3:} In this case, there are one single normal and one single interface for all directions and volume fluxes, which reduces the calculation time. It is also a gain of accuracy since a single-interface (per iteration) model is obviously closer to the exact case than a multiple-interface one. 

Let us finally write the remap equation of the volume:
\begin{equation*}
	Vol_{\alpha c}^{proj} = Vol_{\alpha c}^{lagCF} +  \sum_{faces} \!dVol_{\alpha f} + \sum_{corners} \!\!\!dVol_{\alpha p} 
\end{equation*}

\noindent The rest follows from this formula.

 \clearpage
 \newpage
 \section{Numerical results}

\hspace{1.5em}In order to validate and verify the scheme including Direct remap with Corner Fluxes, several classical test-cases have been ran. They are presented in the following paragraphs. This section also aims at comparing this remap to the AD and Direct remaps. For each test-case, initial state is written, results are shown and some comments are made. The DirectCF implemented is second order at both faces and corners, with diagonal linear reconstruction at the corners.

	\subsection{First observations and convergence}

\hspace{1.5em}Firstly, the accuracy of the three remaps are tested on classical test-cases of linear advection and solid rotation, in both mono-material and multi-material cases.

		\subsubsection{Mono-material}

		\paragraph{Advection}
		
$\:$ \\

Linear advection and return of a square ($\rho = 10\: \mbox{kg.m}^{-3}$) into a field ($\rho = 0.1\: \mbox{kg.m}^{-3}$) with velocity $\mathbf{u} = (5\:; 5) \:\mbox{m.s}^{-1}$. \\
Meshes: $50\!\times\!50$, $100\!\times\!100$, $200\!\times\!200$, $400\!\times\!400$. \\
Mono-material, Perfect Gas $p=(\gamma-1)\rho e$ with $\gamma = 1.4$.

In this test-case, quantities are either constant, either transported. The goal is to measure the diffusivity of the scheme.

\vspace{0.7cm}
\begin{figure}[h]
\begin{center}
\includegraphics[width=0.25\textwidth, height=0.25\textwidth]{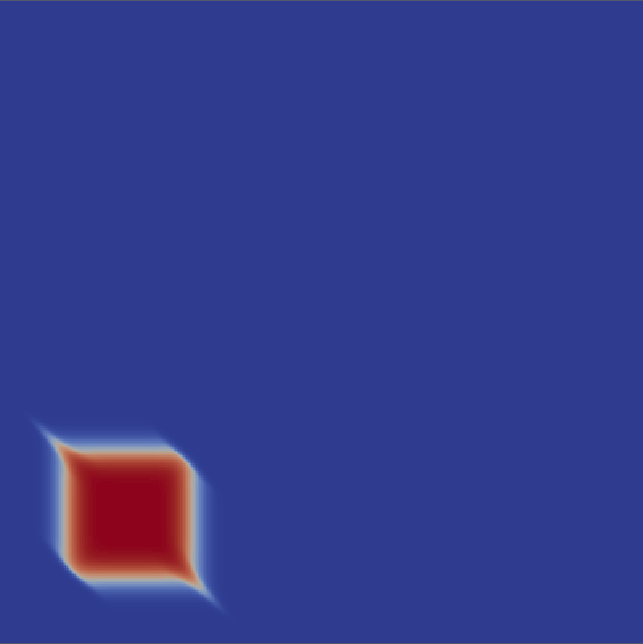}\hspace{0.1\textwidth}\includegraphics[width=0.25\textwidth, height=0.25\textwidth]{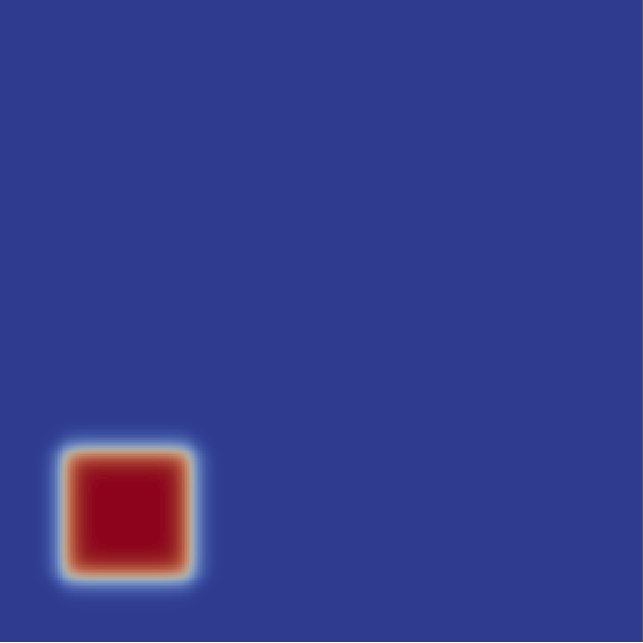}\hspace{0.1\textwidth}\includegraphics[width=0.25\textwidth, height=0.25\textwidth]{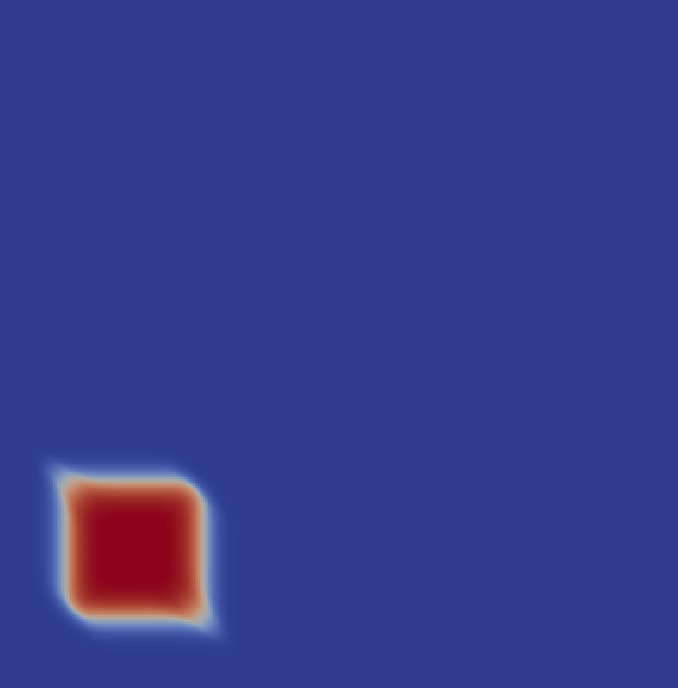}
\caption{Monomat advection. Direct (left), AD (centre) and DirectCF (right).}
\label{fig:adv_return_monomat}
\end{center}
\end{figure}

\begin{figure}[!h]
\begin{center}
\includegraphics[width=0.7\textwidth, height=0.4\textwidth]{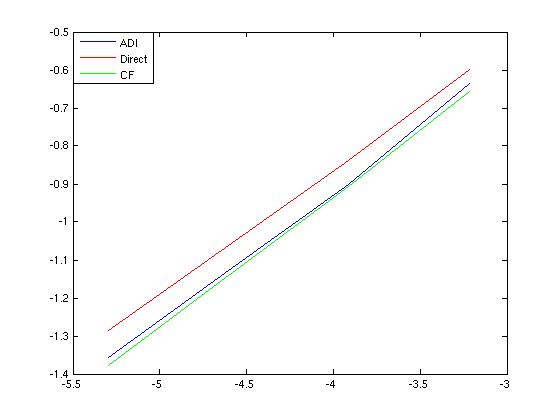}
\caption{Monomat advection. $\log \parallel \!\rho^{final} - \rho^0 \!\parallel_{L^2}$ in function of $\log (\Delta x)$}
\label{fig:error_adv_return_monomat}
\end{center}
\end{figure}

Figure~\ref{fig:adv_return_monomat} shows the density after one advection and return on a mesh of $400\!\times\!400$ cells. The Direct shows the same deformation as the DirectCF, even more marked, whereas ADI draws a nice square. Regarding the $L^2$ error on the density, see Figure~\ref{fig:error_adv_return_monomat}, despite its deformation at the corners along the $x=-y$ direction, the DirectCF error is quite comparable to AD, but the Direct is always less accurate, no matter the mesh.   

		\paragraph{Rotation}

$\:$ \\

Solid rotation ($2\pi$) of a square ($\rho = 10\: \mbox{kg.m}^{-3}$) into a field ($\rho = 0.1\: \mbox{kg.m}^{-3}$) with velocity $\mathbf{u}(r) = r \omega_0 \: \mathbf{e}_{\theta}$ and $\omega_0 = 2\pi \:\mbox{s}^{-1}$ at each iteration. \\
Meshes: $50\!\times\!50$, $100\!\times\!100$, $200\!\times\!200$, $400\!\times\!400$.
Mono-material, Perfect Gas $p=(\gamma-1)\rho e$ with $\gamma = 1.4$.

Since velocities are imposed at each time step, this test-case especially evaluates geometric properties of the remap.
 
\vspace{0.7cm}
\begin{figure}[!h]
\begin{center}
\includegraphics[width=0.25\textwidth, height=0.25\textwidth]{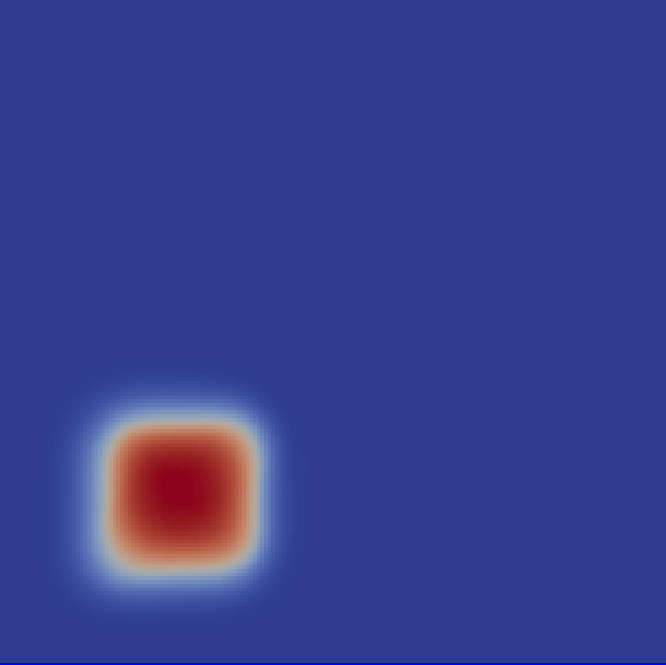}\hspace{0.1\textwidth}\includegraphics[width=0.25\textwidth, height=0.25\textwidth]{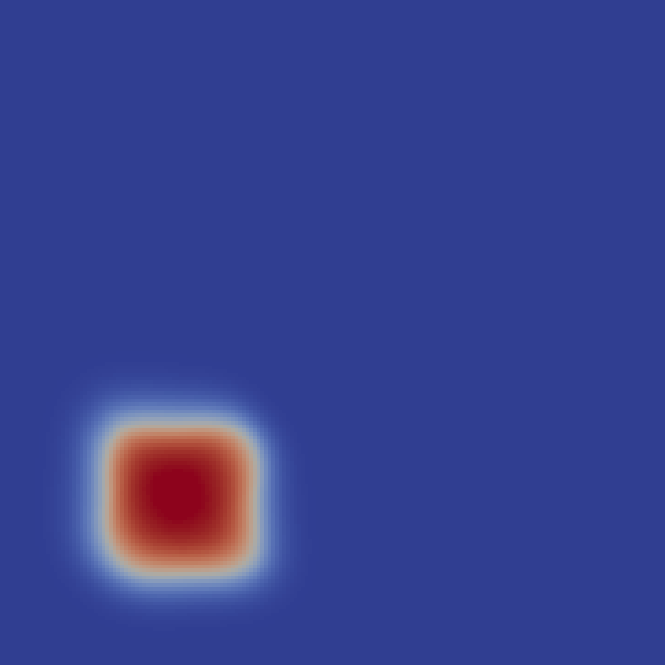}\hspace{0.1\textwidth}\includegraphics[width=0.25\textwidth, height=0.25\textwidth]{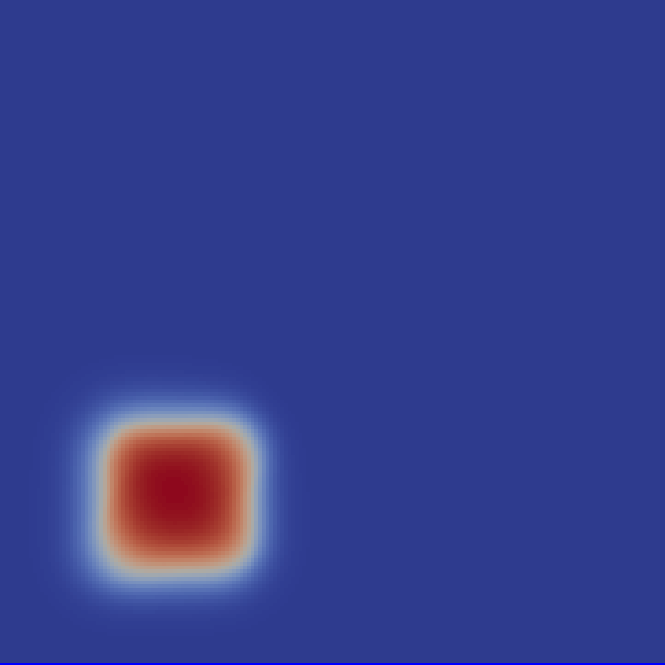}
\caption{Monomat rotation. Direct (left), AD (centre) and DirectCF (right).}
\label{fig:solrot_monomat}
\end{center}
\end{figure}

\begin{figure}[!h]
\begin{center}
\includegraphics[width=0.7\textwidth, height=0.38\textwidth]{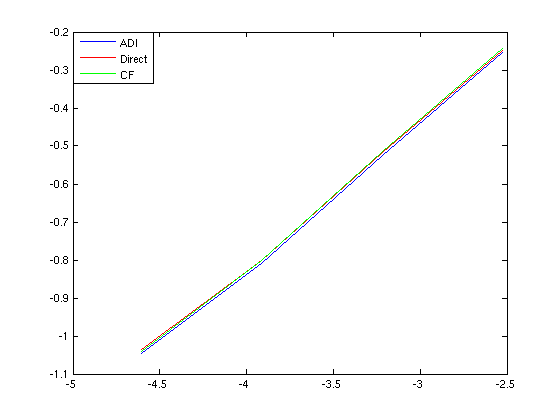}
\caption{Monomat rotation. $\log \parallel \!\rho^{final} - \rho^0 \!\parallel_{L^2}$ in function of $\log (\Delta x)$}
\label{fig:error_solrot_monomat}
\end{center}
\end{figure}

Here again, the density is represented on a $400\!\times\!400$ mesh after one rotation of angle $2\pi$, see Figure~\ref{fig:solrot_monomat}. The Direct has broken the symmetry with respect to $x=y$, which is not the case for the two other remaps, that seem to be comparable to each other. The curves plotted on Figure~\ref{fig:error_solrot_monomat} show that all remaps are equivalent in terms of the $L^2$ error.

		\subsubsection{Multi-material}
	
		\paragraph{Advection}
		
$\:$ \\

Linear advection and return of a square of air into air ($\rho_a = 1.29\: \mbox{kg.m}^{-3}$) with velocity $\mathbf{u} = (5\:; 5) \:\mbox{m.s}^{-1}$. \\
Meshes: $50\!\times\!50$, $100\!\times\!100$, $200\!\times\!200$, $400\!\times\!400$. \\
Multimaterial, air: Perfect Gas $p=(\gamma-1)\rho e$ with $\gamma = 1.4$

\vspace{0.7cm}
\begin{figure}[h]
\begin{center}
\includegraphics[width=0.25\textwidth, height=0.25\textwidth]{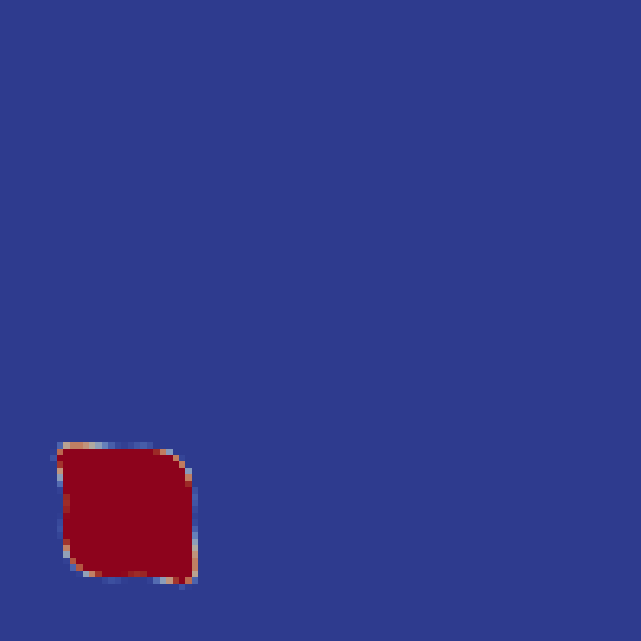}\hspace{0.1\textwidth}\includegraphics[width=0.25\textwidth, height=0.25\textwidth]{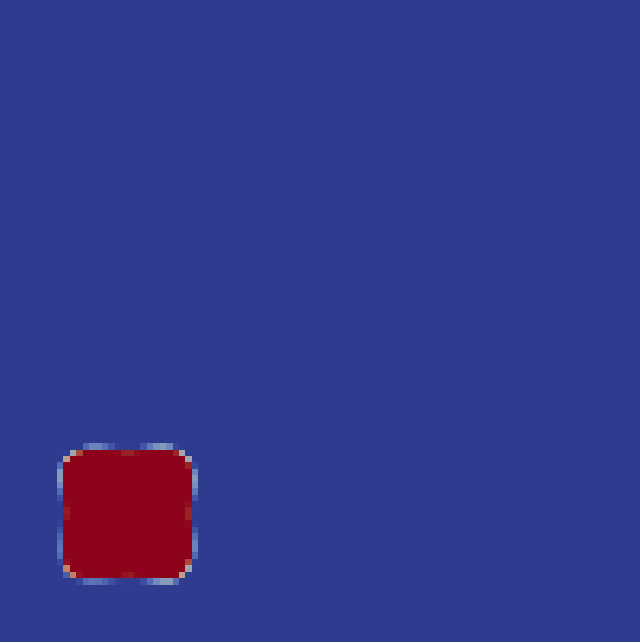}\hspace{0.1\textwidth}\includegraphics[width=0.25\textwidth, height=0.25\textwidth]{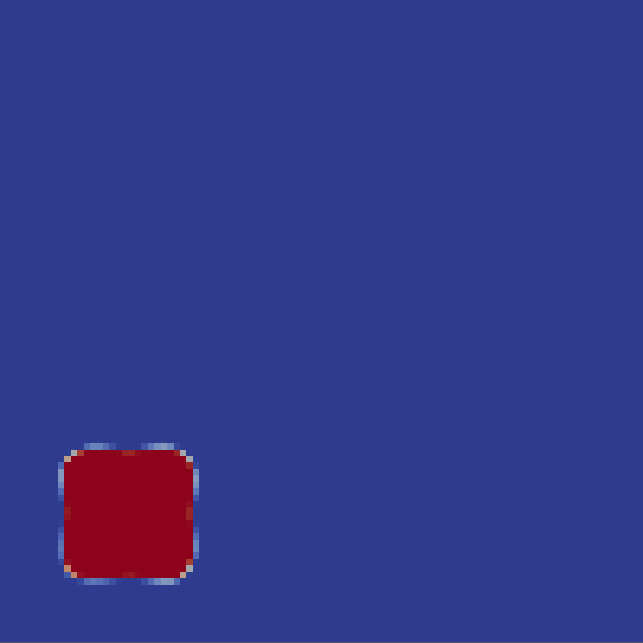}
\caption{Multimat advection. Direct (left), AD (centre) and DirectCF (right).}
\label{fig:adv_return_multimat}
\end{center}
\end{figure}

\begin{figure}[h]
\begin{center}
\includegraphics[width=0.7\textwidth, height=0.4\textwidth]{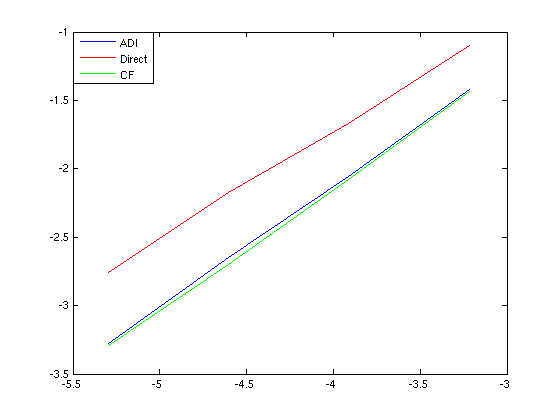}
\caption{Multimat advection. $\log \parallel \!k_1^{final} - k_1^0 \!\parallel_{L^2}$ in function of $\log (\Delta x)$}
\label{fig:error_adv_return_multimat}
\end{center}
\end{figure}

In Figure~\ref{fig:adv_return_multimat}, the volume fraction is represented after the two materials have been advected forwards then backwards along the $x=y$ direction on a $100\!\times\!100$ cells domain. The same difficulties as in the mono-material case are shown by the Direct remap, and they are even emphasized by sharp interface reconstruction. The corners in the direction of propagation are diffused whereas the corners in the orthogonal direction are spread. On the contrary, the AD and DirectCF remaps give nice-shaped results. When having a look at Figure~\ref{fig:error_adv_return_multimat}, it is quite clear that the Direct remap can be ruled out because of its lack of accuracy. The rest of this report will focus on the two other remaps.

	\subsection{Robustness}

\hspace{1.5em}Robustness is a mandatory characteristic for a scheme to be considered valid, especially in the industry, where in general stability and robustness prevail over accuracy. All simulations presented below are multi-material simulations because they strongly test scheme abilities.

		\paragraph{Water-air rotation}
		
$\:$ \\

Solid rotation ($2\pi$) of a square of water into  air with velocity $\mathbf{u}(r) = r \omega_0 \: \mathbf{e}_{\theta}$ and $\omega_0 = 2\pi \: 10^3 \:\mbox{s}^{-1}$ at each iteration. \\
Multimaterial, Air: Perfect Gas $p=(\gamma-1)\rho e$ with $\gamma = 1.4$ \\
Water: Stiffened Gas $p=(\gamma-1)\rho e - \pi$ with $(\gamma \:; \pi) = (7 \:; 2.1 \:10^9 \:\mbox{Pa})$ \\
Mesh: $400\!\times\!400$ on the domain $[0, 4]\times[0,4]$. 

\vspace{0.7cm}
\begin{figure}[h]
\begin{center}
\includegraphics[width=0.3\textwidth, height=0.3\textwidth]{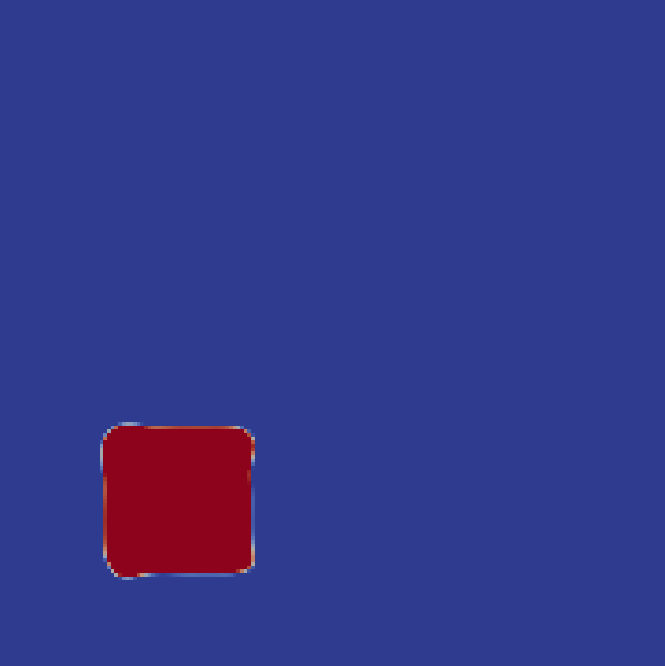}\hspace{0.2\textwidth}\includegraphics[width=0.3\textwidth, height=0.3\textwidth]{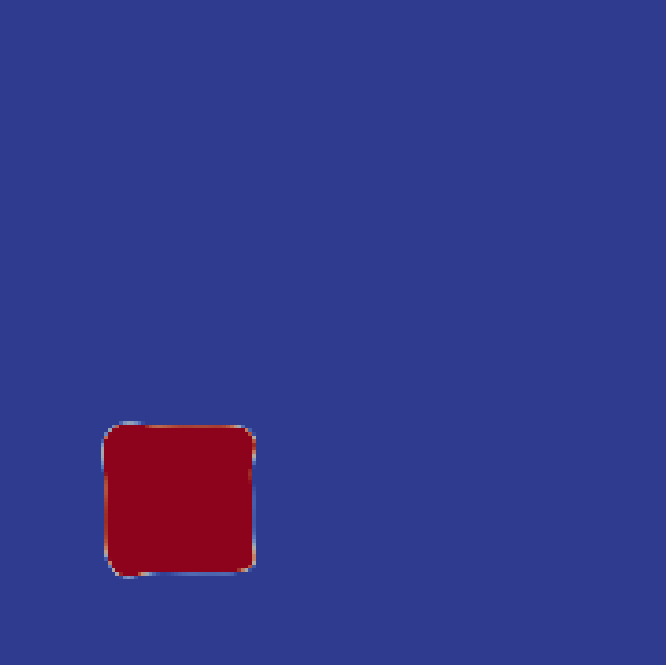}
\caption{Volume fraction at $t=1$ms. AD (left) and DirectCF (right).}
\label{fig:solrot_water_air}
\end{center}
\end{figure}

Even if a simple solid rotation is simulated, this tests robustness because of the gap between the respective densities of air and water, little mistakes can lead to noticeable effects on the variables. By the way, since water sound speed is very high, the time step is very small, which enlightens the defaults of the schemes. However, Figure~\ref{fig:solrot_water_air} shows almost no difference between the two results.

\vspace{2cm}
		\paragraph{HAAS test \cite{haas1987interaction}}
		
$\:$ \\

Interaction between a shock in air and a bubble of Helium. Mesh $1000\!\times\!90$ on the domain $[0, 1000]\!\times\![0, 9]$ cm. \\
Multimaterial, both Helium and air are Perfect Gases. \\
All initial data is gathered in the tabular below.

\vspace{0.5cm}
\footnotesize
\begin{figure}[h]
\begin{center}
\begin{tabular}{|c|c|c|c|}
\hline
Initial state & Left air state (shock) & Right air state & Bubble \\
\hline
Density $\rho$ (kg.m$^3$) & $1.376363$ & $1$ &  $0.18187$ \\
\hline
Velocity $\mathbf{u}.\mathbf{e}_x$ (m.s$^{-1}$) & $124.824$ & $0$ & $0$ \\
\hline
Pressure $p$ (Pa) & $1.5698 \:10^5$ & $10^5$ & $10^5$ \\
\hline
Gamma $\gamma$ & $1.4$ & $1.4$ & $1.66$ \\
\hline
\end{tabular}
\end{center}
\end{figure}
\normalsize

\noindent At the initial state, the bubble is at rest and the shock is propagating with velocity over $100$ m.s$^{-1}$.

\vspace{0.5cm}
\begin{figure}[h]
\begin{center}
\includegraphics[width=0.8\textwidth]{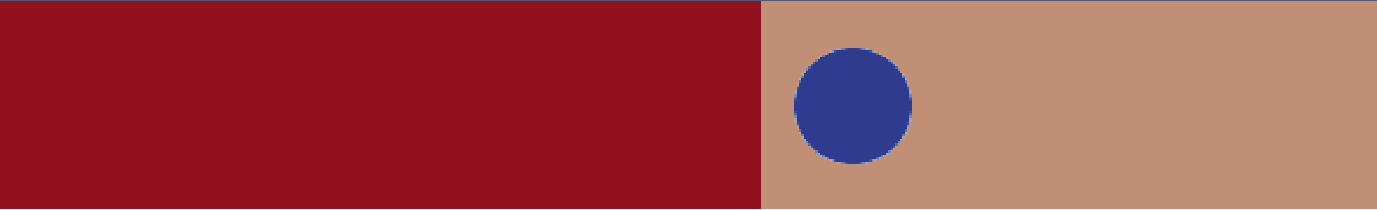}
\end{center}
\begin{center}
\includegraphics[width=0.8\textwidth]{HAAS_init.png}
\end{center}
\caption{Density at initial state, $t=0$ s. AD top, DirectCF bottom.}
\label{fig:HAAS_init}
\end{figure}

\begin{figure}[h]
\begin{center}
\includegraphics[width=0.8\textwidth]{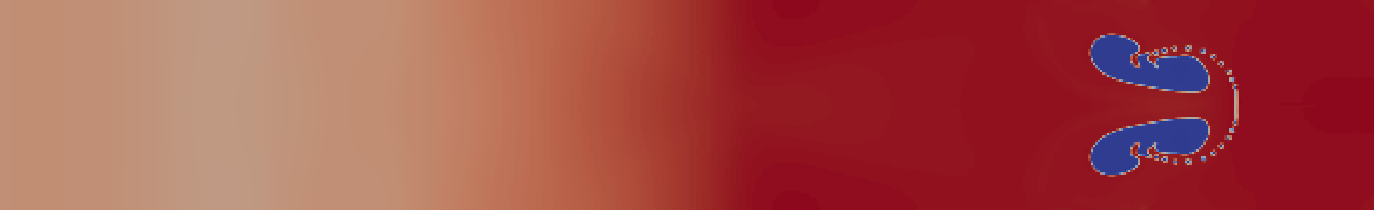}
\end{center}
\begin{center}
\includegraphics[width=0.8\textwidth]{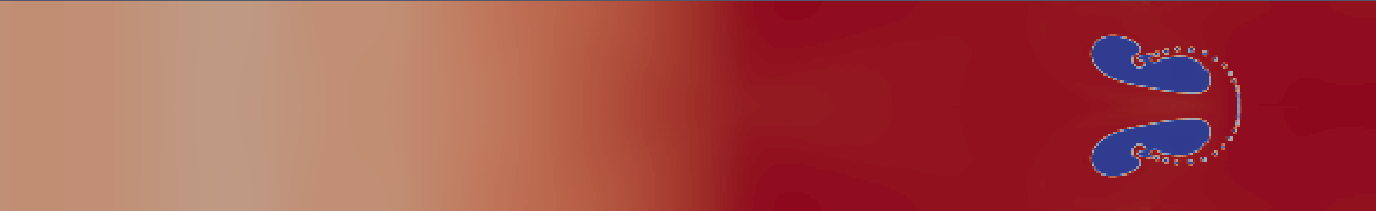}
\end{center}
\caption{Density at final state, $t=1$ ms. AD top, DirectCF bottom.}
\label{fig:HAAS_final}
\end{figure}

\noindent The propagation of the shock through the bubble creates a Richtmyer-Meshkov instability, see Figure~\ref{fig:HAAS_final}.

\vspace{0.5cm}
$\:$ \\
\begin{figure}[h]
\begin{center}
\includegraphics[height=0.3\textwidth]{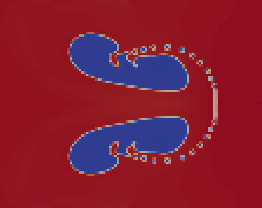}\hspace{0.15\textwidth}
\includegraphics[height=0.3\textwidth]{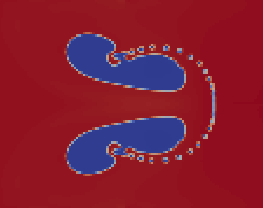}
\end{center}
\caption{Zoom, density at final state, $t=1$ ms. AD left, DirectCF right.}
\label{fig:HAAS_zoom}
\end{figure}

Both schemes give the same mushroom-shaped instability, which is the result expected. However, one can observe some differences, in particular at the extremities of the trickle.

		\paragraph{Water-water impact}
		
$\:$ \\

 Impact Equal Density Test: Impact of a water drop into air on a water wall. Mesh $320\!\times\!160$ on the domain $[0, 10]\!\times\![0, 5]$ cm. \\
Multimaterial, air Perfect Gas, and water Stiffened Gas.
All initial data is gathered in the tabular below.

This test simulates a shock involving materials with very different densities and compressibilities. All figures have been represented using mirror symmetry so that on each one, results given by the AD remap are on the left, and the symmetric of those given by the DirectCF remap are  on the right.

\vspace{0.5cm}
\footnotesize
\begin{figure}[h]
\begin{center}
\begin{tabular}{|c|*{3}{p{1.95cm}|}}
\hline
Initial state & \centering Air & \centering Wall & \centering Drop \tabularnewline
\hline
Density $\rho$ (kg.m$^3$) & \centering $1.29$ & \centering $1000$ &  \centering $1000$ \tabularnewline
\hline
Velocity $\mathbf{u}.\mathbf{e}_x$ (m.s$^{-1}$) & \centering $-1000$ & \centering $0$ & \centering $-1000$ \tabularnewline
\hline
Pressure $p$ (Pa) & \centering $10^5$ & \centering $10^5$ & \centering $10^5$ \tabularnewline
\hline
Gamma $\gamma$ & \centering $1.4$ & \centering $7$ & \centering $7$ \tabularnewline
\hline
Pi $\pi$ (Pa) & \centering $0$ & \centering $2.1 \: 10^9$ & \centering $2.1 \: 10^9$ \tabularnewline
\hline\end{tabular}
\end{center}
\end{figure}
\normalsize

\begin{figure}[h]
\begin{center}
\includegraphics[width=0.23\textwidth, height=0.23\textwidth]{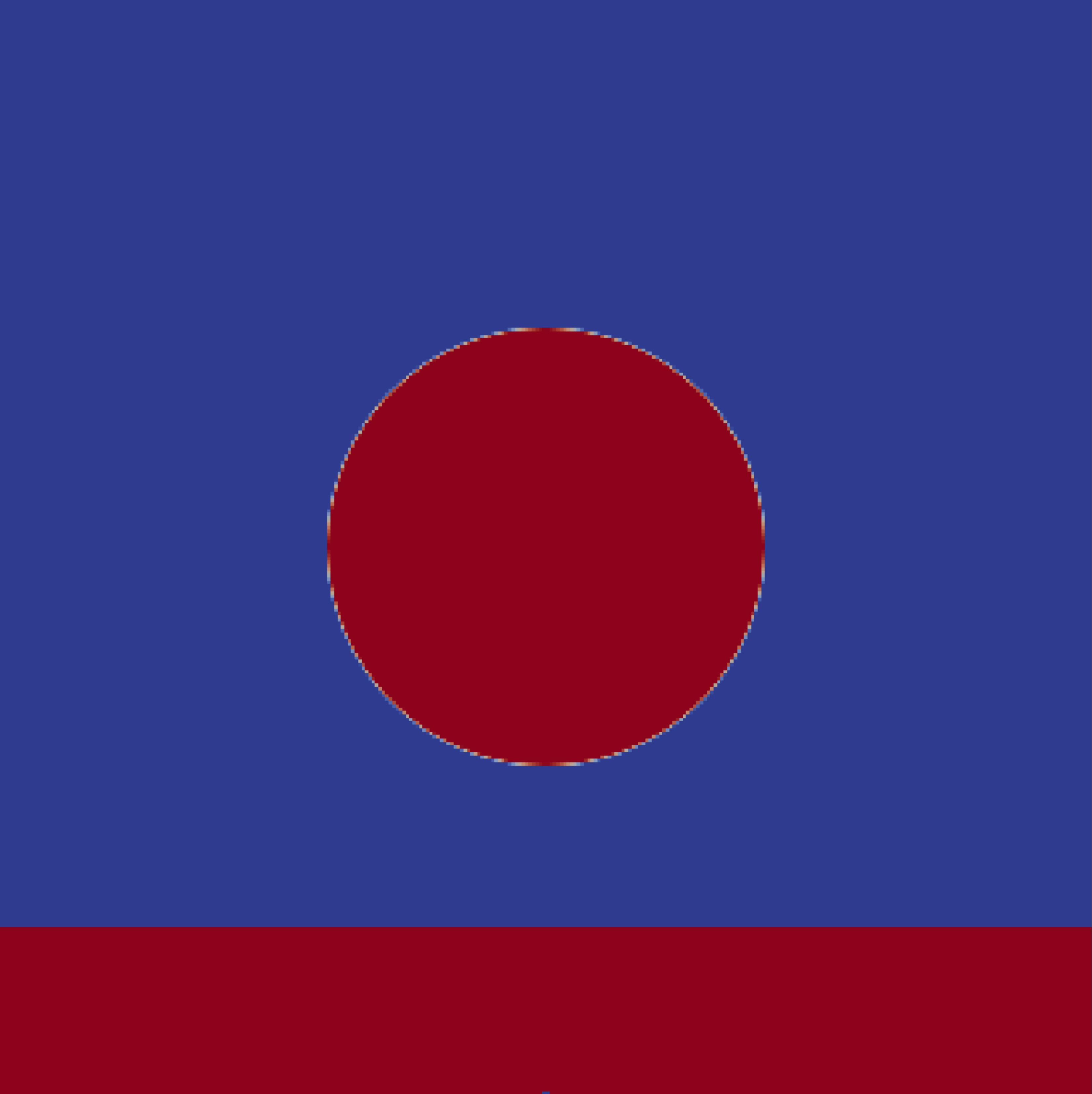}\hspace{0.015\textwidth}
\includegraphics[width=0.23\textwidth, height=0.23\textwidth]{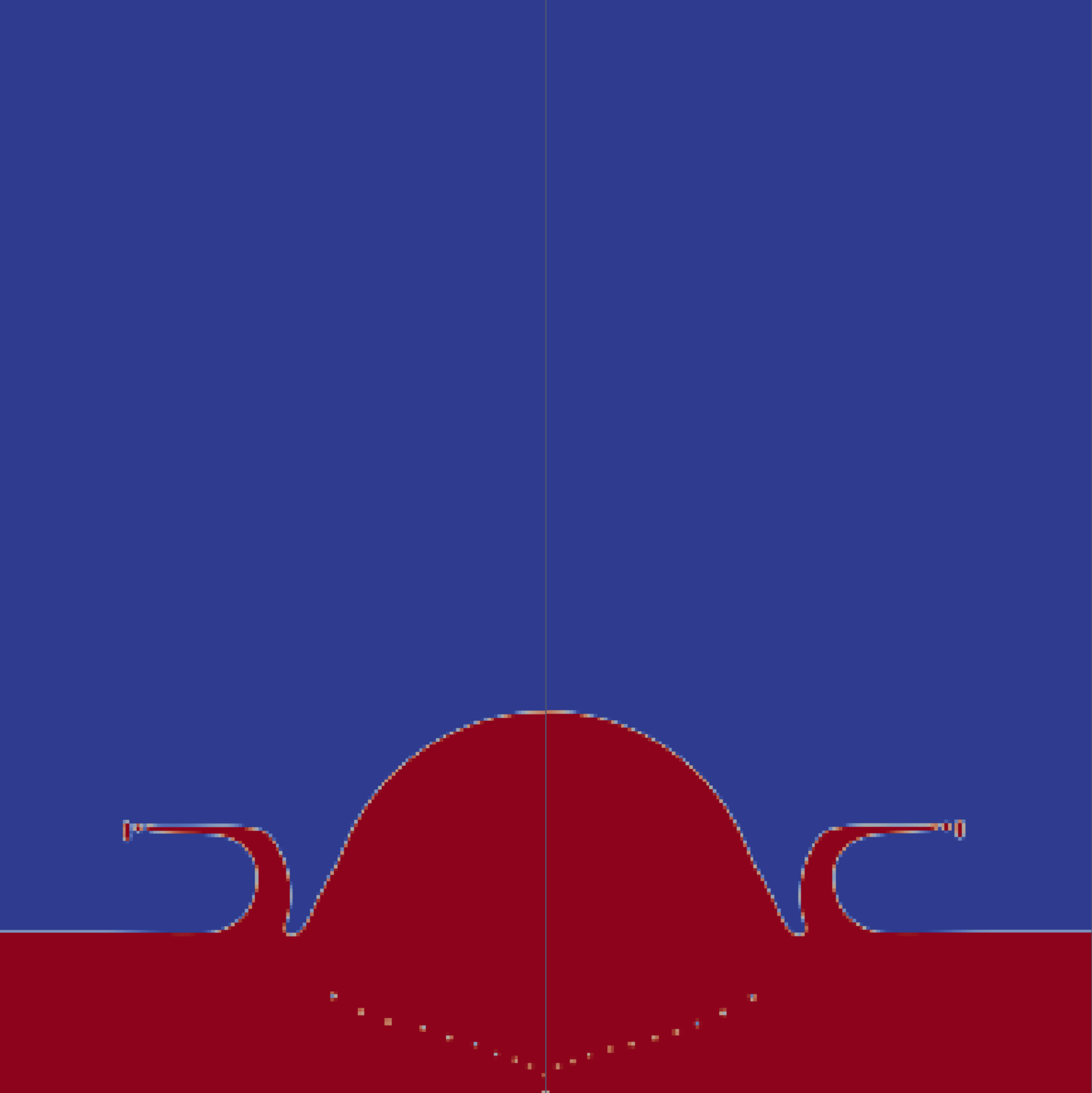}\hspace{0.015\textwidth}
\includegraphics[width=0.23\textwidth, height=0.23\textwidth]{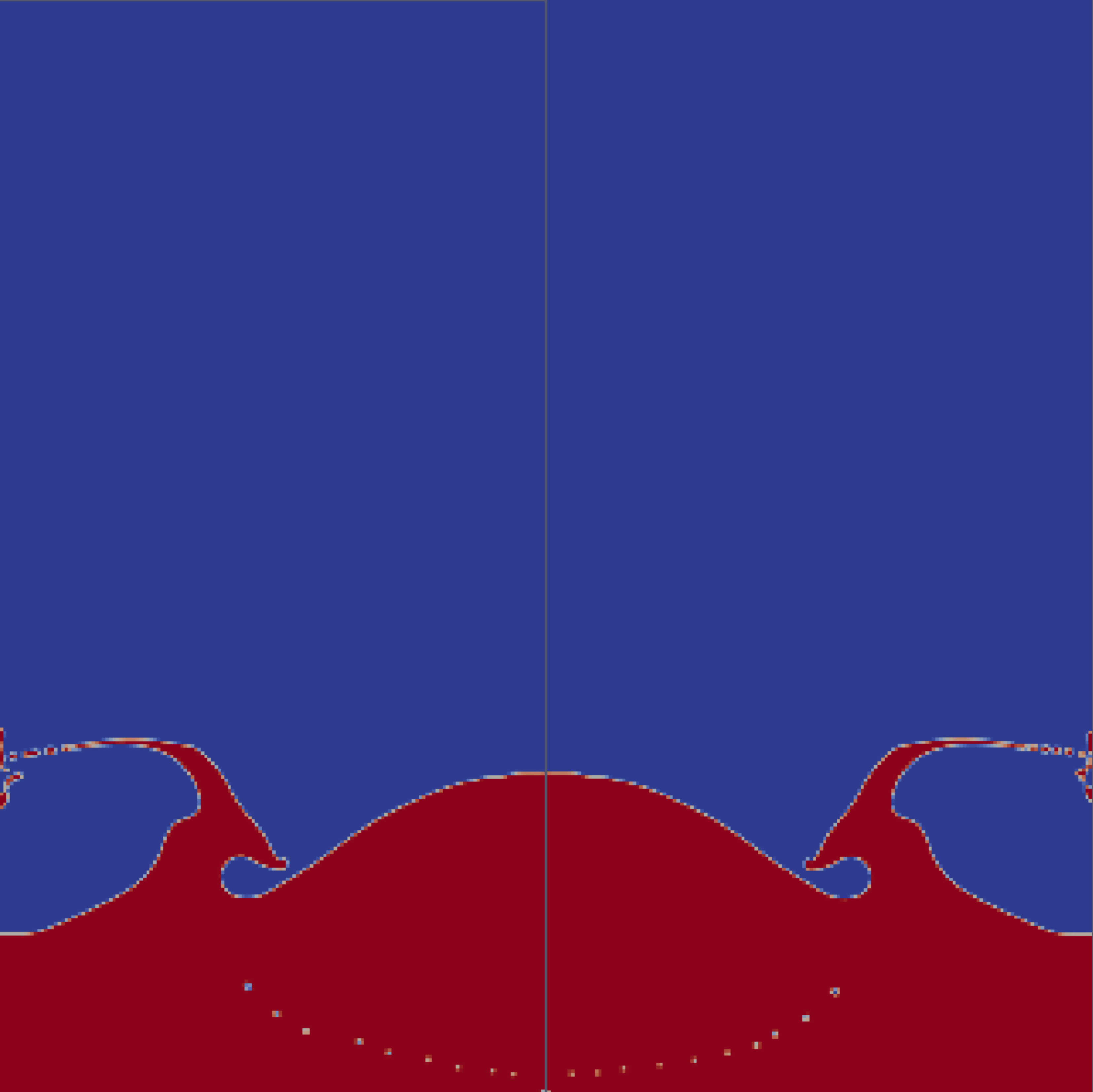}\hspace{0.015\textwidth}
\includegraphics[width=0.23\textwidth, height=0.23\textwidth]{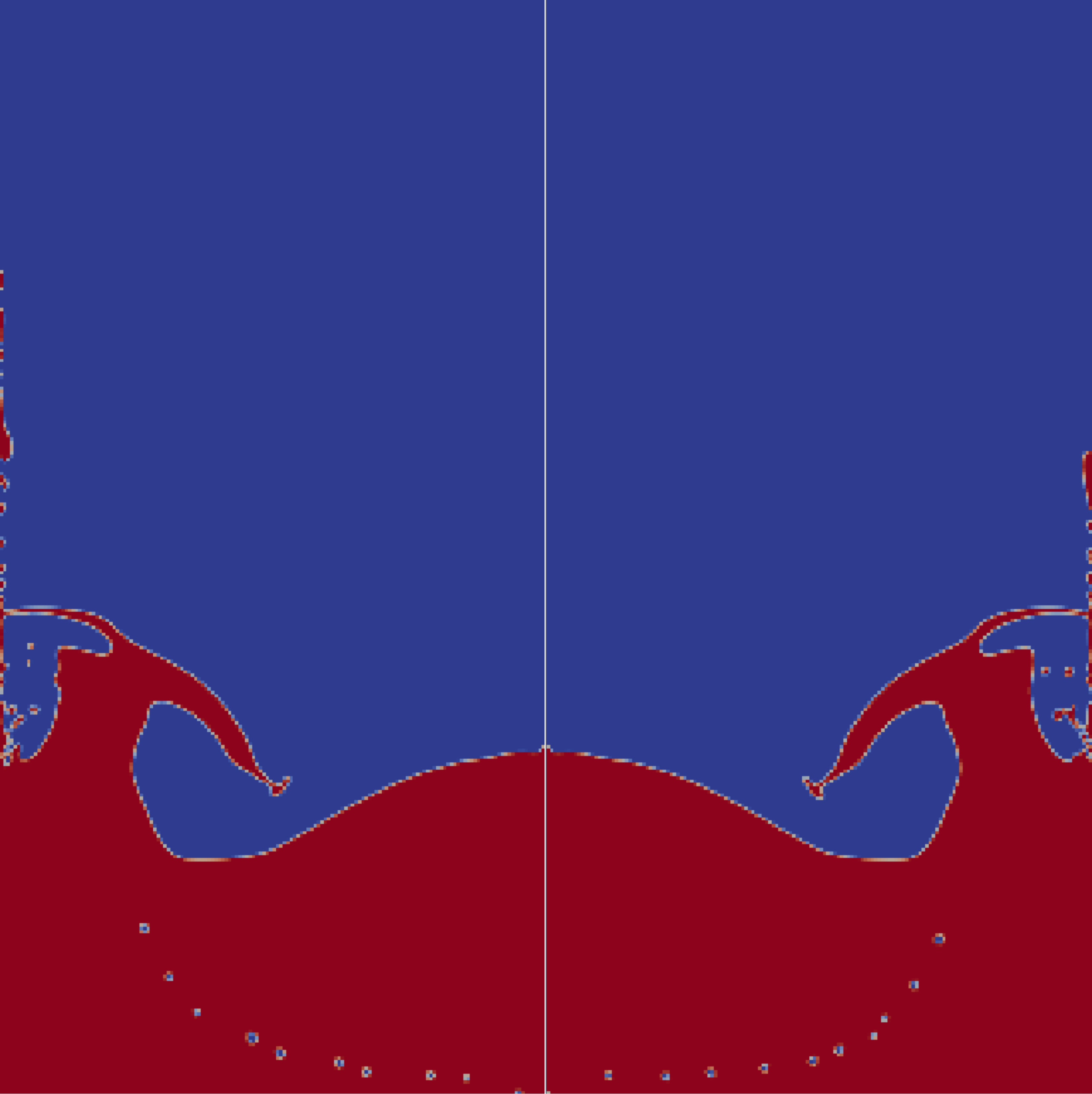}
\end{center}\caption{Volume fraction at $t=0$, $3.5$, $4.5$ and $6$ ms. AD left, DirectCF right.}
\label{fig:Impact_init}
\end{figure}

\noindent At the initial state, the water wall is at rest and the drop is propagating at $1$ km.s$^{-1}$ in the direction of the wall, see Figure~\ref{fig:Impact_init}, $t=0$ ms.

\vspace{0.7cm}

\begin{figure}[h]
\begin{center}
\includegraphics[width=0.8\textwidth]{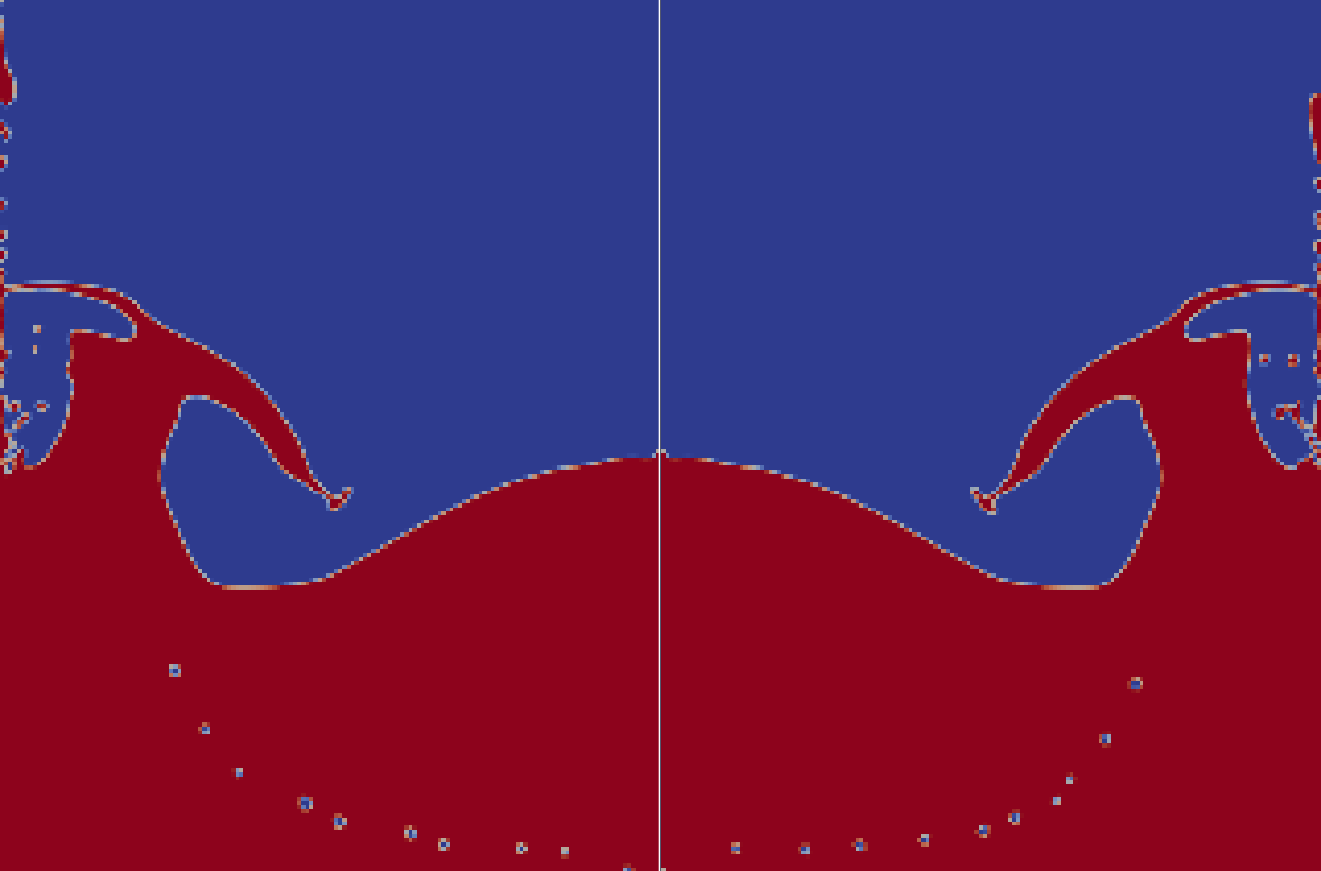}
\end{center}
\caption{Zoom, volume fraction at final state, $t=6$ ms. AD left, DirectCF right.}
\label{fig:Impact_zoom}
\end{figure}

\noindent The drop splashes onto to wall, which creates a jet, see Figure~\ref{fig:Impact_init} at $t=3.5$ ms. At this step, results given by both schemes are very similar. The only hardly visible difference is the shape of the extremity of the jet. Then, at $t=4.5$ ms, craters have appeared on both sides of the drop, and the jets have reached lateral walls (boundary conditions). The second splashing phenomenon shows some differences between the AD and DirectCF remaps, which are more noticeable at $t=6$ ms, see Figure~\ref{fig:Impact_zoom}. The shock wave has propagated up to the boundaries and this creates jets in the upper direction. Droplets in those jets seem to slide higher along lateral walls in the AD case than in the DirectCF case.

Once again, the differences between the two remaps are quite weak compared to the severity of this benchmark. Besides, succeeding in simulating this test is a proof of good robustness. Indeed, at the time of the impact, air is enclosed between two water volumes, which generates extremely high pressure bubbles of air trapped into water, that can be seen on Figure~\ref{fig:Impact_zoom}. There is also high rarefaction of water in the jet, which could create negative pressures of water in some cells, and make the simulation crash.

\clearpage

 \newpage
 \strut
 \newpage

\section*{Conclusion}
\addcontentsline{toc}{section}{Conclusion}	
	
\hspace{1.5em}This internship aimed at proposing and testing an alternative remap to the AD remap, more adapted to exascale computers constraints. The Direct remap with Corner Fluxes has been set up theoretically. This remap enables to catch corner effects even better than the AD remap while performing the remap phase in one step. Those two features are crucial since they impact respectively the accuracy and scalability of the scheme. The remap has also been extended to multi-material flows with sharp interface reconstruction. In SHY code, the Direct and DirectCF remaps have been implemented, first in the mono-material case. Then the code has been enhanced with two different multi-fluid models: a mixing model and a model of non-miscible fluids with sharp interface reconstruction. \\
Numerical results obtained show that the Direct remap is not a good option because of its lack of accuracy, especially in multi-material simulations. They also show that the AD and DirectCF remaps are comparable, in terms of robustness and accuracy, even on severe benchmarks. What follows naturally from this work is the parallelization of the code. Indeed, since the goal was to make the remap scalable in order to be more efficient when running in parallel, the next step will aim at testing this. 

Regarding the scheme itself, many perspectives can be considered so far. The main short-term one is probably to fix problems of second order reconstruction of variables at the corners. Direct remapping enables multidimensional reconstruction. So one can imagine a polynomial interpolation on the 9 points stencil of the DirectCF remap. \\
One could also take into account the rotation of the interface during the Lagrangian phase by computing the vorticity, i.e.\ rot $\mathbf{u}$, at the centres of the cells. Another interesting perspective is the implementation of this remap in a cell-centred Lagrangian scheme (i.e.\ with velocities at the centres of the cells) because first, those schemes are strictly conservative in mass, momentum and total energy, and second, the fact that all variables are centred means that they will be remapped exactly the same way. Thus, one could plan to vectorize the remap in order to run on vector processors. Besides, in the context of multi-material flows simulation, the mixed-cell model could be more advanced, by considering one velocity per material and solving Riemann problems at the interface, see \cite{blais2013dealing}. But in order to do so, it is mandatory for the Lagrangian scheme to be cell-centred. \\

  \clearpage
  \newpage
  \strut
  \newpage
  \appendix
\section{Annexe: \'Etude de la diffusion de la vorticit\'e en hydrodynamique par diff\'erents types de sch\'emas eul\'eriens}
\addtocontents{toc}{\protect\setcounter{tocdepth}{1}}

\vspace{1.5cm}

\subsection*{Introduction}

\hspace{1.5em}Le but de cette \'etude est de comparer trois diff\'erents types de sch\'emas num\'eriques eul\'eriens appliqu\'es \`a l'hydrodynamique. Les trois appartiennent \`a la famille des sch\'emas Lagrange+Projection. On rappelle que ces sch\'emas pr\'esentent deux phases successives pour chaque pas de temps : une phase de d\'eformation du maillage fixe eul\'erien suivant le mouvement lagrangien du fluide (phase lagrangienne), puis une phase de projection des quantit\'es lagrangiennes sur le maillage fixe initial. Nous distinguerons donc trois instants : l'instant initial $t^n$, l'instant correspondant \`a la fin de la phase lagrangienne $t^{lag}$, et celui marquant la fin de la phase de projection $t^{n+1}$, qui est aussi l'instant final de l'it\'eration $n$.

Plus pr\'ecis\'ement, on s'int\'eresse au traitement de la vorticit\'e sur maillage carr\'e. La m\'ethode mise en \oe uvre est d'\'etudier le comportement de plusieurs types de vortex lorsqu'on les laisse \'evoluer sur un pas temps, pour les sch\'emas BBC \cite{saas2014bbc}, MYR (type VNR \cite{vonneumann1950method}) et GLACE \cite{carre2009cell} (ou EUCCLHYD \cite{maire2011contribution}). Pour cela, on cherche \`a obtenir des expressions analytiques du champ de vitesse et des grandeurs thermodynamiques associ\'ees \`a chaque type de vortex, qu'on projette ensuite sur le maillage. Apr\`es un pas de temps, chaque sch\'ema donnera un nouvel ensemble de valeurs discr\`etes pour la vitesse, desquelles on d\'eduira les caract\'eristiques du nouveau vortex. Il s'agit finalement de comparer les comportements des sch\'emas en pr\'esence d'un vortex pur.

\clearpage

\subsection{Cadre th\'eorique}

\hspace{1.5em}On se place dans le cadre de l'hydrodynamique compressible en monomat\'eriau. En n\'egligeant la viscosit\'e du fluide, nous disposons des \'equations d'Euler, dont on rappelle l'expression ci-dessous :

\begin{equation}
\label{euler}
\left \{ \begin{array}{clcr}
\rho \dfrac{d}{dt} \left( \dfrac{1}{\rho} \right) - \mbox{div} (\mathbf{u}) & = & 0 \\
\: \\
\rho \dfrac{d}{dt} \mathbf{u} + \boldsymbol{\nabla}\!P & = & 0 \\
\: \\
\rho \dfrac{d}{dt} E + \mbox{div} (\rho u) & = & 0
\end{array}
\right .
\end{equation}

$\:$ \\
On ajoute la loi d'\'etat des gaz parfaits (EOS) :
\begin{displaymath}
P = (\gamma-1) \rho e
\end{displaymath}

\subsubsection{Principe}

\hspace{1.5em}On cherche une solution analytique de \eqref{euler} sous la forme d'un "vortex pur", en r\'egime stationnaire. En d'autres termes, on impose une vitesse de type $ \mathbf{u}(t,r,\theta) = u_{\theta}(r) \mathbf{e}_{\theta} $ en polaires (on se place dans le cas 2D). Cette condition implique que div$(\mathbf{u})=0$, ce qui nous place dans le cas incompressible.

On peut donc \'ecrire qu'il existe un vecteur $\boldsymbol{\Psi}$ appel\'e vecteur fonction de courant tel que $\mathbf{u} =$ \textbf{rot}$(\boldsymbol{\Psi})$. On a alors :
\begin{displaymath}
\mbox{\textbf{rot rot}}(\boldsymbol{\Psi}) = \boldsymbol{\nabla}(\mbox{div}(\boldsymbol{\Psi}))- \Delta \boldsymbol{\Psi} = \mbox{\textbf{rot}}(\mathbf{u})
\end{displaymath}
En 2D, on a $\mbox{\textbf{rot}}(\mathbf{u}) = \mbox{rot}(\mathbf{u})\mathbf{e}_z$. Si de plus on choisit $\boldsymbol{\Psi}$ \`a divergence nulle, et qu'on pose $\boldsymbol{\Psi} = \psi_z \: \mathbf{e}_z$, il vient :
\begin{equation}
\label{eqn_rot}
\mbox{rot}(\mathbf{u}) = - \Delta \psi_z
\end{equation}
Autrement dit, \'etant donn\'e un champ de vorticit\'e $\mbox{rot}(\mathbf{u})$, $\psi_z$ est solution de l'\'equation de Poisson.

Dans le cadre de notre probl\`eme, comme $\mathbf{u}$ ne d\'epend que de $r$, $\psi_z$ aussi. L'\'equation \eqref{eqn_rot} devient :
 \begin{equation*}
\frac{1}{r}\frac{\partial}{\partial r}\big(r u_{\theta}\big) = - \frac{1}{r}\frac{\partial}{\partial r}\left( r \frac{\partial \psi_z}{\partial r}\right)
\end{equation*}
 \begin{equation*}
\mbox{d'o\`u} \:\:\:\: \: \boxed{u_{\theta} = - \frac{\partial \psi_z}{\partial r}}
\end{equation*}
Ainsi, si on fixe un type de vortex, i.e un champ de vorticit\'e $\mbox{rot}(\mathbf{u})$, on peut retrouver le $u_{\theta}$ associ\'e par l'interm\'ediaire de $\psi_z$, obtenu en r\'esolvant \eqref{eqn_rot}. Les expressions des autres grandeurs du probl\`eme s'obtiennent en injectant $u_{\theta}$ dans le syst\`eme \eqref{euler}.

\subsubsection{Rotationnel num\'erique}

\hspace{1.5em}Dans la mesure o\`u on s'int\'eresse \`a la vorticit\'e, il s'agit de pouvoir la mesurer \`a partir de valeurs num\'eriques. Nous devons donc \`a ce stade d\'efinir un rotationnel discret. On sait d'apr\`es la formule de Stokes que :
\begin{displaymath}
\int \!\!\!\! \int _S \! \mbox{rot} \mathbf{u} \: \mathrm{d}S = \oint_{\partial S} \!\!\! \mathbf{u}\cdot\mathbf{dl}
\end{displaymath}
On veut pouvoir d\'efinir $\mbox{rot} \mathbf{u}$ en un point du maillage. Nous avons a priori 3 lieux possibles : aux centres des mailles primales, aux n\oe uds (centres des mailles duales), ou aux milieux des faces. Usuellement, le contour choisi pour le calcul de la circulation a la m\^eme taille et la m\^eme forme qu'une maille (dans notre cas un carr\'e de c\^ot\'e $\Delta x$), ce qui permet de donner une valeur au centre de ce dernier. Mais pour cela il faut \^etre capable de d\'efinir les vitesses sur les faces du contour. De fait, les vitesses n'\'etant pas localis\'ees aux m\^emes endroits pour les trois sch\'emas, si on veut que notre \'etude comparative ait un sens, il va falloir distinguer les cas pour que le choix du contour soit coh\'erent avec le sch\'ema \'etudi\'e. En termes plus math\'ematiques, on va chercher \`a d\'efinir $\mbox{rot}(\mathbf{u})$ aux centres des mailles d'un maillage conforme dans $H(\mbox{rot})$.

$\:$ \\
\textbf{Sch\'emas d\'ecal\'es, vitesses aux n\oe uds (MYR) :} $\partial S$ contour d'une maille primale

\begin{figure}[h]
\centering
   \includegraphics[width=6cm]{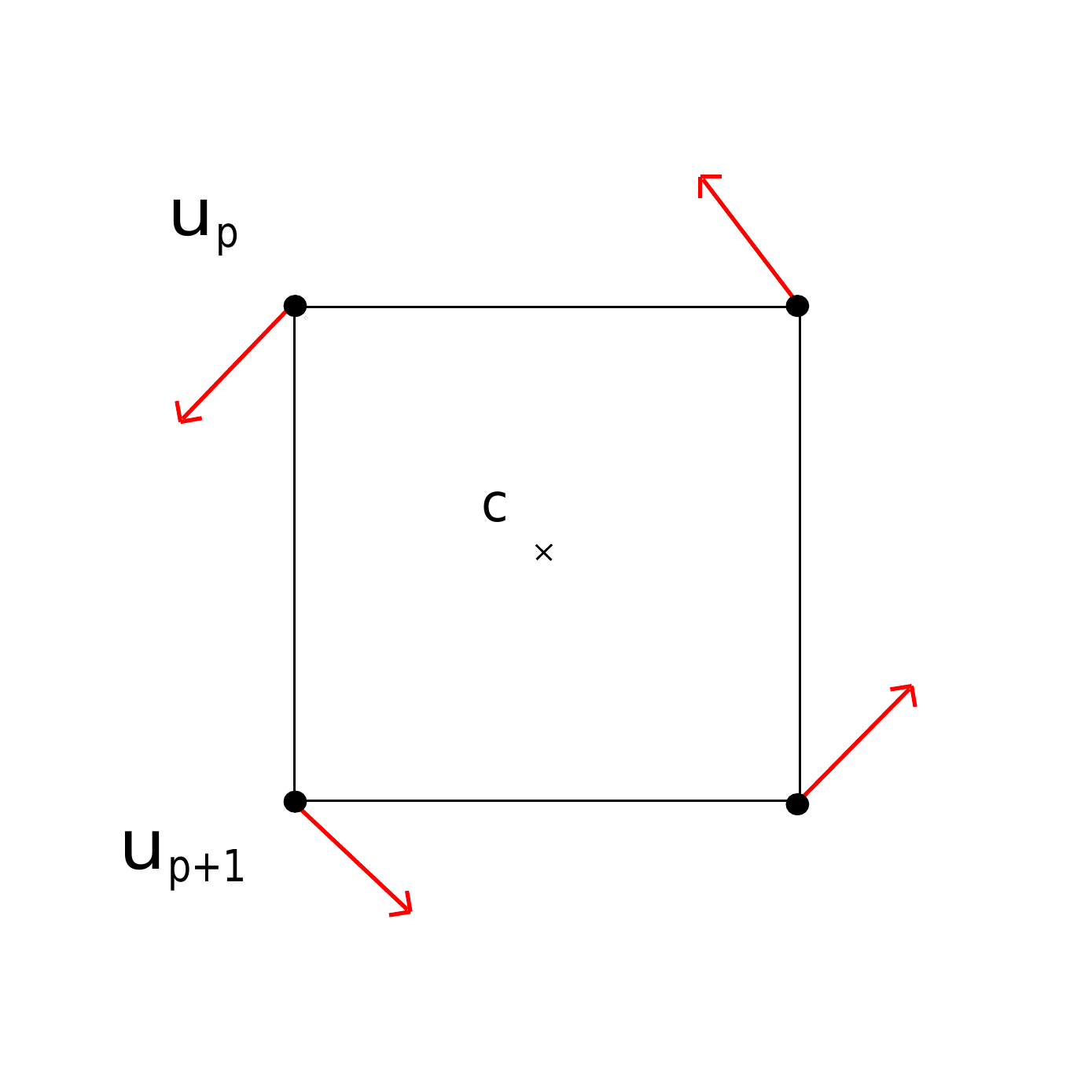}
   \includegraphics[width=6cm]{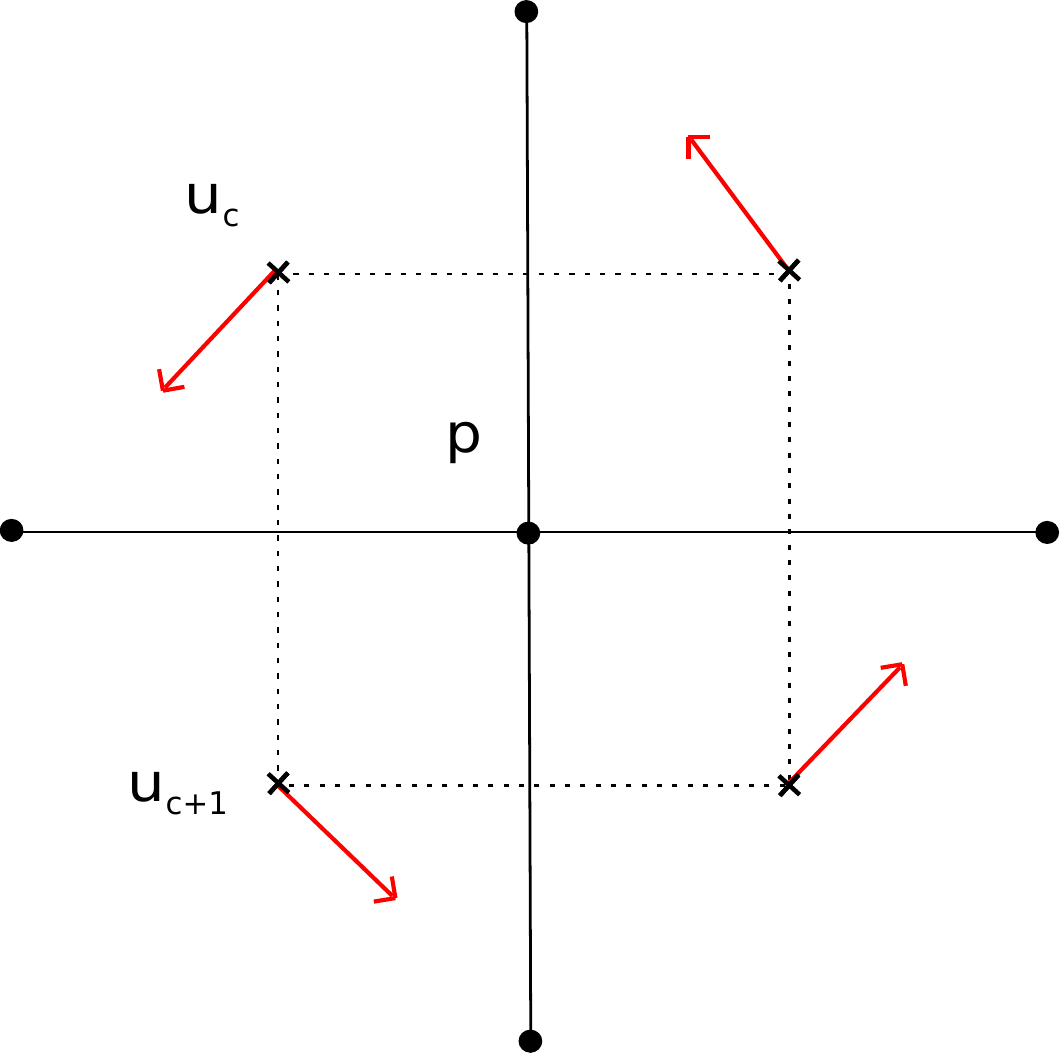}
   \includegraphics[width=6cm]{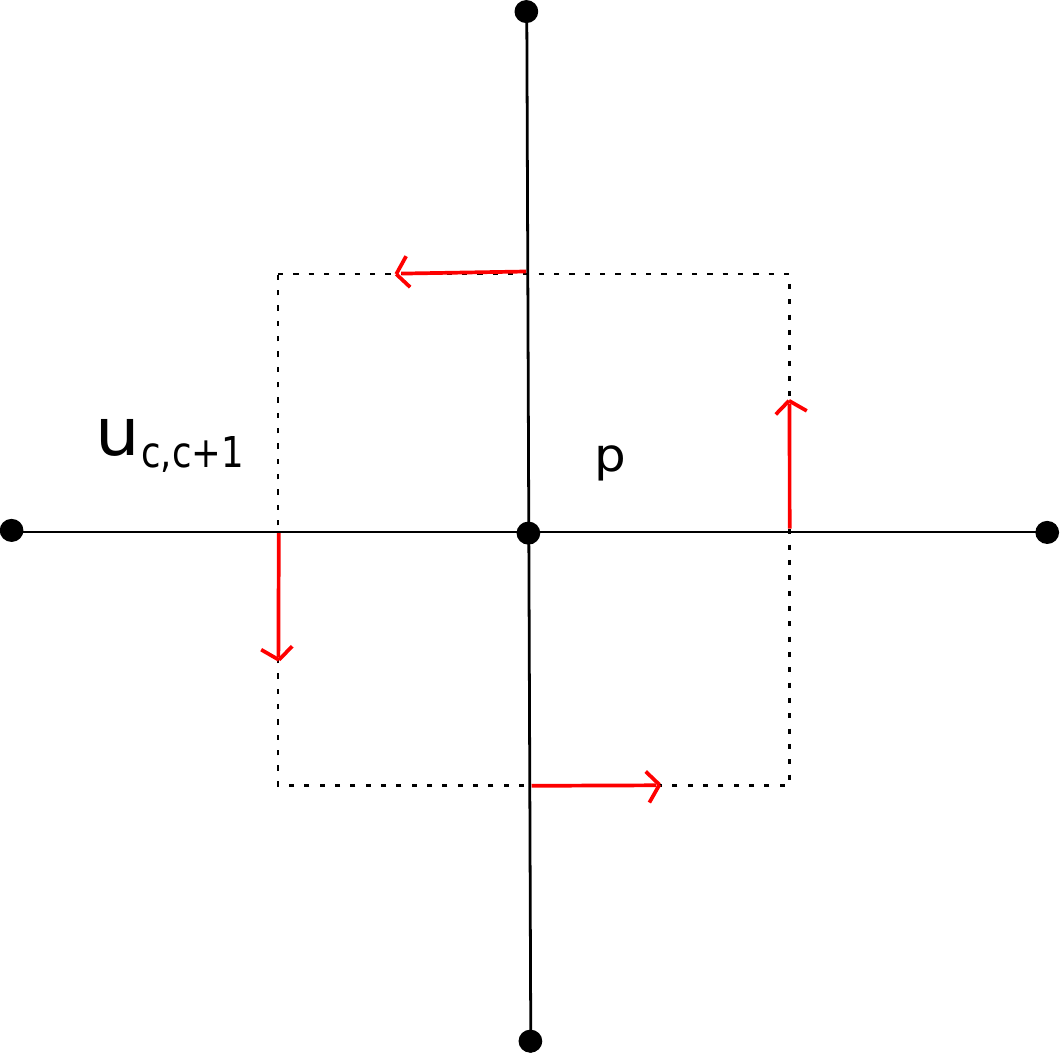} 
   \caption{Sch\'ema de la maille sur laquelle on d\'efinit $\mbox{rot}(\mathbf{u})$ pour les sch\'emas d\'ecal\'e, centr\'e et BBC (de gauche \`a droite et de haut en bas).}  
   \label{fig:schema_rot}
\end{figure}

Pour ce type de sch\'emas, les vitesses sont aux n\oe uds du maillage primal. On peut montrer qu'il s'agit d'un bon maillage pour d\'efinir le rotationnel puisqu'il est conforme dans $H^1\! \supset \!H(\mbox{rot})$. D'un point de vue pratique, en reconstruisant lin\'eairement la vitesse sur une face (comme on le fait pour calculer les flux dans la projection de MYR), on est capable de donner la valeur de l'int\'egrale de $\mathbf{u}\cdot\mathbf{dl}$ sur cette derni\`ere. Si on indice sur $p$ les sommets de la maille $c$ au centre de laquelle on cherche \`a calculer le rotationnel, il vient :
\begin{displaymath}
\mbox{rot} (\mathbf{u})_{c,D} = \frac{1}{2 \Delta x} \sum_{p \in \{c\}} \big(\mathbf{u}_{p}+\mathbf{u}_{p+1}\big)\cdot\mathbf{T}_{p,p+1}
\end{displaymath}
o\`u $\mathbf{T}_{p,p+1}$ repr\'esente le vecteur tangent \`a la face $[p,p+1]$, dirig\'e dans le sens direct (voir Figure~\ref{fig:schema_rot}).

$\:$ \\
\textbf{Sch\'emas centr\'es, vitesses aux centres des mailles (GLACE):} $\partial S$ contour d'une maille duale 

Ces sch\'emas ont les vitesses localis\'ees aux centres des mailles primales, qui correspondent aux n\oe uds du maillage dual. Ici, \`a l'inverse, le maillage primal n'est conforme que dans $L^2$, et c'est le maillage dual, conforme dans $H^1\! \supset \!H(\mbox{rot})$, qui sera adapt\'e \`a l'\'ecriture d'un rotationnel discret. On proc\`ede donc de la m\^eme mani\`ere que pour les sch\'emas d\'ecal\'es, mais sur le maillage dual. Cette fois la vitesse est constante par morceaux sur une face de maille duale. Si on l'int\`egre, on obtient toutefois le m\^eme r\'esultat que pour une reconstruction lin\'eaire, en indi\c ant sur $c$ les mailles primales qui partagent le n\oe ud commun $p$, auquel on \'evalue la valeur de $\mbox{rot} (\mathbf{u})$ en $p$.
\begin{displaymath}
\mbox{rot} (\mathbf{u})_{p,C} = \frac{1}{2 \Delta x} \sum_{c \in \{p\}} \big(\mathbf{u}_{c}+\mathbf{u}_{c+1}\big)\cdot\mathbf{T}_{c,c+1}
\end{displaymath}
o\`u $\mathbf{T}_{c,c+1}$ repr\'esente le vecteur tangent \`a la face duale reliant les centres des mailles $c$ et $c+1$, dirig\'e dans le sens direct (voir Figure~\ref{fig:schema_rot}).

$\:$ \\
\textbf{Sch\'emas d\'ecal\'es, vitesses aux faces (BBC) :} $\partial S$ contour d'une maille duale

Dans ce dernier cas, les vitesses sont d\'efinies aux milieux des faces du maillage primal, et dirig\'ees selon les normales \`a ces faces. On peut montrer dans ce cas que le maillage primal est conforme dans $H(\mbox{div})$, et que le maillage dual est lui conforme dans $H(\mbox{rot})$. En pratique, si on choisissait une maille primale comme contour, on obtiendrait un r\'esultat nul (vitesses orthogonales aux vecteurs tangents aux faces). Il convient donc l\`a encore de placer les valeurs du rotationnel aux n\oe uds du maillage primal. On conna\^it alors la vitesse tangentielle sur chaque face du contour, ce qui permet d'\'ecrire :
\begin{displaymath}
\mbox{rot} (\mathbf{u})_{p,BBC} = \frac{1}{\Delta x} \sum_{c \in \{p\}} \mathbf{u}_{c,c+1}\cdot\mathbf{T}_{c,c+1}
\end{displaymath}
o\`u $\mathbf{u}_{c,c+1}$ correspond \`a la vitesse de la face s\'eparant les mailles $c$ et $c+1$ (voir Figure~\ref{fig:schema_rot}).
$\:$ \\

Il reste maintenant \`a faire \'evoluer la solution sur un pas de temps (entre $t^n$ et $t^{n+1}\!=\!t^n\!+\!\Delta t$) pour chaque type de sch\'emas. On calculera ensuite le rotationnel num\'erique \`a l'origine dans chaque cas, en prenant soin de placer l'origine sur un n\oe ud dans les cas centr\'e et BBC, et au centre d'une maille dans le cas d\'ecal\'e. De plus, comme on s'int\'eresse au caract\`ere diffusif, on calculera \'egalement le rapport entre les rotationnels \'e $t^n$ et $t^{n+1}$, pour s'affranchir des \'eventuelles diff\'erences dues \`a la discr\'etisation par les sch\'emas (par exemple il se peut qu'on ait $\mbox{rot} (\mathbf{u})_{p,C}^n \neq \mbox{rot} (\mathbf{u})_{p,BBC}^n$, voir 2.4). On comparera ce rapport \`a 1, et on fera appara\^itre un coefficient de diffusion num\'erique.

\subsection{Cas d'un vortex ponctuel}

\hspace{1.5em}La repr\'esentation qui semble la plus naturelle pour le champ de vorticit\'e cr\'e\'e par un vortex ponctuel plac\'e \`a l'origine est la masse de Dirac. On note $\sigma_0$ le moment cin\'etique introduit par le vortex.
\begin{displaymath}
\mbox{rot} (\mathbf{u}) = \sigma_0 \delta_0
\end{displaymath}
 Il s'agit bien d'un moment car en dimension 2, $[\delta_0(r)]=m^{-2}$ et $[\mbox{rot} (\mathbf{u})]=s^{-1}$, donc $\sigma_0$ s'exprime en $m^2.s^{-1}$. Un vortex ponctuel correspond en fait \`a un vortex qui conserve le moment cin\'etique.

\subsubsection{Calcul de la solution analytique}

\hspace{1.5em}D'apr\`es le paragraphe 1.1, pour trouver l'expression de la fonction de courant, on est donc amen\'e \`a r\'esoudre l'\'equation :
\begin{equation*}
- \Delta \psi_z = \sigma_0 \delta_0 = \mbox{rot} (\mathbf{u})
\end{equation*}
La solution n'est autre que la fonction de Green du Laplacien, i.e en 2D
\begin{equation*}
\psi_z(r) = - \frac{\sigma_0}{2 \pi}\ln r 
\end{equation*}
On obtient donc finalement :
\begin{equation*}
u_{\theta}(r) = - \frac{\partial \psi_z}{\partial r} = \frac{\sigma_0}{2 \pi r}
\end{equation*}

On va montrer qu'il existe une densit\'e $\rho$ (et une pression $P$ associ\'ee) telle que $(\mathbf{u},\rho,P)$ soit solution du syst\`eme d'Euler.
\\On a d\'ej\`a vu en 1.1 que div$(\mathbf{u}) = 0$ . On s'int\'eresse \`a la deuxi\`eme \'equation de \eqref{euler}, qui donne elle :
\begin{displaymath}
-\frac{{u_{\theta}}^2}{r} = -\frac{1}{\rho} \frac{\partial P}{\partial r}
\end{displaymath}
\begin{displaymath}
\mbox{d'o\`u} \:\:\: \rho (r) = \rho_0 \exp \left( -\frac{K_0}{r^2} \right), \:\:\:  \mbox{avec}  \:\:\: K_0 = \frac{\sigma_0^2}{8\pi^2(\gamma-1)e_0}
\end{displaymath}
On voit que $\rho(r) \rightarrow \rho_0$ lorsque $r \rightarrow +\infty$.

$\:$ \\
On d\'efinit $\alpha_0$ tel que $u_{\theta}\big(\frac{\Delta x}{\sqrt{2}}\big)=\sqrt{2} \: \alpha_0 \Delta x$. $\alpha_0$ est une constante d\'ependante du pas du maillage et de $\sigma_0$, qu'on introduit pour simplifier les calculs.
\begin{displaymath}
\boldsymbol{\alpha_0=\frac{\sigma_0}{2\pi {\Delta x}^2}}
\end{displaymath}
Elle peut \^etre vue comme une vitesse angulaire.
\\On peut consid\'erer que $\rho \simeq \rho_0$ si $\frac{K_0}{R^2} \ll 1$, o\`u $R$ repr\'esente le "rayon" d'une maille. Cette condition implique que :
\begin{equation*}
\boxed{
\Delta x \gg \frac{\sigma_0}{2 \pi \sqrt{(\gamma-1)e_0}} = \frac{\sigma_0}{2\pi c_0}
}
\end{equation*}
Ainsi, pour une valeur de $\sigma_0$ suffisamment petite devant la vitesse du son, et un pas de maillage grand devant ce rapport, on peut consid\'erer que sur notre maillage, la solution analytique de \eqref{euler} associ\'ee \`a un vortex ponctuel est assimilable au triplet $\big(\frac{\sigma_0}{2\pi r}\mathbf{e}_{\theta},\rho_0,(\gamma-1)\rho_0e_0\big)$.

\subsubsection{Etude pour les sch\'emas d\'ecal\'es}

\hspace{1.5em}Le mod\'ele de sch\'ema utilis\'e est le sch\'ema MYR d'ordre 1 en espace. Il pr\'esente la m\^eme discr\'etisation en espace que VNR, mais la discr\'etisation en temps de type saute-mouton a \'ete remplac\'e par un sch\'ema pr\'edicteur-correcteur. Comme on l'a vu, les vitesses sont localis\'ees aux n\oe uds, et le reste des grandeurs aux centres des mailles. Le centre du vortex est plac\'e au centre d'une maille. 
$\:$ \\
\begin{figure}[h]
\centering
   \includegraphics[width=7cm]{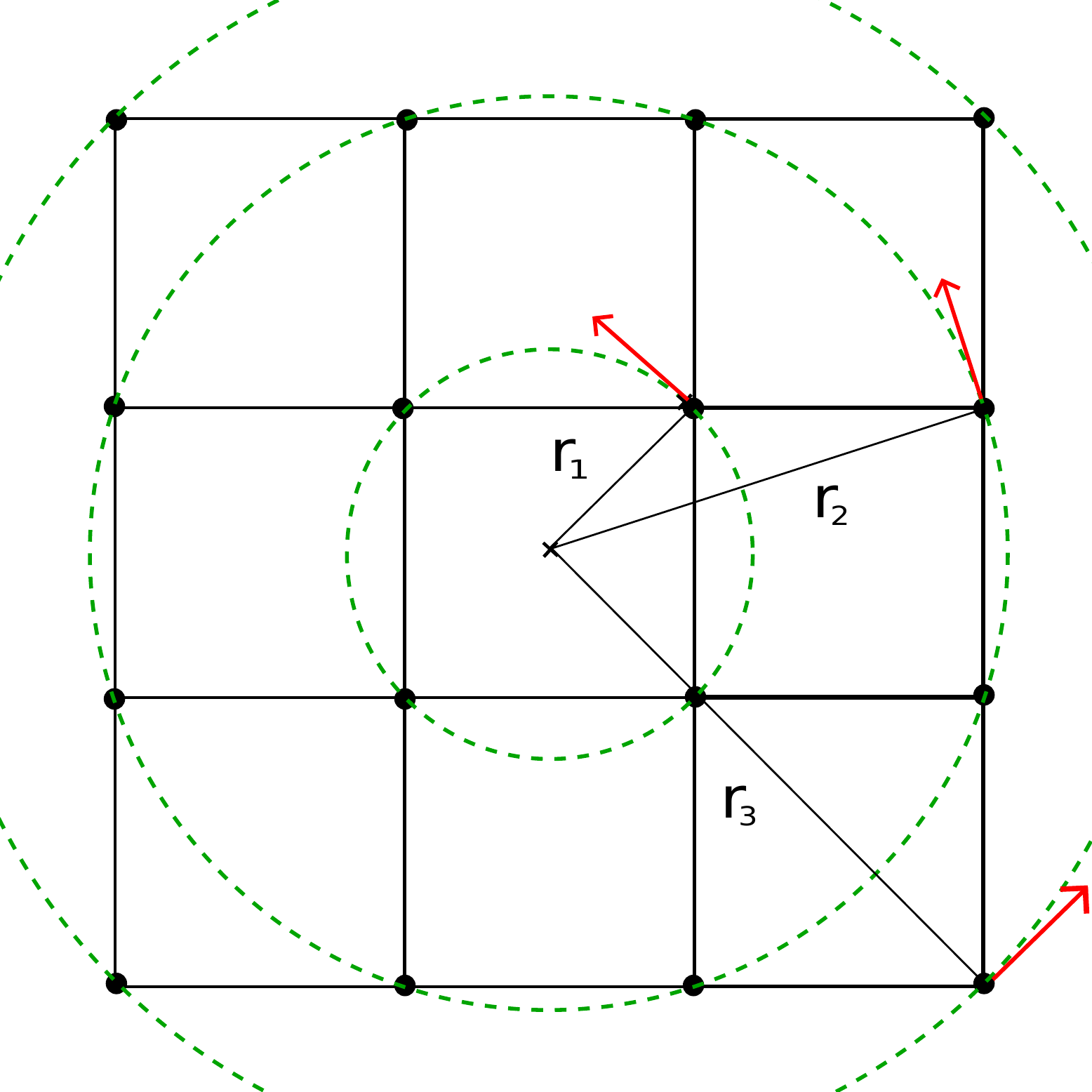}
   \caption{Sch\'ema du champ de vitesses sur le maillage dans le cas d\'ecal\'e.}   			    	   \label{fig:schemadecale}
\end{figure}

\paragraph{D\'etermination de $\mathbf{u}$ sur le maillage}
$\:$

On commence par calculer les vitesses aux n\oe uds \`a $t^n$ (obtenues \`a partir de la solution analytique), dont nous avons besoin pour pouvoir faire \'evoluer le maillage durant la phase lagrangienne.
\begin{displaymath}
r_1 = \frac{\Delta x}{\sqrt{2}},  \:\:\:\: r_2 = \sqrt{\frac{5}{2}}\Delta x,  \:\:\:\:  r_3 = 3\frac{\Delta x}{\sqrt{2}}
\end{displaymath}
Et donc, pour les normes des vitesses associ\'ees \`a ces rayons :
\begin{displaymath}
 u_1 = \sqrt{2}\alpha_0\Delta x,  \:\:\:\: u_2 = \sqrt{\frac{2}{5}}\alpha_0\Delta x,  \:\:\:\:  u_3 = \frac{\sqrt{2}}{3}\alpha_0\Delta x
\end{displaymath}
Il ne reste qu'\`a multiplier par le $\mathbf{e_{\theta}}$ correspondant pour obtenir les vecteurs vitesse.

\paragraph{Phase lagrangienne}
$\:$

\begin{figure}[h]
\centering
   \includegraphics[width=6cm]{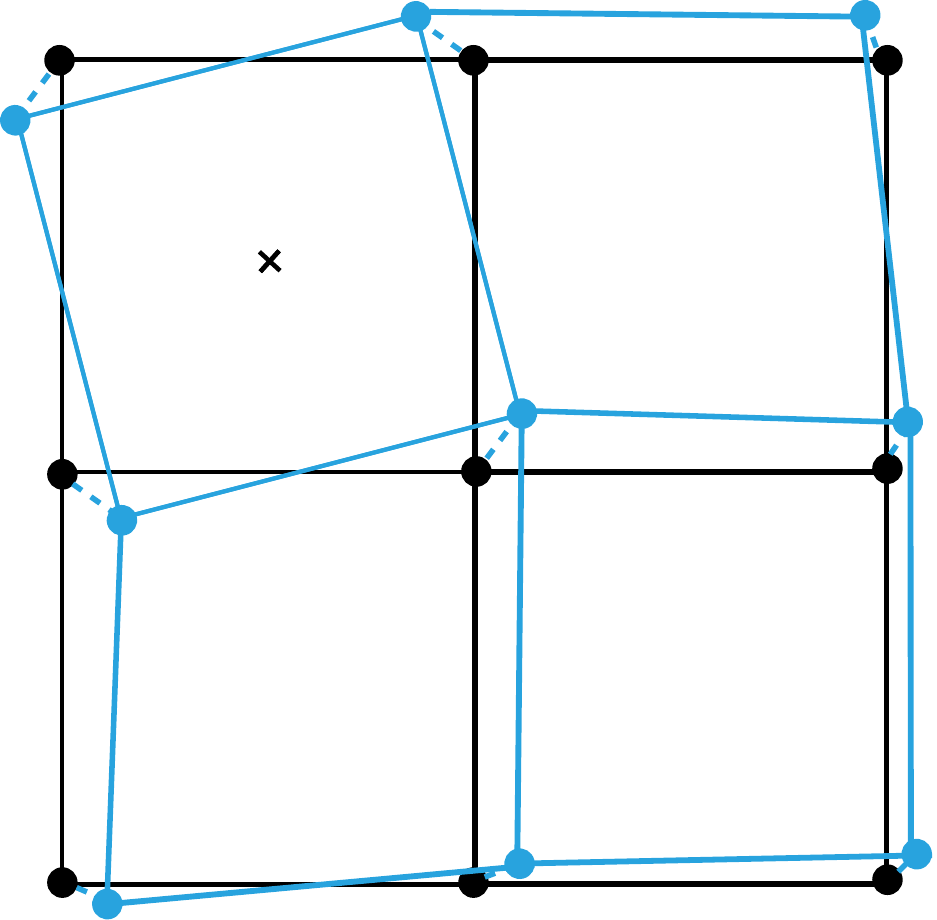}
   \caption{Sch\'ema de la d\'eformation lagrangienne dans le cas d\'ecal\'e pour un vortex ponctuel.}   
   \label{fig:deformation_ponctuel_decale}
\end{figure}

 Puisque les vitesses nodales sont connues \`a $t^n$, on conna\^it directement la d\'eformation lagrangienne du maillage, voir Figure~\ref{fig:deformation_ponctuel_decale}. On garde les notations utilis\'ees dans \cite{courtesinglard2012schemas}. Les cellules sont indic\'ees sur $c$ et les n\oe uds sur $p$. On a, au n\oe ud $p$ :
\begin{displaymath}
\mathbf{u}_p^{lag} = \mathbf{u}_p^{n} + \Delta t \: \boldsymbol{\gamma}_p^{lag}
\end{displaymath}
\begin{displaymath}
\mbox{avec} \:\:\:\: \boldsymbol{\gamma}_p^{lag} = \frac{1}{m_p} \sum_{c \in \{p\}} \! \big( P_c^{lag} + Q_c^{n}\big) L_{cp}^{lag} \mathbf{N}_{cp}^{lag} = 0 
\end{displaymath}
Cela vient du fait que nous ayons un champ de pression uniforme $P_c^{lag} = P_0$, $\forall c$, et que  $ \sum L_{cp}^{lag} \mathbf{N}_{cp}^{lag} = 0$, ce qui traduit la conservation de volume de chaque maille. Cette conservation de volume justifie au passage la nullit\'e du terme de pseudo-viscosit\'e $Q_c^{n}$.
\\ Il vient donc, pour tout noeud $p$, 
\begin{displaymath}
\boxed{
\mathbf{u}_p^{lag} = \mathbf{u}_p^{n}
}
\end{displaymath}

\paragraph{Phase de projection}
$\:$

La projection effectu\'ee est une projection directe. On calcule les flux de masse aux faces de la maille nodale $\phi_p^{\alpha}$, ($\alpha = E,O,N,S$), en faisant la moyenne entre les flux aux faces des 4 mailles voisines du n\oe ud. Ensuite, on \'ecrit la conservation de la quantit\'e de mouvement dans la maille nodale.
\begin{equation}
\label{dual_remap}
m_p^{proj} \mathbf{u}_p^{proj} = m_p^{lag} \mathbf{u}_p^{lag} + \sum_\alpha \phi_p^\alpha \mathbf{u}_p^{\alpha}
\end{equation}
Dans notre cas, on a $m_p^{proj}=m_p^{lag}=\rho_0 {\Delta x}^2$ puisque le volume et la densit\'e sont conserv\'es.
\\Pour illustrer cette \'etape, prenons le n\oe ud $p$ correspondant au n\oe ud en bas \`a droite (les autres s'en d\'eduisent par sym\'etrie). A l'ordre 1, les vitesses $\mathbf{u}_p^{\alpha}$ sont choisies comme suit :
\begin{displaymath}
\mathbf{u}_p^O = \alpha_0 \Delta x \left( \!\! \begin{array}{clcr}1\\-1 \end{array} \!\!\right), \:\:\: \mathbf{u}_p^S = \frac{\alpha_0 \Delta x}{5} \left( \!\! \begin{array}{clcr}3\\1 \end{array} \!\!\right), \:\:\: \mathbf{u}_p^E = \frac{\alpha_0 \Delta x}{5} \left( \!\! \begin{array}{clcr}1\\3 \end{array} \!\!\right), \:\:\: \mathbf{u}_p^N = \alpha_0 \Delta x \left( \!\! \begin{array}{clcr}-1\\1 \end{array} \!\!\right)
\end{displaymath}
Ce qui donne donc, \`a l'ordre 1 :
\begin{displaymath}
\boxed{
\mathbf{u}_{p,1}^{n+1}= \alpha_0 \Delta x \left( \! \begin{array}{clcr}1+ \frac{68}{75}\alpha_0 \Delta t \\ \: \\1- \frac{68}{75}\alpha_0 \Delta t \end{array} \!\right)
}
\end{displaymath}
Pour obtenir une projection d'ordre 2, il suffit de proc\'eder \`a une reconstruction lin\'eaire pour les $\mathbf{u}_p^{\alpha}$. Ce qui revient \`a poser :
\begin{displaymath}
\mathbf{u}_p^O = \alpha_0 \Delta x \left( \!\! \begin{array}{clcr}1\\0 \end{array} \!\!\right), \:\:\: \mathbf{u}_p^S = \frac{\alpha_0 \Delta x}{5} \left( \!\! \begin{array}{clcr}4\\3 \end{array} \!\!\right), \:\:\: \mathbf{u}_p^E = \frac{\alpha_0 \Delta x}{5} \left( \!\! \begin{array}{clcr}3\\4 \end{array} \!\!\right), \:\:\: \mathbf{u}_p^N = \alpha_0 \Delta x \left( \!\! \begin{array}{clcr}0\\1 \end{array} \!\!\right)
\end{displaymath}
On obtient donc, pour l'ordre 2 :
\begin{displaymath}
\boxed{
\mathbf{u}_{p,2}^{n+1}= \alpha_0 \Delta x \left( \! \begin{array}{clcr}1+ \frac{34}{75}\alpha_0 \Delta t\\ \: \\1- \frac{34}{75}\alpha_0 \Delta t \end{array} \!\right)
}
\end{displaymath}

\paragraph{Calcul des rotationnels discrets}
$\:$

D'apr\`es 1.2, comme ici $\frac{1}{2}\big(\mathbf{u}_p^n+\mathbf{u}_{p+1}^n\big)\cdot\mathbf{T}_{p,p+1} = \alpha_0 \Delta x$,
\begin{displaymath}
\mbox{rot} (\mathbf{u})_{c,D}^n = 4 \: \alpha_0 = \frac{\sigma_0}{\pi {\big(\Delta x / \sqrt{2}\big)}^2}
\end{displaymath}
Notons que $\mbox{rot} (\mathbf{u})_d \rightarrow +\infty$ lorsque $\Delta x \rightarrow 0$, ce qui est coh\'erent avec notre choix de repr\'esentation par une masse de Dirac.
\\De plus, \`a $t^{n+1}$, compte tenu des r\'esultats obtenus dans la partie pr\'ec\'edente, on a quel que soit l'ordre :
\begin{displaymath}
\frac{1}{2}\big(\mathbf{u}_p^{n+1}+\mathbf{u}_{p+1}^{n+1}\big)\cdot\mathbf{T}_{p,p+1} = \alpha_0 \Delta x
\end{displaymath} 
\begin{displaymath}
\mbox{et donc,} \:\:\:\: \boxed{\frac{\mbox{rot} (\mathbf{u})_{c,D}^{n+1}}{\mbox{rot} (\mathbf{u})_{c,D}^{n}} = 1} 
\end{displaymath}
$\:$ 

Le sch\'ema d\'ecal\'e ne diffuse donc pas la vorticit\'e des vortex ponctuels, quel que soit l'ordre choisi pour la projection.
\\En revanche, on remarque que $\|\mathbf{u}_p^{n+1}\|>\|\mathbf{u}_p^{n}\|$, ce qui peut \^etre une source d'instabilit\'e. De plus, la direction des vitesses change elle aussi, on a apparition d'une composante strictement positive selon $\mathbf{e_r}$, \'egale \`a $\frac{68}{75}\alpha_0 \Delta t \big(\sqrt{2} \: \alpha_0 \Delta x\big)$ pour la projection d'ordre 1, et \`a $\frac{34}{75}\alpha_0 \Delta t \big(\sqrt{2} \: \alpha_0 \Delta x\big)$ pour l'ordre 2. On observe donc que le vortex a tendance \`a s'\'etirer, et le ph\'enom\'ene est logiquement att\'enu\'e \`a l'ordre 2.

\subsubsection{Etude pour les sch\'emas centr\'es}

\hspace{1.5em}Le mod\'ele de sch\'ema utilis\'e est dans ce cas le sch\'ema GLACE d'ordre 1, voir \cite{carre2009cell}. Toutes les grandeurs sont centr\'ees, et on place le centre du vortex sur un n\oe ud du maillage.
$\:$\\
\begin{figure}[h]
\centering
   \includegraphics[width=7cm]{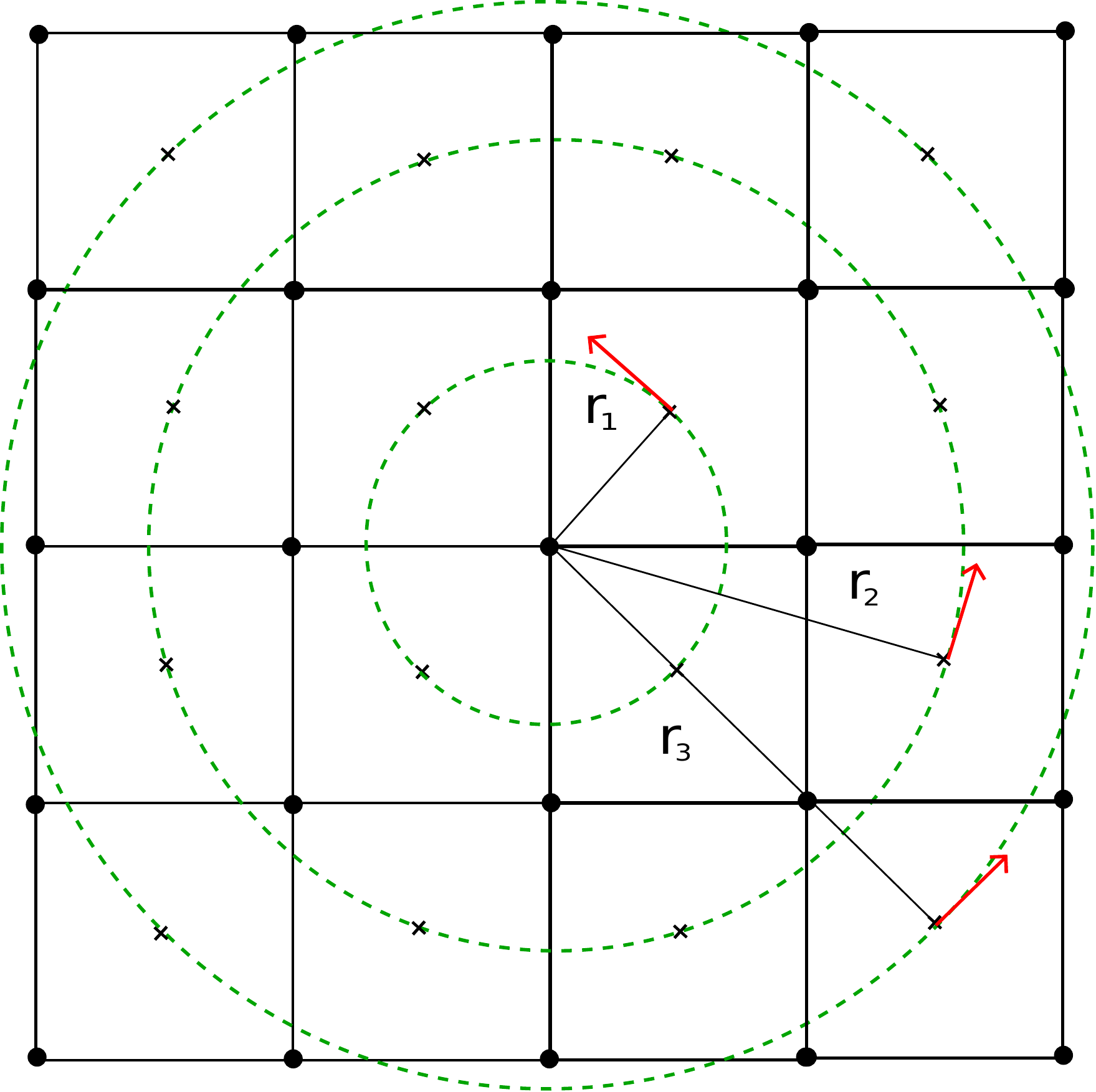}
   \caption{Sch\'ema du champ de vitesses sur le maillage dans le cas centr\'e.} 
   \label{schemacentre}  
\end{figure}

\paragraph{D\'etermination de $\mathbf{u}$ sur le maillage}
$\:$

Commen\c cons par d\'eterminer les valeur des vitesses centr\'ees sur le maillage. On a cette fois :
\begin{displaymath}
r_1 = \frac{\Delta x}{\sqrt{2}},  \:\:\:\:  r_2 = \sqrt{\frac{5}{2}} \Delta x,  \:\:\:\:  r_3 = \frac{3\sqrt{2}}{2} \Delta x
\end{displaymath}
\begin{displaymath}
\mbox{et donc :}\:\:\:\: u_1 = \sqrt{2} \:\alpha_0\Delta x,  \:\:\:\: u_2 = \sqrt{\frac{2}{5}}\alpha_0\Delta x,  \:\:\:\:  u_3 = \frac{\sqrt{2}}{3}\alpha_0\Delta x  
\end{displaymath}
$\:$\\
Comme pr\'ec\'edemment, il suffit ensuite de multiplier par le bon $\mathbf{e_{\theta}}$.

\paragraph{Phase lagrangienne}
$\:$

Contrairement au cas d\'ecal\'e, on ne dispose pas des vitesses aux noeuds, il nous faut donc les calculer pour pouvoir d\'eformer le maillage. Cette \'etape cruciale repose sur la construction d'un solveur aux noeuds, qui permet de d\'eterminer non seulement les $\mathbf{u}_p$, mais aussi les pressions aux demi-faces voisines du noeud $p$. Ce solveur provient de la r\'e\'ecriture des \'equations de conservation qui, pour un champ de pression $P_0$ uniforme donne :
\begin{equation}
\label{nodal_solver}
\sum_{c\in \{ p \}}\left(P_{c-1/2}^{c-1}-P_{c-1/2}^{c}\right)\Delta x \: \mathbf{N}_{c-1}^{c} = 0
\end{equation}
\begin{displaymath}
\left \{ \begin{array}{clcr}
P_{c-1/2}^{c} & = & P_0 + \rho_0 c_0 \big(\mathbf{u}_p-\mathbf{u}_c\big)\cdot\mathbf{N}_{c-1}^{c} \\
P_{c-1/2}^{c-1} & = & P_0 - \rho_0 c_0 \big(\mathbf{u}_p-\mathbf{u}_{c-1}\big)\cdot\mathbf{N}_{c-1}^{c}
\end{array}
\right .
\end{displaymath}

\begin{figure}[h]
\centering
   \includegraphics[width=6cm]{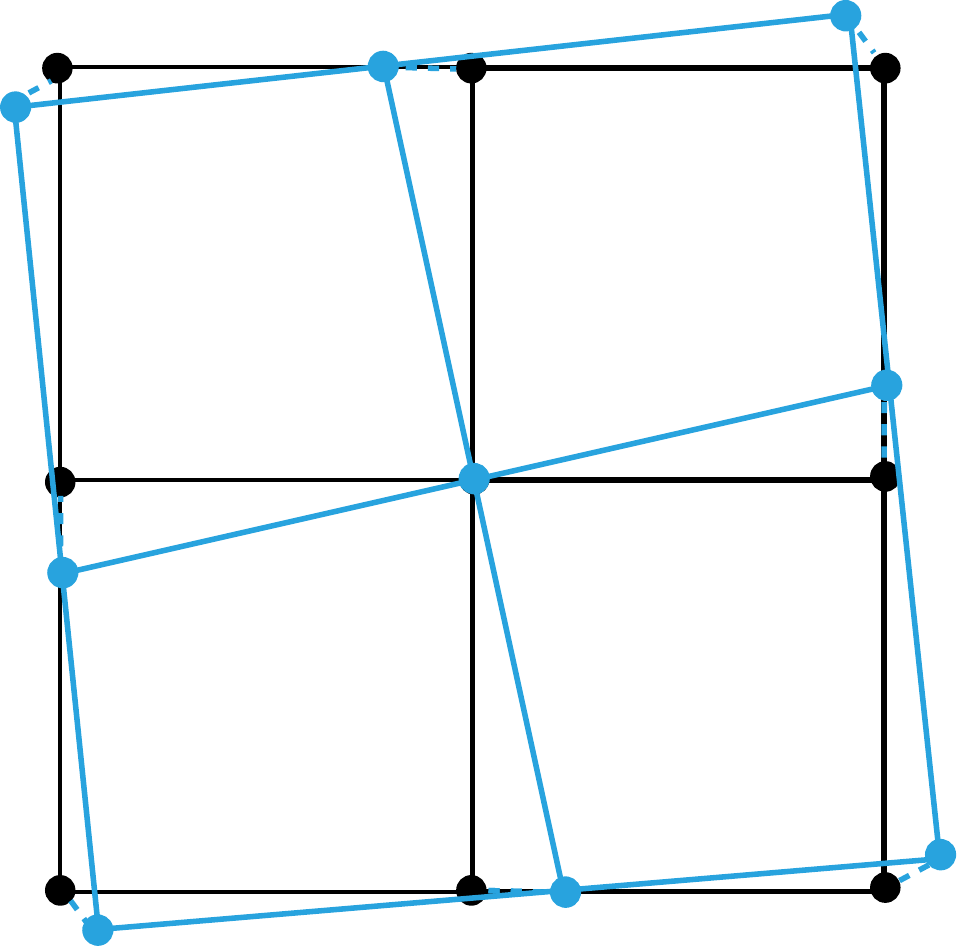}
   \caption{Sch\'ema de la d\'eformation lagrangienne dans le cas centr\'e pour un vortex ponctuel.}
   \label{fig:deformation_ponctuel_centre}   
\end{figure}

Une fois les vitesses aux noeuds, et donc la d\'eformation lagrangienne, connues (voir Figure~\ref{fig:deformation_ponctuel_centre}), on calcule les nouvelles vitesses centr\'ees. Pour ce type de sch\'emas, on obtient lors de la phase lagrangienne :
\begin{equation}
\label{lagrangian_velocity_centred}
\mathbf{u}_c^{lag} = \mathbf{u}_c^n - \frac{\Delta t}{m_c}\sum_{p\in \{ c \}} \frac{1}{2}\left(P_p^{c,p+1}+P_{p+1}^{c,p}\right)\Delta x \: \mathbf{N}_{p,p+1}
\end{equation}
o\`u $P_p^{c,p+1}$ repr\'esente la pression exerc\'ee sur la demi-face de $[p,p+1]$ correspondant au noeud $p$ dans la maille $c$, et $\mathbf{N}_{p,p+1}$ la normale sortante \`a la face $[p,p+1]$.
\\Pour la maille en bas \`a droite, les calculs donnent :
\begin{displaymath}
\boxed{
\mathbf{u}_{c}^{lag}= \alpha_0 \Delta x \left( \! \begin{array}{clcr}1- \frac{4}{3}\frac{\Delta t}{\Delta x}c_0 \\ \: \\1- \frac{4}{3}\frac{\Delta t}{\Delta x}c_0 \end{array} \!\right)
}
\end{displaymath}

\paragraph{Phase de projection}
$\:$

Comme pour le cas d\'ecal\'e, on va calculer les flux de masse alg\'ebriques, mais cette fois ce sont les flux \`a travers les faces de la maille r\'eelle. Puis, comme pour \eqref{nodal_solver}, on \'ecrit l'equation de la conservation de la quantit\'e de mouvement, en utilisant ici un d\'ecentrement upwind pour le calcul des $\mathbf{u}_c^\alpha$.
On obtient finalement pour la maille du bas, \`a l'ordre 1 en $\Delta t$ :
\begin{displaymath}
\boxed{
\mathbf{u}_{c}^{n+1}= \alpha_0 \Delta x \left( \! \begin{array}{clcr}1 - \frac{4}{3}\frac{\Delta t}{\Delta x}c_0 - \frac{8}{75}\alpha_0 \Delta t\\ \: \\1 - \frac{4}{3}\frac{\Delta t}{\Delta x}c_0 - \frac{76}{75}\alpha_0 \Delta t \end{array} \!\right)
}
\end{displaymath}

\paragraph{Calcul des rotationnels discrets}
$\:$

D'apr\`es 1.2, comme ici $\frac{1}{2}\big(\mathbf{u}_c^n+\mathbf{u}_{c+1}^n\big)\cdot\mathbf{T}_{c,c+1} = \alpha_0 \Delta x$,
\begin{displaymath}
\mbox{rot} (\mathbf{u})_{p,C}^n = 4 \: \alpha_0 = \frac{\sigma_0}{\pi {\big(\Delta x / \sqrt{2}\big)}^2}
\end{displaymath}
On remarque qu'on a \'egalit\'e des valeurs du rotationnel num\'erique \`a $t^n$, ce qui confirme le bon choix des contours.
\\Et \`a $t^{n+1}$, \`a pertir des r\'esultats pr\'ec\'edents, il vient :
\begin{displaymath}
\boxed{\frac{\mbox{rot} (\mathbf{u})_{p,C}^{n+1}}{\mbox{rot} (\mathbf{u})_{p,C}^{n}} = 1 - \frac{4}{3}\frac{\Delta t}{\Delta x}c_0 - \frac{14}{25}\alpha_0 \Delta t} 
\end{displaymath}
$\:$ 

Le sch\'ema centr\'e a lui tendance \`a diffuser la vorticit\'e des vortex ponctuels.
\\En ce qui concerne les vitesses, contrairement au cas d\'ecal\'e, on a toujours  $\|\mathbf{u}_p^{n+1}\|<\|\mathbf{u}_p^{n}\|$. On observe \'egalement le ph\'enom\`ene d'\'etirement vu dans le cas d\'ecal\'e, mais dans une moindre mesure.

\subsubsection{Etude pour les sch\'emas de type BBC}

\hspace{1.5em}Le mod\`ele de sch\'ema utilis\'e est le sch\'ema BBC, pour lequel les vitesses sont localis\'ees aux centres des faces du maillage, et dirig\'ees uniquement selon les normales \`a ces faces. Les grandeurs tehrmodynamiques, quant \`a elles, sont toujours localis\'ees aux centres des mailles. Conform\'ement \`a ce qui a \'et\'e vu en 1.2, on place le centre du vortex sur un n\oe ud du maillage.
$\:$\\
\begin{figure}[h]
\centering
   \includegraphics[width=7cm]{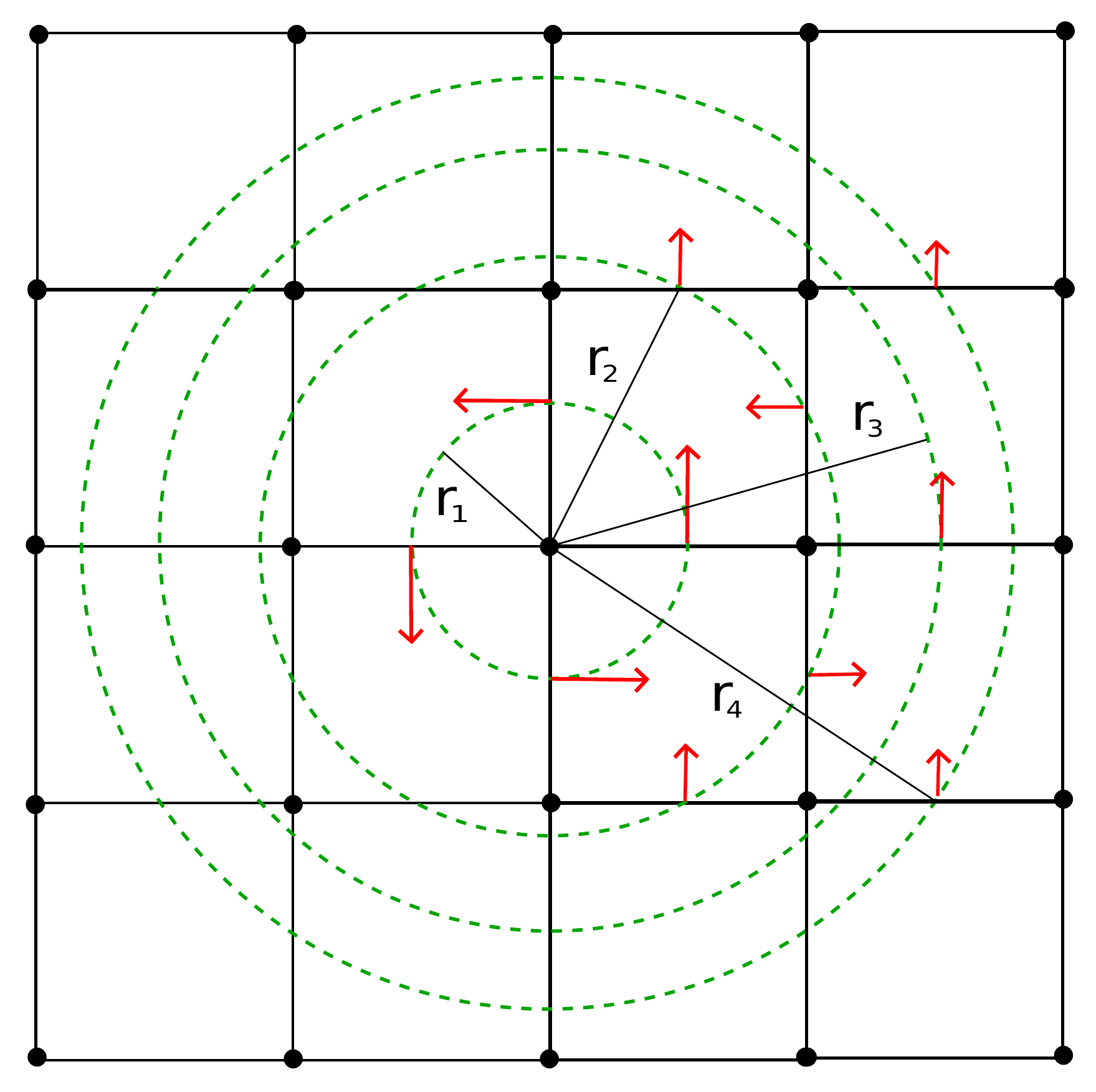}
   \caption{\label{fig:schemabbc} Sch\'ema du champ de vitesses sur le maillage pour les sch\'emas BBC.}   
\end{figure}

\paragraph{D\'etermination de $\mathbf{u}$ sur le maillage}
$\:$

Le calcul des vitesses aux faces n\'ecessite de conna\^itre les rayons associ\'es.
\begin{displaymath}
r_1 = \frac{1}{2}\Delta x,  \:\:\: r_2 = \frac{\sqrt{5}}{2}\Delta x,  \:\:\:  r_3 = \frac{3}{2}\Delta x, \:\:\: r_4 = \frac{\sqrt{13}}{2}\Delta x
\end{displaymath}
Puis, on a pour tout $i \in [1,4]$, $u_i=u_\theta(r_i) \mathbf{e_{\theta_i}}.\mathbf{e_x}$, pour $\theta_i \in \left[0,\frac{\pi}{4}\right]$, d'o\`u
\begin{displaymath}
u_1 = 2\alpha_0\Delta x,  \:\:\: u_2 = \frac{2}{5}\alpha_0\Delta x,  \:\:\:  u_3 = \frac{2}{3}\alpha_0\Delta x, \:\:\: u_4 = \frac{6}{13}\alpha_0\Delta x
\end{displaymath}

\paragraph{Phase lagrangienne}
$\:$

 Pour la phase lagrangienne, la discr\'etisation en temps est effectu\'ee en trois \'etapes, de $t^n$ \`a $t^{n+1/4}$, de $t^{n+1/4}$ \`a $t^{n+1/2}$, et de $t^{n+1/2}$ \`a $t^{lag}$. Pour une maille repr\'esent\'ee par $(i,j)$, ses faces verticales sont d\'esign\'ees par $(i\pm 1/2,j)$, et ses faces horizontales par $(i,j\pm 1/2)$.
\\Premi\`ere \'etape : On calcule les vitesses aux faces \`a $t^{n+1/4}$ \`a partir des forces de pression et des termes de pseudo-viscosit\'e. Dans notre cas, le maillage carr\'e, le champ de pression uniforme et l'absence de termes de pseudo-viscosit\'e (conservation du volume de chaque maille) donnent, pour la face de gauche :
\begin{displaymath}
ux_{i-1/2,j}^{n+1/4}=ux_{i-1/2,j}^{n}
\end{displaymath}
 Il en est de m\^eme pour toutes les faces, ainsi le champ de vitesse \`a $t^{n+1/4}$ est le m\^eme qu'\`a $t^n$.
 \\Deuxi\`eme \'etape : On calcule la variation de volume 
\begin{displaymath}
dV_{i,j}^{n+1/2}=\frac{\Delta x}{\Delta t/2}\left(ux_{i+1/2,j}^{n+1/4} - ux_{i-1/2,j}^{n+1/4} + uy_{i,j+1/2}^{n+1/4} - uy_{i,j-1/2}^{n+1/4}\right)=0
\end{displaymath}
Ce qui donne finalement pour la face de gauche (il en va l\`a encore de m\^eme pour les autres) :
\begin{displaymath}
ux_{i-1/2,j}^{n+1/2}=ux_{i-1/2,j}^{n}
\end{displaymath}
Le champ de vitesse, et les autres grandeurs, n'\'evoluent pas entre $t^n$ et $t^{n+1/2}$.
\\Troisi\`eme \'etape : On calcule la nouvelle variation de volume
\begin{displaymath}
dV_{i,j}^{lag}=\frac{\Delta x}{\Delta t}\left(ux_{i+1/2,j}^{n+1/2} - ux_{i-1/2,j}^{n+1/2} + uy_{i,j+1/2}^{n+1/2} - uy_{i,j-1/2}^{n+1/2}\right)=0
\end{displaymath}
Ce qui assure finalement la conservation des grandeurs thermodynamiques entre $t^n$ et $t^{lag}$. Pour les vitesses, on utilise un sch\'ema centr\'e en temps :
\begin{displaymath}
ux_{i+1/2,j}^{lag} = 2 \: ux_{i+1/2,j}^{n+1/2} - ux_{i+1/2,j}^{n} = ux_{i+1/2,j}^{n}
\end{displaymath}
Le champ de vitesse \`a $t^{lag}$ est donc lui aussi le m\^eme qu'\`a $t^n$.

\paragraph{Phase de projection}
$\:$

La phase de projection est effectu\'ee par splitting d'op\'erateur par directions altern\'ees. On commence par exemple par la projection en $X$. On calcule les flux de masse traversant chaque maille des maillages primal (pour les $ux$) et dual (pour les $uy$). On en d\'eduit les flux de quantit\'e de mouvement, en utilisant un d\'ecentrement upwind (\`a l'ordre 1) pour les vitesses. Et la conservation de la quantit\'e de mouvement donne les nouvelles vitesses aux faces $X$ et $Y$. Ces nouvelles vitesses induisent une nouvelle d\'eformation des maillages, \`a partir de laquelle on effectue la projection en $Y$.

On cherche \`a calculer les vitesses aux faces situ\'ees respectivement en-dessous et \`a droite du n\oe ud o\`u est plac\'e le centre du vortex. On se place pour cela dans la maille en bas \`a droite de ce sommet, que l'on note $(i,j)$. On s'int\'eresse donc aux faces $(i+1/2,j)$ et $(i,j+1-2)$.
\\Les calculs donnent, \`a l'ordre 1 en $\Delta t$ :
\begin{displaymath}
\boxed{
ux_{i-1/2,j}^{n+1} = uy_{i,j+1/2}^{n+1} = 2 \alpha_0 \Delta x \left( 1-\frac{48}{25} \alpha_0 \Delta t \right)
}
\end{displaymath}

\paragraph{Calcul des rotationnels discrets}
$\:$

Comme pour les sch\'emas centr\'es, on calcule la circulation sur le contour de la maille duale. Ainsi, comme on a $\mathbf{u}_{c,c+1}^n\cdot\mathbf{T}_{c,c+1} = 2 \alpha_0 \Delta x$,
\begin{displaymath}
\mbox{rot} (\mathbf{u})_{p,BBC}^n = 8 \alpha_0 = \frac{8 \sigma_0}{2 \pi {(\Delta x)}^2}
\end{displaymath}
On remarque qu'il y a une diff\'erence d'un facteur 2 entre $\mbox{rot} (\mathbf{u})_{c,D}^n = \mbox{rot} (\mathbf{u})_{p,C}^n$ et $\mbox{rot} (\mathbf{u})_{p,BBC}^n$. En effet, la vitesse $u_{\theta}$ \'etant une fonction d\'ecroissante de $r$, les premi\'eres valeurs dont dispose le sch\'ema BBC pour traduire la vorticit\'e (celles aux milieux des faces de la maille duale centrale, donc \`a $\Delta x/2$ du centre) sont plus \'elev\'ees que celles dont disposent les sch\'emas d\'ecal\'e et centr\'e (respectivement celles aux noeuds de la maille primale centrale et de la maille duale centrale, donc \`a $\Delta x/\sqrt{2}$ du centre).
\\Par ailleurs, compte tenu des valeurs des $\mathbf{u}_{c,c+1}^{n+1}$, on obtient :
\begin{displaymath}
\boxed{\frac{\mbox{rot} (\mathbf{u})_{p,BBC}^{n+1}}{\mbox{rot} (\mathbf{u})_{p,BBC}^{n}} = 1 - \frac{48}{25}\alpha_0 \Delta t} 
\end{displaymath}
Le sch\'ema BBC diffuse, comme le sch\'ema centr\'e, la vorticit\'e des vortex ponctuels.

\subsection{Cas d'un vortex "id\'eal"}

\hspace{1.5em}On entend par vortex "id\'eal" un vortex qui induit un mouvement de la mati\`ere dans son ensemble, c'est-\`a-dire qui conserve la vitesse angulaire (rotation solide). Si la vitesse angulaire est la m\^eme en tout point, on a donc un champ de vorticit\'e uniforme :
\begin{displaymath}
\mbox{rot} (\mathbf{u}) = \omega_0
\end{displaymath}

\subsubsection{Calcul de la solution analytique}

\hspace{1.5em}On doit cette fois r\'esoudre l'\'equation suivante pour trouver $\psi_z$ :
\begin{equation*}
- \Delta \psi_z = \omega_0
\end{equation*}
Ce qui donne directement pour la vitesse :
\begin{equation*}
u_{\theta}(r) = - \frac{\partial \psi_z}{\partial r} = \frac{\omega_0}{2}r
\end{equation*}
On cherche alors la densit\'e $\rho$ (et la pression $P$ associ\'ee) telle que $(\mathbf{u},\rho,P)$ soit solution du syst\`eme d'Euler. On rappelle la deuxi\`eme \'equation de \eqref{euler} :
\begin{displaymath}
-\frac{{u_{\theta}}^2}{r} = -\frac{1}{\rho} \frac{\partial P}{\partial r}
\end{displaymath}
\begin{displaymath}
\mbox{qui donne} \:\:\:\: \rho (r) = \rho_0 \exp \big( M_0 r^2 \big), \:\:\:  \mbox{avec}  \:\:\: M_0 = \frac{\omega_0^2}{8(\gamma-1)e_0}
\end{displaymath}
On voit que $\rho(r) \rightarrow \rho_0$ lorsque $M_0 r^2 \rightarrow 0$.
\\On peut supposer que $\rho \simeq \rho_0$ si $M_0 r^2 \ll 1$, pour les $r$ consid\'er\'es ici. Pour les sch\'emas que nous utilisons, nous ne d\'epassons pas $r_{max}=2\sqrt{2}\:\Delta x$. La condition pr\'ec\'edente revient alors \`a :
\begin{equation*}
\boxed{
\Delta x \ll \frac{c_0}{\omega_0}
}
\end{equation*}
Ainsi, pour un maillage suffisamment fin, ou pour un $\omega_0$ tel que $\omega_0 \Delta x$ petit devant la vitesse du son, le triplet $\big(\frac{\omega_0}{2}r\:\mathbf{e}_{\theta},\rho_0,(\gamma-1)\rho_0e_0\big)$ repr\'esente sur notre maillage la solution analytique des \'equations d'Euler incompressible associ\'ee \`a un champ de vorticit\'e uniforme.
$\:$ \\

Dans la suite, pour simplifier les calculs et adopter une notation similaire au cas du vortex ponctuel, on va plut\^ot travailler avec $\boldsymbol{\alpha_0=\omega_0/2}$ de sorte que $u_{\theta}(r)=r\alpha_0$, qui repr\'esente donc la vitesse angulaire r\'eelle en tout point.

\subsubsection{Etude pour les sch\'emas d\'ecal\'es}

\paragraph{D\'etermination de $\mathbf{u}$ sur le maillage}
$\:$

Comme pour le cas du vortex ponctuel, on calcule les vitesses aux n\oe uds \`a $t^n$ (projet\'ees \`a partir de la solution analytique). On a \'evidemment les m\^emes rayons puisque le maillage n'a pas chang\'e.
\begin{displaymath}
r_1 = \frac{\Delta x}{\sqrt{2}},  \:\:\:\: r_2 = \sqrt{\frac{5}{2}}\Delta x,  \:\:\:\:  r_3 = 3\frac{\Delta x}{\sqrt{2}}
\end{displaymath}
\begin{displaymath}
\mbox{et donc :}\:\:\:\: u_1 = \frac{1}{\sqrt{2}}\alpha_0\Delta x,  \:\:\:\: u_2 = \sqrt{\frac{5}{2}}\alpha_0\Delta x,  \:\:\:\:  u_3 = \frac{3}{\sqrt{2}}\alpha_0\Delta x
\end{displaymath}
L\`a encore il ne reste qu'\`a multiplier par le $\mathbf{e}_{\theta}$ correspondant pour obtenir les vecteurs vitesse.

\paragraph{Phase lagrangienne}
$\:$

\begin{figure}[h]
\centering
   \includegraphics[width=6cm]{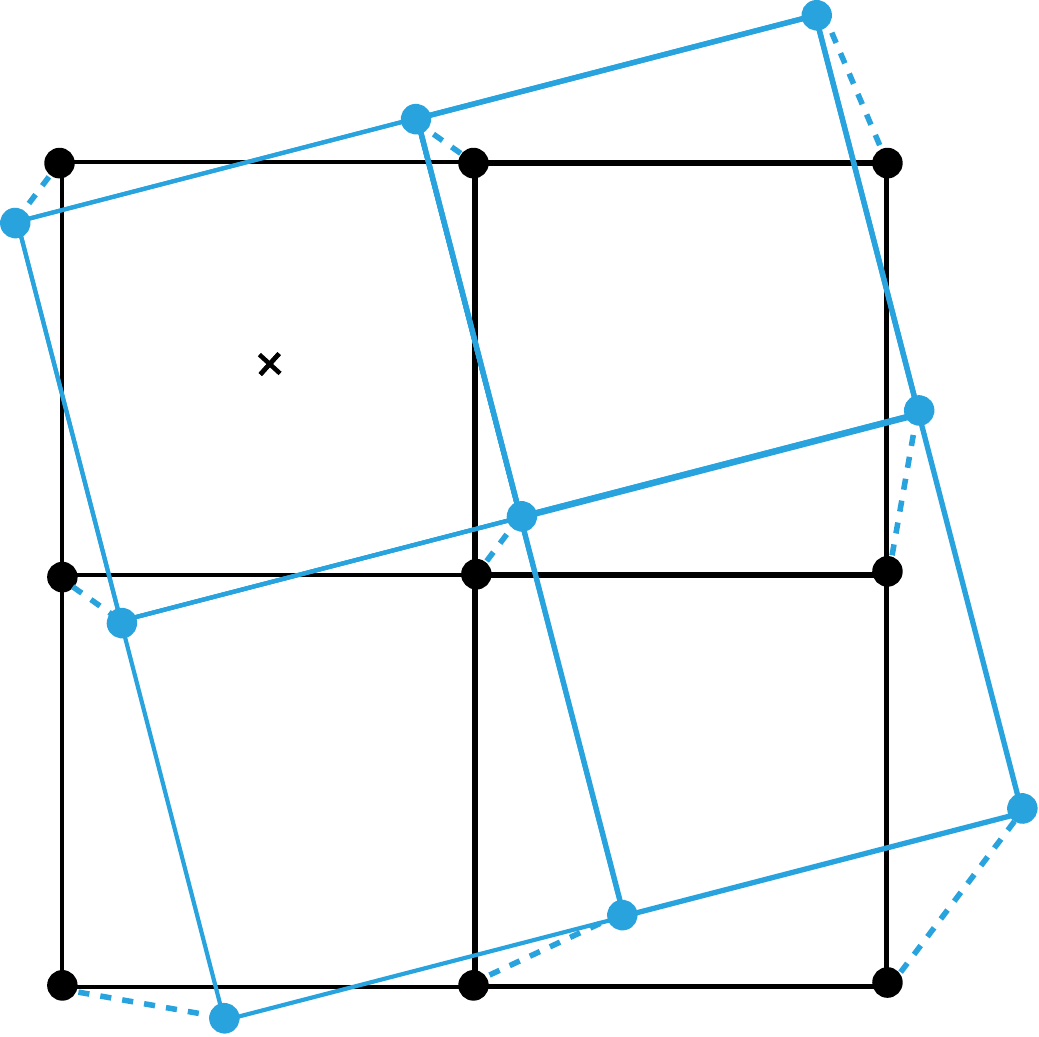}
   \caption{Sch\'ema de la d\'eformation lagrangienne dans le cas d\'ecal\'e pour un vortex id\'eal.}   
   \label{fig:deformation_ideal_decale}
\end{figure}

La conservation du volume de chaque maille et le champ de pression uniforme impliquent que, pour tout n\oe ud $p$ :
\begin{displaymath}
\boxed{
\mathbf{u}_p^{lag} = \mathbf{u}_p^{n}
}
\end{displaymath}
On remarque \'egalement que la forme carr\'ee des mailles est conserv\'ee, voir Figure~\ref{fig:deformation_ideal_decale}. Toute la mati\`ere \'etant d\'eplac\'ee \`a la m\^eme vitesse par le vortex, il est normal que l'on observe une simple rotation de l'ensemble du maillage.

\paragraph{Phase de projection}
$\:$

Comme en 2.2.3, on d\'etermine les flux aux faces de la maille nodale, et on applique la conservation de la quantit\'e de mouvement \eqref{dual_remap}, en prenant des valeurs pour les $u_p^\alpha$ correspondant \`a chaque ordre.
\\ \`A l'ordre 1, on a les vitesses suivantes :
\begin{displaymath}
\mathbf{u}_p^O = \frac{\alpha_0 \Delta x}{2} \left( \!\! \begin{array}{clcr}1\\-1 \end{array} \!\!\right), \:\:\: \mathbf{u}_p^S = \frac{\alpha_0 \Delta x}{2} \left( \!\! \begin{array}{clcr}3\\1 \end{array} \!\!\right), \:\:\: \mathbf{u}_p^E = \frac{\alpha_0 \Delta x}{2} \left( \!\! \begin{array}{clcr}1\\3 \end{array} \!\!\right), \:\:\: \mathbf{u}_p^N = \frac{\alpha_0 \Delta x}{2} \left( \!\! \begin{array}{clcr}-1\\1 \end{array} \!\!\right)
\end{displaymath}
Ce qui donne donc, pour le n\oe ud $p$:
\begin{displaymath}
\boxed{
\mathbf{u}_{p,1}^{n+1}= \frac{1}{2}\alpha_0 \Delta x \left( \! \begin{array}{clcr}1+ 2\alpha_0 \Delta t \\ \: \\1- 2\alpha_0 \Delta t \end{array} \!\right)
}
\end{displaymath}
Puis, pour l'ordre 2 :
\begin{displaymath}
\mathbf{u}_p^O = \alpha_0 \Delta x \left( \!\! \begin{array}{clcr}\frac{1}{2}\\0 \end{array} \!\!\right), \:\:\: \mathbf{u}_p^S = \alpha_0 \Delta x \left( \!\! \begin{array}{clcr}1\\\frac{1}{2} \end{array} \!\!\right), \:\:\: \mathbf{u}_p^E = \alpha_0 \Delta x \left( \!\! \begin{array}{clcr}\frac{1}{2}\\1 \end{array} \!\!\right), \:\:\: \mathbf{u}_p^N = \alpha_0 \Delta x \left( \!\! \begin{array}{clcr}0\\\frac{1}{2} \end{array} \!\!\right)
\end{displaymath}
On obtient donc \`a l'ordre 2 :
\begin{displaymath}
\boxed{
\mathbf{u}_{p,2}^{n+1}= \frac{1}{2}\alpha_0 \Delta x \left( \! \begin{array}{clcr}1+ \alpha_0 \Delta t\\ \: \\1- \alpha_0 \Delta t \end{array} \!\right)
}
\end{displaymath}

\paragraph{Calcul des rotationnels discrets}
$\:$

On a dans ce cas-l\`a $\frac{1}{2}\big(\mathbf{u}_p^n+\mathbf{u}_{p+1}^n\big)\cdot\mathbf{T}_{p,p+1} = \frac{1}{2} \alpha_0 \Delta x$, ce qui donne :

\begin{displaymath}
\mbox{rot} (\mathbf{u})_{c,D}^n = 2 \: \alpha_0 = \omega_0
\end{displaymath}
Ce qui n'est pas \'etonnant compte tenu du choix d'un champ de vorticit\'e uniforme et \'egal \`a $\omega_0$.
\\ Puis, \`a $t^{n+1}$, on a quel que soit l'ordre :
\begin{displaymath}
\frac{1}{2}\big(\mathbf{u}_p^{n+1}+\mathbf{u}_{p+1}^{n+1}\big)\cdot\mathbf{T}_{p,p+1} = \frac{1}{2} \alpha_0 \Delta x
\end{displaymath} 
\begin{displaymath}
\mbox{et donc,} \:\:\:\: \boxed{\frac{\mbox{rot} (\mathbf{u})_{c,D}^{n+1}}{\mbox{rot} (\mathbf{u})_{c,D}^{n}} = 1} 
\end{displaymath}
$\:$ \\

Les sch\'emas d\'ecal\'es ne diffusent donc pas les vortex id\'eaux non plus. Et comme c'\'etait le cas pour les vortex ponctuels, on observe une amplification de la norme des vitesses nodales.

\subsubsection{Etude pour les sch\'emas centr\'es}

\paragraph{D\'etermination de $\mathbf{u}$ sur le maillage}
$\:$

Les rayons \'etant les m\^emes qu'en 2.3, on a, confrom\'ement aux notations pr\'ec\'edentes :
\begin{displaymath}
u_1 = \frac{1}{\sqrt{2}} \: \alpha_0\Delta x,  \:\:\:\: u_2 = \sqrt{\frac{5}{2}}\alpha_0\Delta x,  \:\:\:\:  u_3 = \frac{3\sqrt{2}}{2}\alpha_0\Delta x,  
\end{displaymath}

\paragraph{Phase lagrangienne}
$\:$

\begin{figure}[h]
\centering
   \includegraphics[width=6cm]{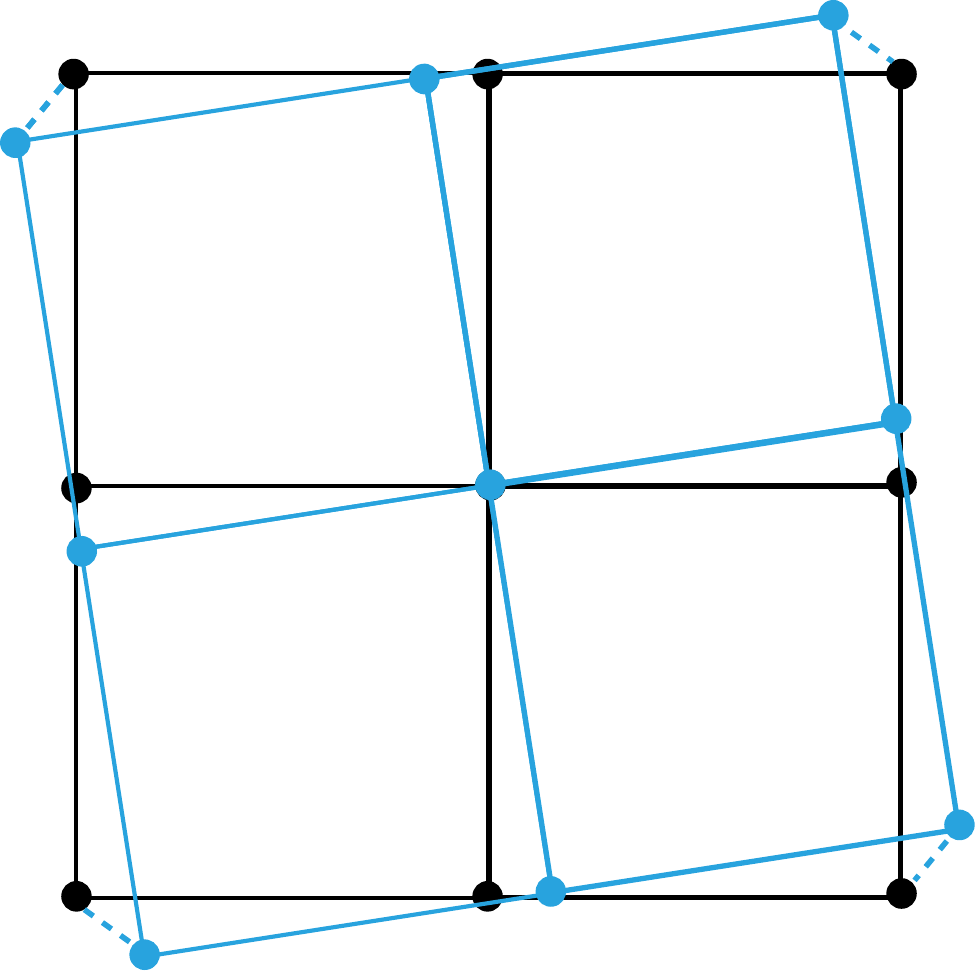}
   \caption{Sch\'ema de la d\'eformation lagrangienne dans le cas centr\'e pour un vortex id\'eal.}
   \label{fig:deformation_ideal_centre}   
\end{figure}

On acc\`ede aux vitesses nodales par l'interm\'ediare de \eqref{nodal_solver}, et on obtient la d\'eformation lagrangienne suivante (voir Figure~\ref{fig:deformation_ideal_centre}). On remarque qu'il n'y a aucune diff\'erence avec le cas d\'ecal\'e, l'ensemble du maillage est en rotation solide.

Ensuite, il reste \`a calculer les vitesses centr\'ees \`a $t^{lag}$ avec \eqref{lagrangian_velocity_centred}. Les calculs m\`enent au r\'esultat suivant pour la vitesse au centre de la maille $c$ en bas \`a droite :
\begin{displaymath}
\boxed{
\mathbf{u}_c^{lag} = \mathbf{u}_c^{n} = \frac{1}{2}\alpha_0 \Delta x \left( \!\! \begin{array}{clcr} 1\\1\end{array} \!\!\right)
}
\end{displaymath}

\paragraph{Phase de projection}
$\:$

Enfin, pour la projection, on proc\`ede exactement comme en 2.3. Apr\`es calcul des flux de masse aux faces de la maille  $c$, et application de la conservation de qantit\'e de mouvement, il vient :
\begin{displaymath}
\boxed{
\mathbf{u}_c^{n+1} = \frac{1}{2}\alpha_0 \Delta x \left( \!\! \begin{array}{clcr} 1+\alpha_0 \Delta t\\ \: \\ 1-\alpha_0 \Delta t \end{array} \!\!\right)
}
\end{displaymath}

\paragraph{Calcul des rotationnels discrets}
$\:$

Comme cette fois aussi $\frac{1}{2}\big(\mathbf{u}_c^n+\mathbf{u}_{c+1}^n\big)\cdot\mathbf{T}_{c,c+1} = \frac{1}{2}\alpha_0 \Delta x$,
\begin{displaymath}
\mbox{rot} (\mathbf{u})_{p,C}^n = 2 \: \alpha_0 = \omega_0
\end{displaymath}
De plus, \`a $t^{n+1}$, on obtient :
\begin{displaymath}
\frac{1}{2}\big(\mathbf{u}_c^{n+1}+\mathbf{u}_{c+1}^{n+1}\big)\cdot\mathbf{T}_{c,c+1} = \frac{1}{2}\alpha_0 \Delta x
\end{displaymath} 
\begin{displaymath}
\mbox{d'o\`u} \:\:\:\: \boxed{\frac{\mbox{rot} (\mathbf{u})_{p,C}^{n+1}}{\mbox{rot} (\mathbf{u})_{p,C}^{n}} = 1}
\end{displaymath}
$\:$ \\

Contrairement au cas des vortex ponctuels, les sch\'emas centr\'es ne vont pas diffuser la vorticit\'e induite par des vortex id\'eaux. En effet, pour ce type de vortex, les sch\'emas d\'ecal\'es (avec projection d'ordre 2) et centr\'es donnent exactement les m\^emes r\'esultats : m\^emes d\'eformations lagrangiennes et m\^emes vitesses \`a $r=\Delta x/\sqrt{2}$ du centre du vortex, \`a $t^{n+1}$. Il est donc logique qu'on observe une conservation de la vorticit\'e pour les deux. On note \'egalement le m\^eme \'etirement du vortex et la m\^eme amplification de la norme des vitesses.

\subsubsection{Etude pour les sch\'emas de type BBC}

\paragraph{D\'etermination de $\mathbf{u}$ sur le maillage}
$\:$

Le maillage \'etant le m\^eme qu'en 2.4, et le centre du vortex ayant le m\^eme position, on garde les valeurs des rayons, avec cette fois les vitesses
\begin{displaymath}
u_1 = \frac{1}{2}\alpha_0\Delta x,  \:\:\: u_2 = \frac{1}{2}\alpha_0\Delta x,  \:\:\:  u_3 = \frac{3}{2}\alpha_0\Delta x, \:\:\: u_4 = \frac{3}{2}\alpha_0\Delta x
\end{displaymath}

\paragraph{Phase lagrangienne}
$\:$

La phase lagrangienne est rigoureusement identique au cas du vortex ponctuel (voir 2.4). On retrouve \`a $t^{lag}$  le m\^eme champ de vitesse qu'\`a $t^n$.

\paragraph{Phase de projection}
$\:$

Conform\'ement aux notations adopt\'ees en 2.4, on a cette fois apr\`es projection par splitting des vitesses aux faces de la maille en bas \`a droite, \`a l'ordre 1 en $\Delta t$ :

\begin{displaymath}
\boxed{
ux_{i-1/2,j}^{n+1} = uy_{i,j+1/2}^{n+1} = \frac{3}{2} \alpha_0 \Delta x
}
\end{displaymath}
Comme pour le vortex ponctuel, on voit que la sym\'etrie radiale est conserv\'ee \`a l'ordre 1.

\paragraph{Calcul des rotationnels discrets}
$\:$

\`A l'instant $t^n$, on a $\mathbf{u}_{c,c+1}^n\cdot\mathbf{T}_{c,c+1} = \frac{1}{2} \alpha_0 \Delta x$,
\begin{displaymath}
\mbox{rot} (\mathbf{u})_{p,BBC}^n = 2 \alpha_0 = \omega_0
\end{displaymath}
Pour tous les types de sch\'emas, dans le cas d'une rotation solide, on a donc bien un rotationnel discret \'egal \`a $\omega_0$ \`a $t^n$.
\\De plus, compte tenu des valeurs des $\mathbf{u}_{c,c+1}^{n+1}$, on obtient \`a l'ordre 1 :
\begin{displaymath}
\boxed{
\frac{\mbox{rot} (\mathbf{u})_{p,BBC}^{n+1}}{\mbox{rot} (\mathbf{u})_{p,BBC}^{n}} = 1
} 
\end{displaymath}
Le sch\'ema BBC ne diffuse pas la vorticit\'e des vortex ponctuels (\`a l'ordre 1). En effet, si on pousse le calcul des $\mathbf{u}_{c,c+1}^{n+1}$ jusqu'\`a l'ordre 2 en $\Delta t$, on voit appara\^itre un terme qui brise la sym\'etrie du syst\`eme (qui d'ailleurs existe aussi pour le cas du vortex ponctuel) d\^u \`a la m\'ethode de projection par splitting directionnel. Il en r\'esulte l'apparition d'un terme diffusif \'egal \`a $\frac{1}{2}{\big(\alpha_0 \Delta t\big)}^2$ dans le rapport des rotationnels. Ce qui n'est pas le cas pour les sch\'emas centr\'es et d\'ecal\'es, qui conservent exactement la vorticit\'e des vortex id\'eaux.

\clearpage

\subsection*{Conclusion}

\hspace{1.5em}Effectuons pour terminer un bilan des r\'esultats obtenus. On rassemble dans un tableau les diff\'erentes valeurs du rapport $\frac{\mbox{rot} (\mathbf{u})^{n+1}}{\mbox{rot} (\mathbf{u})^{n}}$ pour chacun des cas trait\'es.

\begin{displaymath}
  \begin{array}{|c|c|c|}
\hline
    \: & \mbox{Vortex ponctuel} & \: \: \: \mbox{Vortex id\'eal} \: \: \: \\
	\hline
    \mbox{Sch\'ema d\'ecal\'e} & 1 & 1 \\
\hline
    \mbox{Sch\'ema centr\'e} & 1-\frac{4}{3}\frac{\Delta t}{\Delta x}c_0-\frac{14}{25}\alpha_0\Delta t & 1 \\
\hline
    \mbox{Sch\'ema BBC} & 1-\frac{48}{25}\alpha_0 \Delta t & 1 \\
	\hline
  \end{array}
\end{displaymath}
$\:$ \\

Il en ressort que seul le sch\'ema d\'ecal\'e ne diffuse ni la vorticit\'e des vortex ponctuels, ni celle des vortex id\'eaux. 
\\Une premi\`ere explication peut venir du fait que le maillage primal n'est conforme dans $H(\mbox{rot})$ que pour ces sch\'emas-l\`a. En effet, pour les autres sch\'emas, bien que le choix du point de calcul du rotationnel (i.e le choix du centre du vortex) ait \'et\'e fait en prenant ce fait en consid\'eration, il se trouve que le champ de vitesse \'evolue sur le maillage primal, qui lui n'est pas conforme dans $H(\mbox{rot})$. Ce qui peut \^etre source de diffusion num\'erique.
\\Une deuxi\`eme explication est li\'ee aux stencils de chaque type de sch\'emas. Dans le cas d\'ecal\'e, pour permettre le calcul des vitesses \`a $t^{n+1}$ n\'ecessaires \`a la d\'etermination de $\mbox{rot}(\mathbf{u})$, il faut utiliser les valeurs des vitesses nodales de la maille centrale et des 4 mailles directement voisines (sch\'ema \`a 5 points), ce qui fait au total 12 vitesses. Alors que pour les sch\'emas BBC et centr\'es, il faut utiliser les vitesses aux noeuds de la maille duale centrale et des 8 mailles duales voisines (sch\'ema \`a 9 points), d'o\`u au total 16 vitesses. On fait les calculs avec des vitesses plus \'eloign\'ees du centre du vortex, ce qui explique, dans le cas ponctuel, comme $u_\theta$ est une fonction d\'ecroissante de $r$, une diffusion plus importante pour ces sch\'emas-l\`a.

Enfin, on observe qu'on a un crit\`ere portant sur la CFL \`a partir duquel le sch\'ema centr\'e devient moins diffusif que le sch\'ema BBC. Pour commencer, le sch\'ema BBC est stable si :
\begin{displaymath}
\alpha_0\Delta t = \frac{\Delta t}{\Delta x} \frac{\sigma_0}{2 \pi \Delta x} \leq  \frac{25}{48}
\end{displaymath}
De plus, si BBC diffuse plus la vorticit\'e que les sch\'emas centr\'es alors cela signifie que 
\begin{displaymath}
\frac{4}{3} \frac{\Delta t}{\Delta x} c_0+\frac{14}{25} \alpha_0 \Delta t \: \leq \:\frac{48}{25}\alpha_0\Delta t ,
\end{displaymath}
ce qui donne finalement :
\begin{displaymath}
\boxed{
\mbox{CFL} = \frac{\Delta t}{\Delta x}c_0 \leq \frac{17}{32} = 0,53125
}
\end{displaymath}
Ce qui le cas puisque $\mbox{CFL}_{BBC} = 0,34$ classiquement.

  \clearpage
  \newpage
  \strut
  \newpage	
  \bibliography{bibliographie}
  \addcontentsline{toc}{section}{References}
  
\end{document}